\documentclass{amsart}
\newtheorem{Theorem}{Theorem}[section]

\newtheorem{Lemma}[Theorem]{Lemma}

\newtheorem{Corollary}[Theorem]{Corollary}

\newtheorem{Remark}[Theorem]{Remark}

\newtheorem{Definition}[Theorem]{Definition}

\begin{document}

\title{Resolution of singularities for
 3-folds in positive characteristic}

\author{Steven Dale Cutkosky}
\thanks{Research    partially supported by NSF}

\maketitle

\begin{abstract}  A concise, complete proof of resolution of singularities
of 3-folds in positive characteristic $>5$ is given. Abhyankar first proved this theorem in 1966. 
\end{abstract}

\section{Introduction}
\markboth{Steven Dale Cutkosky}{Resolution of singularities for
3-folds}

The  purpose of this paper is to give a complete proof of
resolution of singularities of 3-folds in positive characteristic
which is concise.  The paper is essentially self contained, and includes proofs
of resolution of singularities of curves and surfaces in positive characteristic.

Resolution of singularities of algebraic varieties over fields of
characteristic zero was solved in all dimensions in a very
complete form by Hironaka \cite{H1} in 1964.

Resolution of singularities of varieties over fields of positive characteristic is significantly harder.
  One explanation for this is the lack of ``hypersurfaces of maximal contact'' in positive characteristic,
   which is the main technique used in  characteristic zero proofs to reduce to a lower dimension (``Maximal contact and approximate manifolds'' of this introduction and Sections 6.2 and 7.4 \cite{C2}).

There are several proofs of resolution of singularities of surfaces in positive characteristic.
 The first one is by Abhyankar \cite{Ab1} which appeared in 1955. There are later proofs by Hironaka,
  outlined in his  notes \cite{H2}, and  by Lipman \cite{L2} for excellent surfaces.

Abhyankar proved resolution of singularities of positive characterisitic 3-folds in 1966. The proof  appears in his book \cite{Ab7} and
the papers \cite{Ab5}, \cite{Ab6}, \cite{Ab8} and \cite{Ab9}. The entire proof is extremely long and difficult. It encompasses
508 pages. 

Within the last few years there has been a resurgence of interest in resolution in
positive characteristic.  Some  papers making progress on
resolution in positive characteristic are Cossart \cite{Co1}, \cite{Co2}, \cite{Co3},
\cite{Co4}, Giraud \cite{G}, Hauser 
\cite{Hau1}, de Jong \cite{Jo}, Hironaka \cite{H9}, Kawanoue \cite{Ka},
Kawanoue and Matsuki \cite{KM}, Kuhlmann \cite{Ku}, Moh \cite{M}, 
Piltant \cite{P},
Spivakovsky \cite{S} and Teissier \cite{T}. Several simplifications of Hironaka's proof of
characteristic zero resolution have appeared, making the proof quite accessible now, including Abramovich and de Jong \cite{AJo},
Bierstone and Milman \cite{BrM}, Bogomolov and Pantev  \cite{BP},
Bravo, Encinas and Villamayor  \cite{BEV}, Encinas and Hauser \cite{EH}, Encinas and Villamayor \cite{EV}, Hauser  \cite{Hau2}, Koll\'ar \cite{Ko}, Villamayor \cite{V}, Wlodarczyk  \cite{W}.

Most recently, Cossart and Piltant \cite{CP1}, \cite{CP2}  have proved resolution of singularities
for a 3-fold over an arbitrary field. The proof is  extremely long
and difficult.

We give a self contained and concise proof of the following theorem in this paper.

\begin{Theorem}\label{Theorem82} Suppose that $V$ is a projective variety of dimension 3 over an algebraically closed field ${\bf k}$ of characteristic $\ne 2,3$ or $5$.  Then there exists a
nonsingular projective variety $W$ and a birational morphism $\phi:W\rightarrow V$, which is an isomorphism above the nonsingular locus of $V$. 
\end{Theorem}

The restrictions on the field $\bf k$ in the statement of Theorem \ref{Theorem82} are those of Abhyankar's original statement. 
 
It is our hope that this paper will motivate study of the fundamental open problems in resolution.
Our proof could have been written in 1967, as it only makes use of methods known  at that time.

Essential ingredients in the proof are the following  theorems on embedded resolution of surface singularities.

\begin{Theorem}\label{Theorem1}(Embedded resolution of surface singularities)
Suppose that  $V$ is a nonsingular  3-dimensional variety over an algebraically closed field ${\bf k}$,  $S$ is a reduced surface (a pure 2-dimensional reduced closed subscheme) in $V$
and $E$ is a simple normal crossings divisor on $V$.
 Then there exists a sequence of blow ups 
 $$
\pi:V_n\rightarrow V_{n-1}\rightarrow \cdots \rightarrow V
$$
  such that the strict transform $S_n$ of $S$ on $V_n$ is nonsingular, and the divisor $\pi^{*}(S+E)$ is a simple normal crossings divisor on $S_n$.
Further, each $V_i\rightarrow V_{i-1}$ is the blow up of a point or nonsingular curve in the locus in $V_{i-1}$ where the preimage of $S+E$ is not a simple normal crossings divisor.
\end{Theorem}

An effective divisor is a simple normal crossings divisor if its components 
are nonsingular and intersect transversely (Section \ref{SecPre}).
The morphisms $V_i\rightarrow V_{i-1}$ of Theorem \ref{Theorem1} are required to be permissible (Definition \ref{Def101}). In particular, the centers blown up  are required to be transversal to the preimage of $E$.

\begin{Theorem}\label{Theorem23} (Principalization of ideals)
 Suppose that $V$ is a nonsingular 3 dimensional variety over an algebraically closed field
  $\bf k$ and ${\mathcal I}\subset {\mathcal O}_V$ is a nonzero  ideal sheaf on $V$. Then there exists a sequence of blow ups
$$
V_n\rightarrow V_{n-1}\rightarrow \cdots \rightarrow V
$$
such that ${\mathcal I}{\mathcal O}_{V_n}$ is locally principal. Further,
each $V_i\rightarrow V_{i-1}$ is the blow up of a point
 or nonsingular curve in the locus in $V_{i-1}$ where ${\mathcal I}{\mathcal O}_{V_{i-1}}$ is not invertible.
\end{Theorem}

Theorems \ref{Theorem1} and \ref{Theorem23} were first proven  by Abhyankar in the book
\cite{Ab7}, relying on material proven in  \cite{Ab5}, \cite{Ab6}, \cite{Ab8} and \cite{Ab9}.
 This is the most difficult part of Abhyankar's proof. 

We give self contained proofs of these theorems in Sections \ref{Section4} - \ref{SecERSI}.
The algorithm that we follow for resolution is essentially that of
Levi \cite{Le} to resolve the singularities of a two dimensional hypersurface, embedded in a nonsingular 3-fold. This method  was given a completely rigourous proof in characteristic zero by 
Zariski in \cite{Z4}.   The algorithm will be discussed at greater length in 
 Section \ref{Intro2} on resolution of surfaces. The first part of this
algorithm extends without much trouble to positive characteristic. 
Hironaka constructs in \cite{H2} an invariant   and  outlines a proof 
that in arbitrary characteristic, the
invariant drops under the resolution algorithm of  Levi.
This involves overcoming new major obstacles in positive characteristic.
Our proof requires an extension of Hironaka's proof for a 2 dimensional hypersurface to resolution of ideal sheaves on a nonsingular 3-fold of arbitrary characteristic,
with the further condition that the resolution procedure produces transversality with  the exceptional divisor.  Our main resolution theorem for ideal sheaves
on a nonsingular 3-fold is stated in Theorem \ref{Theorem6}. We use the 
philosophy of Bravo, Encinas and Villamayor in  \cite{EV} and \cite{BEV},
resolving by taking the weak transform of an ideal rather  than the strict
transform. Our main resolution invariant is then the order of an ideal, rather
than the multiplicity. This leads to vast simplifications in the proofs of resolution theorems.

The most fundamental problem of resolution in positive characteristic (open at the time of this writing)
is to construct an embedded resolution of a 3 dimensional hypersurface in a
nonsingular 4 dimensional variety.

To circumvent the difficulty of constructing an embedded resolution of a 3-fold, Abhyankar proves and makes use of the following theorem, which is a generalization of
a method of Albanese. 

\begin{Theorem}\label{Theorem19}  Suppose that $K$ is an algebraic function field of dimension $d$ over an algebraically closed field ${\bf k}$. Then there exists  a normal projective variety $V$ such that $K$ is the function field of $V$, and all points of $V$ have multiplicity $\le d!$.
\end{Theorem}

Abhyankar then shows that 3-fold singularities of multiplicity less than $p$  can be resolved essentially using characteristic zero techniques. This is where the restriction $p>5$ comes in,
as 3!=6. Abhyankar uses a generalization of Jung's method of local resolution,
resolving the branch locus of a finite projection 
 to a nonsingular 3-fold, using Theorem \ref{Theorem1}.

Our proof of Theorem \ref{Theorem82} can be broken down into 5
steps, which we enumerate.
\begin{enumerate}
\item[1.] Prove Theorems \ref{Theorem1} and \ref{Theorem23}
(embedded resolution of surface singularities and principalization
of ideals). \item[2.] Prove that a  projective variety of dimension
$n$ is birationally equivalent to a normal projective variety all
of whose points have multiplicity $\le n!$. \item[3.] Prove that
points of  multiplicity  $<p$ can be resolved locally. \item[4.] Patch local resolutions  to
produce a nonsingular projective variety, birationally equivalent
to $V$.
\item[5.] Modify a resolution of singularities $W\rightarrow V$ to produce
a resolution $W_1\rightarrow V$ which is an isomorphism over the nonsingular locus of $V$.
\end{enumerate}

Steps 1 - 4 appear in Abhyankar's original proof \cite{Ab7}.  They draw on
almost all classical approaches to resolution. Steps 1 and 4  already appear in
the characteristic zero proof of resolution of algebraic 3-folds of Zariski in 
\cite{Z4}. Step 5 appears in Cossart's article \cite{Co4}. 

Step 2 (Theorem \ref{Theorem19})  produces a variety $V_0$, birationally equivalent to our given 3-fold
$V$, such that all points of $V_0$ have multiplicity $\le 3!=6$. Thus all
points of $V_0$ have multiplicity less than the characteristic $p$ of the
ground field ${\bf k}$ (if ${\bf k}$ has positive characteristic).
Step 3 produces for all $p\in V_0$ local resolutions of singularities $W_p\rightarrow U_p$, where $U_p$ is an affine neighborhood of $p$ in $V_0$ and $W_p\rightarrow U_p$ is a birational projective morphism such that $W_p$
is nonsingular. We have thus proven local uniformization of the function
field of $V$. The concept of local uniformization is explained in ``Valuation rings and local uniformization'' in Section \ref{Intro3}. Now in Step 4, we show that we can construct from these local
resolutions a projective morphism $V_1\rightarrow V$ such that $V_1$ is
nonsingular. In Step 5, we show that we can construct from $V_1$ a resolution 
of singularies $V_2\rightarrow V$ such that this morphism is an 
isomorphism away from the singular locus of $V$. This will complete the proof of  Theorem
\ref{Theorem82}.

In Section \ref{Intro2} we outline the proof of Step 1,
(embedded resolution of surfaces) and in Section \ref{Intro3}, we outline the proof of Steps 2-5.  
The actual  proof of Theorem \ref{Theorem82} starts with Section \ref{SecPre}.

There is  another, later, proof of Step 3 of Theorem \ref{Theorem82}, making 
 use of contributions of Hironaka, Giraud and Cossart. We outline the proof in ``Hypersurfaces of maximal contact for singularities of low multiplicity'' in Section \ref{Intro3}.
 In this proof,  the essential point is also that 3-fold singularities of multiplicity less than $p$  can be resolved essentially using characteristic zero techniques, but in this case it is accomplished by showing that 
under this restriction on multiplicity, a hypersurface of maximal contact exists.

The author would like to thank the three reviewers of this paper for their
careful and critical reading, and for their constructive comments.

I dedicate this paper to Zariski, Abhyankar, Hironaka and Lipman.

\section{On  Resolution of Surfaces}\label{Intro2}

In this section, we give an overview of the proof of  the embedded resolution theorem, Theorem  \ref{Theorem6}, and its corollaries, 
Theorem \ref{Theorem1} and Theorem \ref{Theorem23}. These theorems are proven in Sections \ref{Section4} - \ref{SecERSI} of this paper. This will complete the proof of the first step of the proof
of Theorem \ref{Theorem82}.
These theorems are generalizations of the proof of resolution of two dimensional
hypersurfaces outlined in Hironaka's lectures \cite{H2}, which is itself a generalization of the characteristic zero proof of Levi \cite{Le} and Zariski  \cite{Z4} to arbitrary characteristic.   Hironaka's invariant 
 is given an elegant  geometric interpretation by Hauser in \cite{Hau1}, which we will discuss later in this section. All particulars of Hironaka's proof are
worked out in Chapter 7, ``Resolution of Surfaces in Positive Characteristic'',
in the book ``Resolution of Singularities'' \cite{C2}. 
This method of
resolution for a hypersurface on a nonsingular 3-fold extends to resolve an ideal sheaf on a nonsingular 3-fold.  To accomplish this, we 
make a slightly different definition of  the $\tau$ invariant from the classical one of Hironaka.

\vskip .2truein
\noindent {\bf Maximal contact   and approximate manifolds.}  Suppose that $V$ is a nonsingular variety over an algebraically closed field ${\bf k}$.
We use the definition of order, $\nu_q({\mathcal I})$, given in Section  \ref{SecPre}. $\nu_q({\mathcal I})$ is the largest power of the maximal ideal
of ${\mathcal O}_{V,q}$ containing ${\mathcal I}_q$.
Let 
$$
r=r({\mathcal I})=\mbox{max}\{\nu_q({\mathcal I})\mid q\in V\},
$$
and let 
$$
\mbox{Sing}_r({\mathcal I})=\{q\in V\mid \nu_q({\mathcal I})=r\}.
$$

We summarize a few basic properties of the behavior of the order $\nu$ under blow ups. 
This topic is discussed in depth in \cite{C2}, especially in Chapter 4. Some definitions and basic facts about resolution are given in Sections \ref{SecPre} and \ref{Section4} of this paper.

$\mbox{Sing}_r({\mathcal I})$
is a closed subset of $V$.  Suppose that $W\subset \mbox{Sing}_r({\mathcal I})$ is a  nonsingular subvariety.
Let $\pi_1:V_1\rightarrow V$ be the blow up of $W$. The weak transform ${\mathcal I}_1$ of ${\mathcal I}$,
 is defined to be ${\mathcal I}{\mathcal I}_E^{-r}$, where $E=\pi_1^{-1}(W)$ is the exceptional divisor 
of $\pi_1$, and ${\mathcal I}_E$ is the ideal sheaf of $E$ in $V_1$. We have that $r_1=r({\mathcal I}_1)\le r$.

Suppose that $q\in \mbox{Sing}_r({\mathcal I})$.
Let $R={\mathcal O}_{V,q}$. 
A nonsingular subvariety  $Y\subset \mbox{Spec}(R)$ 
 is called a manifold of maximal 
contact   for ${\mathcal I}$ at $q$ if 
$Y$ contains the germ of $\mbox{Sing}_r({\mathcal I})$ at $q$, and if 
$W\subset \mbox{Sing}_r({\mathcal I})$ is a nonsingular subvariety, 
and $\pi_1:V_1\rightarrow V$ is the blow up of $W$, then
$\pi_1^{-1}(q)\cap \mbox{Sing}_r({\mathcal I}_1)$ is contained in the 
 strict transform $Y_1$ of $Y$ on $V_1$. Further, if $q_1\in \pi_1^{-1}(q)\cap \mbox{Sing}_r({\mathcal I}_1)$, then $Y_1$ is a manifold of maximal contact
for ${\mathcal I}_1$ at $q_1$.

The fundamental theorem of resolution of singularities in characteristic
zero is that manifolds of maximal contact always exist (for any dimension of $V$). This allows all
theorems on resolution to  by reduced by induction to 
appropriate resolution statements in lower dimension.

We  associate a positive  integer $\tau=\tau_{\mathcal I}(q)$
to $q\in \mbox{Sing}_r({\mathcal I})$ 
such that if a manifold of maximal contact $Y$ exists at $q$ for ${\mathcal I}$, we have 
$\mbox{dim }Y=\tau$.  The number $\tau$ can be thought of
as the minimal number of variables in the tangent cone of ${\mathcal I}$.

Unfortunately, manifolds of maximal contact do not generally exist 
in positive characteristic. This can be  considered as the major 
obstacle to resolution.
We can however approximate a manifold of maximal contact to first order,
at least over algebraically closed fields of arbitrary characteristic. Such a (nonsingular) subvariety
 will be called an approximate manifold.
The dimension of an approximate manifold  is $\tau$.

The method of approximate manifolds, which he formulates slightly differently and calls the directrix, was  introduced by Hironaka
\cite{H2}.  
\vskip .2truein
\noindent 
{\bf Resolution of an ideal.} 
In this section, we give an overview of our proof of principalization of  an
ideal sheaf ${\mathcal I}$ on a nonsingular 3-fold $V$. 
We adapt the resolution algorithm of Levi \cite{Le} and Zariski  \cite{Z4}, and 
Hironaka's resolution invariant and termination proof in \cite{H2} to this situation.
The main theorem is stated precisely in Theorem \ref{Theorem6} of 
Section \ref{Section4}.
We use the 
philosophy of Bravo, Encinas and Villamayor in  \cite{EV} and \cite{BEV},
resolving by taking the weak transform of an ideal rather than the strict
transform. Our main resolution invariant is then the order of an ideal, rather
than the multiplicity. This leads to vast simplifications in the proofs of resolution theorems.

Let $r=r({\mathcal I})$.
 We will
construct 
 a sequence of blow ups of points and nonsingular curves $W_n$
\begin{equation}\label{In1}
\cdots\rightarrow V_n\rightarrow V_{n-1}\rightarrow \cdots\rightarrow V
\end{equation}
such that if ${\mathcal I}_n$ is the weak transform of ${\mathcal I}_{n-1}$ on $V_n$, 
then $W_n$ is contained in $\mbox{Sing}_r({\mathcal I}_n)$. 
Further, we will have that
two dimensional irreducible components of $\mbox{Sing}_r({\mathcal I})$
are isolated, and at a point of such a component, ${\mathcal I}$ is locally principal, and is generated by $f^r$
where $f=0$ is a local equation of a nonsingular surface $S$. We may thus assume that  2 dimensional components 
never appear in $\mbox{Sing}_r({\mathcal I}_n)$, as they may be eliminated by blowing up the nonsingular surface $S$.

Our goal is to construct a
sequence which terminates after a finite number of steps with $\mbox{Sing}_r({\mathcal I}_n)=\emptyset$.

As shown in Section \ref{SecERIII}, without too much difficulty, we can reach a $V_n$ where
$\mbox{Sing}_r({\mathcal I}_n)$ is  a finite set of points,
so we may as well assume that $\mbox{Sing}_r({\mathcal I})$ is  a finite set of points.  We then apply 
(in Section \ref{SecERIII}) the following algorithm to construct (\ref{In1}).

\begin{enumerate}
\item[1.] If $\mbox{Sing}_r({\mathcal I}_n)$ contains a nonsingular curve $C$, then blow up $C$.
\item[2.] If $\mbox{Sing}_r({\mathcal I}_n)$ is a finite set of points,
then blow up a point in $\mbox{Sing}_r({\mathcal I}_n)$.
\end{enumerate}

We will show that the sequence (\ref{In1}) is finite, and terminates with $\mbox{Sing}_r({\mathcal I}_n)=\emptyset$.

The number $\tau$ can never go down in such a sequence, so we are reduced to
 assuming that there is an infinite sequence (\ref{In1}), with points
$q_n\in V_n$,  such that $q_n$ maps to $q_{n-1}$ for all $n$, and $\nu_{q_n}({\mathcal I}_n)=\nu_{q}(\mathcal I)=r$,  $\tau(q_n)=\tau(q)$ for all $n$.

With our assumptions, in a neighborhood of $q_n$, $\mbox{Sing}_r({\mathcal I}_n)$ is an isolated point,
or a nonsingular curve contained in the exceptional divisor.  We can thus
find approximate manifolds $D_{q_n}$ containing  the germ of 
$\mbox{Sing}_r({\mathcal I}_n)$ at $q_n$.   Let our blow up be $\pi_{n+1}:V_{n+1}\rightarrow V_n$.
By the theory of approximate manifolds,
we have that  $\pi_{n+1}^{-1}(q_n)\cap D_{q_{n+1}}$ is the intersection of 
$\pi_{n+1}^{-1}(q_n)$ and the strict transform of $D_{q_n}$.

We now distinguish several cases, according to the value of $\tau$.
The case $\tau(q)=3$ is trivial. In this case $q$ is an approximate manifold and the strict transform of $q$ under the blow up of $q$ is the empty set,
so $\pi_1^{-1}(q)\cap \mbox{Sing}_r({\mathcal I}_1)=\emptyset$.

In a similar way, if $\tau(q)=2$ and there exists a (nonsingular) curve $C$ through $q$, then $C$ is an approximate manifold. The blow up $\pi_1$ of $C$ then satisfies $\pi_1^{-1}(q)\cap \mbox{Sing}_r({\mathcal I}_1)=\emptyset$.

The case of $\tau(q)=2$ and $q_n$ isolated in $\mbox{Sing}_r({\mathcal I}_n)$ for all $n$ is the first nontrivial case.  In this case we can actually find a manifold
(curve) of maximal contact, by the process of well preparation.
This is a generalization of the method  used by Newton, to compute the branches of an algebraic function
of two variables. This case is worked out in Section \ref{Sect2}. The simpler case of an ideal sheaf  on a nonsingular surface is
in Section \ref{SecERS}.

The most difficult case is when $\tau(q)=1$. This case is in Section \ref{Sect1}.
We start with a regular system of parameters $x,y,z$ at $q$, where $z=0$ is an 
 approximate manifold of ${\mathcal I}$ at $q$,  and make a formal change of variables, replacing $z$
with $z'=z-\Psi(x,y)$ for some appropriate series $\Psi(x,y)$, to well prepare
(Subsection \ref{Wp}).  We will have that $z'=0$ is an approximate manifold for ${\mathcal I}$ at $q$.
Even if we start with regular parameters
$x,y,z$ in ${\mathcal O}_{V,q}$, our resulting well prepared parameter
$z'$ might be formal, existing only in $\hat {\mathcal O}_{V,q}$. 
We may also need to very well prepare (Subsection \ref{VWp}), which involves making a change of
variables replacing $y$ with an expression $y'=y-\phi(x)$, for some appropriate series $\phi(x)$.

   We associate  elements $\Omega(q_n)=(\beta_n,\frac{1}{\epsilon_n},\alpha_n)\in {\bf Q}_+^3$ to $q_n$,
in an ordered set, and show that $\Omega(q_{n+1})<\Omega(q_n)$ for all $n$,
and that $\Omega(q_n)$ cannot decrease indefinitely. This is sufficient to show that (\ref{In1})
has finite length.

Let $R_n$ be the completion of ${\mathcal O}_{V_n,q_n}$, and $I_n={\mathcal I}_{n,q_n}R_n$.

We inductively construct a series of (formal) approximate manifolds $D_n$ in $\mbox{Spec}(R_n)$ by 
taking a local equation of the strict transform of $D_{n-1}$, and modifying
if necessary by very well preparation.  Let $x_n,y_n,z_n$ be the resulting regular parameters in $R_n$,
where $V(z_n)=D_n$. We insist that our parameters $x_n,y_n,z_n$
are such that  
$$
\mbox{Sing}(I_n) \subset V(x_ny_n,z_n)=V(x_n,z_n)\cap V(y_n,z_n)
$$ for all
$n$.  This may prevent very well preparation, if $\mbox{Sing}(I_n) = V(y_n,z_n)$.

$\Omega$ will depend on our choice of regular parameters $x_n,y_n,z_n$ in $R_n$. 
We construct a polygon $|\Delta|=|\Delta(I_n,x_n,y_n,z_n)|$ contained in the real upper half plane, from 
a projection of the nonzero coefficients of monomials in $x_n,y_n,z_n$ in $R_n$ of elements $f\in I_n$.

We choose the parameters $x_n,y_n,z_n$ in $R_n$ to obtain the smallest polygon possible
under substitutions $z'=z_n-\Psi(x_n,y_n)$ (to well prepare) or 
under substitutions $z'=z_n-\Psi(x_n,y_n)$  and $y'=y_n-\phi(x_n)$ (to very well prepare).
 
$(\alpha_n,\beta_n)$ is the upper left hand corner of $|\Delta|$, and $-\epsilon_n$ is the slope of the non vertical boundary line segment through $(\alpha_n,\beta_n)$.

There are only 4 possible types of maps $R_n\rightarrow R_{n+1}$
which can occur in (\ref{In1}), which are denoted by  Tr1 - Tr4. In each case, $\mbox{Sing}(I_n)$ is a curve or point
(Subsection \ref{SubsecPB}).

Most of the cases of Tr1 - Tr4 are just monomial substitutions in the parameters,
except for the case $\eta\ne0$ of Tr1. Local equations of Tr1 are recalled
two paragraphs below. 

The monomial substitutions cause combinatorially simple transformations
of the Newton polygon, and do not involve characteristic $p$ phenomena.

The really interesting case, and only non monomial substitution, is when a point is blown up by a Tr1 transformation,
with $\eta\ne 0$ (Theorem \ref{TheoremRSp27}). Local equations of the map are then $x_n=x', y_n=x'(y'+\eta), z_n=x'z'$.
The blow up is not a monomial substitution in $x_n,y_n,z_n$ in this case, but becomes a monomial substitution after
a  linear change of variables (a ``translation'') $y_n'=y_n-\eta x_n$.
Then the proof becomes a detailed study of the binomial theorem in positive
characteristic.  

In characteristic zero the proof is much simpler and reduces to the statement
that the term $na^{n-1}b$ is nonzero in the expansion
$$
(a+b)^n=a^n+na^{n-1}b+\cdots.
$$
This observation was used by
Zariski in his proof \cite{Z1} of resolution
of algebraic surface singularities, and in his subsequent proofs of resolution.
\vskip .2truein
\noindent{\bf Flags.} Hauser \cite{Hau1} gives a beautiful formulation of Hironaka's invariant used in the proof of termination of the resolution algorithm, which we will
outline here, as it makes  the construction more intuitive.
Hauser assumes that ${\mathcal I}$ is the ideal sheaf of a reduced surface.

We define a flag through $q_n\in \mbox{Sing}_r({\mathcal I}_n)$ to be  the germs
at $q_n$  of a nonsingular surface $F_2$ in $V_n$ and a nonsingular curve $F_1$ in $F_2$.
We assume that $F_1$ is transversal to an approximate manifold $D_{q_n}$ at $q_n$.
If the connected component of $\mbox{Sing}_r({\mathcal I}_n)$ is a (nonsingular)
curve $C$, then we require that either $F_1$ and $F_2$ are both transversal to $C$, or $F_1$ is transversal to $C$ and $F_2$ contains $C$.

We say that parameters $x_n,y_n,z_n$ at $q_n$ are subordinate coordinates for the flag if $x_n=0$ is a local equation of $F_2$ and $x_n=y_n=0$ is a local equation
of $F_1$.  All subordinate coordinates for the flag are (essentially) obtained
from $x_n,y_n,z_n$ by making changes of variables
$$
z_n'=z_n-\Psi(x_n,y_n), \mbox { and }y_n'=y_n-\phi(x_n),
$$ 
for appropriate series $\Psi$ and $\phi$. These are the type of change of variables which we consider in the local resolution algorithm, to well prepare and very well prepare. Hauser shows that  $\Omega(q_n)$ then becomes
an actual invariant of the flag. We maximize $(\beta,\frac{1}{\epsilon},\alpha)$
for subordinate changes of variables for the flag, by well preparing and very well preparing.

Suppose that we have a sequence of points $\{q_n\}$ in (\ref{In1}) with
$q_n\in \mbox{Sing}_r({\mathcal I}_n)$ and $q_n$ mapping to $q_{n-1}$ for all $n$.

Hauser shows that a flag at $q_{n-1}$ transforms to a flag at $q_{n}$ under the sequence (\ref{In1}),  that we can always choose  subordinate coordinates
so that $R_{n-1}\rightarrow R_n$
is monomial in these coordinates, and the invariants computed from these flags
satisfy  $\Omega(q_n)<\Omega(q_{n-1})$. The descending invariants $\Omega(q_n)$
on the sequence are  then computed by
taking the transforms of an initial choice of a flag at $q$.

\vskip .2truein

\noindent{\bf Embedded resolution.}

The concepts of simple normal crossings (SNCs) and of transversal intersection are  defined in Section \ref{SecPre}.
It is discussed in more depth in \cite{C2}, especially in Chapter 4 of \cite{C2}. The geometric idea is that all components
are nonsingular and intersect transversely. Our notion of simple normal crossings is sometimes called
``strict normal crossings'' in the literature. 

To give the full proof of embedded resolution of singularities of an ideal sheaf ${\mathcal I}$
on a nonsingular 3-fold $V$ (Theorem \ref{Theorem6}) we must start with a SNC divisor $E$ on
$V$ as
well as our ideal sheaf ${\mathcal I}$. We then want to construct a sequence (\ref{In1}) so that not only is
$\mbox{Sing}_r({\mathcal I}_n)=\emptyset$,  but the divisor $E_n$ on $V_n$
consisting of the sum of the exceptional divisor of $V_n\rightarrow V$ and
the preimage of $E$ on $V_n$ is a SNC divisor. 
 Fortunately, essentially no new 
characteristic $p$ problems occur in this proof, beyond those we have already encountered.

 A little care must
be taken however in making sure that all blow ups are of centers which are
transversal to the exceptional divisor. The problem is that we cannot deviate from the resolution algorithm, which is to blow
up a curve if there is a curve in the singular locus (we can ensure that the curve will always be nonsingular after a little preliminary blowing up
before starting the algorithm), and to blow up a point
if there is no curve in the singular locus. Hironaka's invariant which we use to show termination does not behave well 
 if there is deviation from this
procedure.

In order for $E_n$ to be a SNC divisor on $V_n$, we require that all blow ups $V_{i}\rightarrow V_{i-1}$
(with $i\le n$) be permissible (Definition \ref{Def101}). This means that they are the blow up of a nonsingular subvariety of  $\mbox{Sing}_r(\mathcal{I}_{i-1})$ which is transversal to $E_{i-1}$.

To accomplish this, we make use of the $\eta$ invariant. This invariant is used by Abhyankar in his 
``good point'' proof
of embedded resolution of (characteristic zero) surface singularities, which is presented in \cite{Ab10}, in the lectures \cite{O} of Orbanz, and in Lecture 3 of  Lipman's Arcata notes \cite{L1}. We define $\eta$ as follows.

On $V$, set $E^+=\emptyset$ and $E^-=E$.
At a point $q_n\in V_n$ where $\nu_{q_n}({\mathcal I}_{n})=\nu_{q}({\mathcal I})=r$, define $E_n^-$ to be the strict transform of $E^-=E$, and $E_n^+$ to be the preimage of the exceptional divisor of $V_n\rightarrow V$.  
If $\nu_{q_n}({\mathcal I}_{n})$ drops to less than $r$, we define 
$E_n^-$ to be $E_n$ and $E_n^+$ to be the empty set.

For $q_n\in\mbox{Sing}_r({\mathcal I}_n)$, let $\eta(q_n)$ be the number of irreducible components of $E_n^-$ containing $q_n$.
We have $0\le \eta(q_n)\le\eta(q_{n-1})\le 3$ if $\nu_{q_n}({\mathcal I}_{n})=\nu_{q_{n-1}}({\mathcal I}_{n-1})=r$.

In Section \ref{SecERIII}, we reduce by descending induction on $\eta$  to the case of $\eta=0$
($E_n^-$ is disjoint from $\mbox{Sing}_r({\mathcal I}_n)$) and $\mbox{Sing}_r({\mathcal I}_n)\subset E_n^+$. This is achieved by noticing
that our resolution problem reduces to a resolution problem on a nonsingular
subvariety of $V$ of dimension $3-\eta(q)\le 2$ (an intersection of components of $E^-$).

Further, we show that
 the exceptional divisor of $V_n\rightarrow V$ is transversal to an approximate hypersurface of ${\mathcal I}_n$.

Now any blow up that we make in our algorithm (given after (\ref{In1})) to reduce the
order of ${\mathcal I}_n$ will automatically be permissible. The point is that
one dimensional components of $\mbox{Sing}_r({\mathcal I}_n)\cap E_n^-$ are not well controlled, while those of $\mbox{Sing}_r({\mathcal I}_n)\cap E_n^+$ are projective lines transversal to $E_n^+$.
\vskip .2truein
\noindent{\bf Further generalization of the algorithm.} The most interesting generalization would be to make the above reduction argument
work for a hypersurface $X$ in a nonsingular variety $V$ of dimension $d\ge 4$.
The major problems occur in making the local reduction argument work for small $\tau$. The cases of
\begin{equation}\label{In2}
0\le  d-\tau({\mathcal I}_{X,q})\le 2
\end{equation}
are analogous to the cases of $\tau({\mathcal I}_{q}) = 3,2,1$ which we analyze in this paper for ideal sheaves on
a nonsingular variety of dimension 3.
However, the case of $d-\tau({\mathcal I}_{X,q})>2$ appears to be very difficult.
\vskip .2truein

\noindent{\bf Resolution of excellent surfaces.} Hironaka  announced in his notes \cite{H2} a proof of resolution
of singularities of an arbitrary excellent surface. In the notes
he gives a few comments about the realization of the general proof, and gives some parts of the proof in \cite{H3}, \cite{H4} and \cite{H5}  (see also  Cossart's paper \cite{Co1}). 

We outline this approach, when
$S$ is a surface which is  embedded in a nonsingular variety $V$  of
arbitrary dimension $d$, over an algebraically closed field ${\bf k}$. 

Making the proof work over a non perfect base field, or for an arbitrary excellent surface requires an excursion into a further series of complexities
and treacherous pitfalls which we will not address here. Some details on this
are given in \cite{H4}, \cite{H5}, \cite{H6} and \cite{Co1}. A complete proof of resolution of excellent surfaces using a different method, is given by Lipman in \cite{L2}.

When $d=3$, the  order and multiplicity of ${\mathcal I}_S$ are the same. Also, under permissible blow ups,
the weak transform and strict transforms of ${\mathcal I}_S$ are the same. For larger $d$, however, this is not the case. We now take $r$ to be the largest multiplicity of ${\mathcal I}_S$, and define $\mbox{Sing}_r({\mathcal I})$
to be the locus of points of maximal multiplicity. Under permissible transforms, we take ${\mathcal I}_n$ to be the strict transform of ${\mathcal I}_{n-1}$.  ${\mathcal I}_{n}$ is the ideal sheaf of the strict transform $S_n$ of $S$ on $V_n$.

The first part of the proof when $d=3$ (in ``Resolutions of ideals'' above) now goes through without trouble for general $d$.
We reduce to the assumption that the locus of maximal multiplicity $\mbox{Sing}_r({\mathcal I})$ is a finite set of points, and we have an infinite sequence of blow ups of points and nonsingular curves (\ref{In1}),
constructed by following the algorithm stated after ({\ref{In1}).
As in the analysis for the case $d=3$, we reduce to consideration of a sequence of points $q_n\in S_n$
such that $q_n$ maps to $q_{n-1}$ for all $n$, and the multiplicity of $S_n$ at $q_n$ is $r$ for all $n$. 

Let $R=\hat{\mathcal O}_{V,q}$, with maximal ideal $m$, and let $I={\mathcal I}_{S,q}R$.

We can obtain any value  $\tau(q)=\tau(I)\ge 1$, and 
as we have commented above in ``Further generalization of the algorithm'', our
local  argument for the reduction of the order of an ideal  does not extend to
$d-\tau(I)>2$.  We can however get around this by using $\overline\tau$
 instead of 
$\tau$, which we now define. We are in fact already forced to use this new concept, since  when we are using
multiplicity and taking strict transform, instead of  using the simpler concepts of order and weak transform.
We point out that Hironaka defines $\tau$ to be what we call $\overline\tau$ in this paper, and does not consider what
we call $\tau$, although $\tau=\overline\tau$ for locally principal ideals.

Let  $\mbox{in}(I)$ be the ideal of $\mbox{gr}_{m}(R)=\oplus_{j\ge 0} m^j/m^{j+1}$ generated by  initial forms of 
elements of $I$.
We define $\overline\tau(I)$ to be the dimension of the smallest $R/m$ linear subspace $T$ of $\mbox{gr}_{m}(R)$ such that the initial forms of
all elements of $I$ are contained in the polynomial ring $R/m[T]\subset \mbox{gr}_{m}(R)$.  

Since $\mbox{in}(I)$ is contained in the ideal $(T)$ of $\mbox{gr}_m(R)$ generated by
$T$, we see that 
$$
2=\mbox{dim }R/I=\mbox{dim }\mbox{gr}_{mR/I}(R/I)
=\mbox{dim }\mbox{gr}_{m}(R)/\mbox{in}(I)\ge\mbox{dim }\mbox{gr}_{m}(R)/(T).
$$
Thus $\overline \tau(I)\ge \mbox{dim }R-2$.

The problem becomes to extend the reduction argument used in ``Resolution of an ideal'' in the case  when $d=3$, using  $\overline\tau$ instead of $\tau$. This is feasible  as $\overline\tau(I)\ge d- 2$ for the ideal of a surface. The price of this is that we must then use the full range of sophisticated methods developed by Hironaka in \cite{H1} to control the reduction of
multiplicity of the strict transform of an ideal under a permissible blow up
(instead of just the order of the weak transform of an ideal),
and requires a more 
sophisticated analysis of the polyhedra associated to a singularity
which incorporates information from an appropriate standard basis of $I$.
Some of this analysis is done in \cite{H3}, in a very general setting. 
 
\vskip .2truein

\noindent{\bf Abhyankar's proof of embedded resolution of surfaces.}
We conclude this section with a quick overview of Abhyankar's proof in \cite{Ab7} and \cite{Ab5}, \cite{Ab6}, \cite{Ab8}, \cite{Ab9}.

The main resolution theorem in Abhyankar's proof is Theorem 1.1 of \cite{Ab6}, 
which is restated in Theorem 5.1 in  \cite{Ab7}.

Suppose that ${\bf k}$ is an algebraically closed field of characteristic $p$, and
$A={\bf k}[x,y]_{(x,y)}$ is a regular local ring. Suppose that
\begin{equation}\label{AP1}
f(Z)=Z^n+f_1Z^{n-1}+\cdots+f_n
\end{equation}
is a monic polynomial with coefficients in $A$, and $\omega$ is a
rank 1 nondiscrete valuation which dominates $A$. 
Let
\begin{equation}\label{AP2}
A\rightarrow A_1\rightarrow \cdots\rightarrow A_t\rightarrow \cdots
\end{equation}
be the sequence of blow ups of the maximal ideals of the local rings $A_i$
which are dominated by $\omega$. Theorem 1.1 of \cite{Ab6} states that $f(Z)$ has a nice form
for all $t\gg 0$. The idea of an approximate manifold is implicit in this argument.

Local uniformization can be easily deduced from Theorem 1.1 of \cite{Ab6}. The essential case is
when our germ of a singular surface $S$  has the local equation $f(Z)=0$, and
the $f_i$ have order $\ge i$. Let $\nu$ be a valuation which dominates $f(Z)=0$.
Let $\omega$ be the restriction of $\nu$ to $A$. The most difficult case where 
the worst complexities of characteristic $p$ occur is when $\omega$ is  nondiscrete of rational rank 1. We may as well assume then that $\omega$ satisfies the assumptions of Theorem 1.1 of \cite{Ab6}. We then easily deduce from the theorem that for $t\gg0$, there exists
a substitution in $Z$ which gives a birational extension such that the center of
$\nu$ has multiplicity less than $n$ on the strict transform of $S$.

Local uniformization is discussed in ``Valuation rings and local uniformization'' in Section \ref{Intro3}.
In his characteristic zero proof of resolution of surface singularities in \cite{Z4}, Zariski introduces the notion of the dominant character of a
normal sequence. This is the statement that if embedded local uniformization 
 is true by possibly quite different local resolution
sequences along each valuation centered on a two dimensional  hypersurface, 
then a reasonable global algorithm of resolution will lead to a global resolution. Zariski proves this in \cite{Z4}.

Abhyankar generalizes this notion to resolution of ideals supported on
an arbitrary characteristic  nonsingular 3-fold, and requires that the resolution procedure generate a resolved object
which makes SNCs with the exceptional divisor. This is accomplished in
a very difficult proof in  
Chapter 1 of \cite{Ab7}. In
this way, Abhyankar deduces in \cite{Ab7} the embedded resolution theorems
for surfaces, stated  as  Theorems \ref{Theorem1} and \ref{Theorem23} in this paper.

The proof of Theorem 1.1 of \cite{Ab6} is realized in the papers \cite{Ab6}
and \cite{Ab9}, with the aid of a technical result from \cite{Ab8}.

The idea is to consider a splitting field $L$ of the polynomial $f(Z)$ over
the quotient field $K$ of $A$.  In \cite{Ab6}, making ingenious use of
the norm and discriminant of a field extension, as well as the observation that the coefficients of a monic polynomial are the elementary symmetric functions in its roots, the proof of Theorem 1.1 of \cite{Ab6} is reduced to the case where $f(Z)$ is irreducible in $K[Z]$, $n=p^m$ for some integer $m$, and the valuation $\omega$ does not split under extension to $L$.
The point is that this isolates the part of the resolution procedure where
characterstic $p$ is an obstacle. Now that we have made the reduction,
we have obtained nonsplitting type and ramified type for $f(Z)$
(Section 2 of \cite{Ab9}).  Lemma 2.6 \cite{Ab6} reduces the proof of  Theorem 1.1 of \cite{Ab6} in this case to the technical result of  Theorem 5.3 of \cite{Ab9}.  The goal is to find $t$ sufficiently large in 
(\ref{AP2}) such that $f_n$ has small positive order, and that this cannot
be changed after translation of $Z$ by an element of  $A_t$.
Recall that $n$ is a power of $p$. If we translate by an element $r$ of $A_t$, replacing $Z$ with $Z+r$, we get
$$
Z^n+ f_{1,t} Z^{n-1}+\cdots + f_{n,t}
$$
where $ f_{i,t}\in A_t$, and
$$
 f_{n,t}=f_n+r^n+\mbox{ contributions from the coefficients of $Z^i$ with $0<i<n$.}
$$
The non splitting type and ramified type assumptions tell us that
the contributions from the coefficients of $Z^i$ with $0<i<n$ cannot interfere too much in low degree with $f_n+r^n$. We write
$$
f_{n,t}=x_t^{a_t}y_t^{b_t}H_t(x_t,y_t)
$$
for appropriate regular parameters $x_t,y_t$ in $A_t$, where
$x_t,y_t$ do not divide $H_t$. Let  
$$\mbox{ord}(H_t(0,y_t))=c_t.$$
We want to find a $t$ such that $c_t$ is small, and we cannot remove the term $x_t^{a_t}y_t^{b_t+c_t}$ from
$f_{n,t}$ by translation of $Z$. We may reduce to the case where the $x_t^{a_t}y_t^{b_t+c_t}$ term cannot be removed if at least one of $a_t$ and $b_t+c_t$ is not divisible by $n$.
Sections 6 through 9 of \cite{Ab9} give an algorithm for achieving this.
The idea is that after translation of $Z$ to remove appropriate monomials in $x_t$ and $y_t$ which are powers of $n$, if we blow up we have that $c$ can only go up under a very special circumstance (this is where Lemma 27 of \cite{Ab8} comes in), and after enough blow ups, there must actually be a drop in $c_t$.

The problem with order going up under resolution also occurs in resolution
of vector fields (Seidenberg \cite{Se}, Cano \cite{Ca} and Panazzolo \cite{Pa}) and in monomialization of morphisms (Cutkosky \cite{C1}, \cite{C3}; and   Cutkosky
and Piltant \cite{CP}).
In the case of dimension two (as in the case of resolution of surfaces) this problem is not overwhelming, but  the level of complexity goes up tremendously in dimension three (\cite{Ca}, \cite{C1}, \cite{C3}).

To compare Abhyankar's proof with Hironaka's invariant  used in this paper, we point out that the well preparations of Section \ref{Sect1} of this paper are
essentially Abhyankar's translations, and the numbers $(\alpha,\beta)$
of Section \ref{Sect1} of this paper are analogous to 
Abhyankar's $(a,b+c)$.

\section{On   Resolution of 3-folds}\label{Intro3}

In this section we give some indication of the proofs of the different
steps, discuss some of the history,  review some of the
background material, and point out a  few subtle considerations which come up
in the proof.  We translate
Abhyankar's proof into  language which will be more easily accessible to
the reader familar with the language of schemes.  We have also simplified some
of the proofs, partly by simplifying the statements to the case of an
algebraically closed ground field, which is all that we need.
We assume throughout this section that ${\bf k}$ is an algebraically closed field.

\vskip .2truein

\noindent{\bf Step 2.}
In Section \ref{SecProj}, we prove the projection theorem, Theorem \ref{Theorem19}.
The simple idea of the proof is as follows. Suppose that
 $V_0$ is a projective variety of dimension $d$, with function field $K={\bf k}(V_0)$, and 
$q\in V_0$. The first observation (in Theorem \ref{Theorem47}) is that the multiplicity $e({\mathcal O}_{V,q})$ of the local ring ${\mathcal O}_{V,q}$ is less than or equal to the degree $\mbox{deg}(V_0)$ of $V_0$.

The idea of the proof of Theorem \ref{Theorem19} is to project from points of high multiplicity until we get
a rational map from $V_0$ to a variety $V_n$ such that $V_0$ and $V_n$ have
the same dimension, and if $[V_0:V_n]$ is the degree of the 
extension of function fields ${\bf k}(V_0)$ over ${\bf k}(V_n)$, then 
every point of $V_n$ has multiplicity $\le \frac{d!}{[V_0:V_n]}$.
Then we take the normalization  of $V_n$ in ${\bf k}(V_0)$, to get
a  variety $V$  birational to $V_0$, such that  $\mbox{deg}(V)\le d!$, and thus every point of $V$
has multiplicity $\le d!$.

Albanese \cite{Alb} used this method in  resolution of algebraic surfaces. Abhyankar gives an algebraic proof of this theorem, valid in 
all dimensions in Section 12 of his book \cite{Ab7}.
Lipman also gives a survey of  the proof 
in all dimensions, in geometric language closer to that of Albanese, in Section 5 of
Lecture 1 of his Arcata notes \cite{L1}.

In Section \ref{SecProj}, we essentially present Lipman's argument in his survey \cite{L1}, but fill in a few details.

Zariski analyzes Albanese's proof on pages 21 and 22 of his 1935 book \cite{Z6}.
Zariski points out that Albanese's proof makes extensive use of intersection
properties of curves on a surface, which are well defined for nonsingular
surfaces, but for surfaces with singularites, present new difficulties.
Zariski cautions us that ``with this distinction in mind, Albanese's proof should be closely scrutinized''.  

The first instance where this issue is apparent is in the necessity to
use multiplicity instead of the more intuitive  local intersection number
$(H^d\cdot V_0)_q$ 
at a point $q\in V_0$, obtained by computing the length of
${\mathcal O}_{V_0,q}/(x_1,\ldots,x_d)$, where $x_i=0$ are local equations
of general  hyperplane sections of $V_0$ through $q$.  
This invariant is related to multiplicity by the inequality 
$$
e({\mathcal O}_{V_0,q})\le (H^d\cdot V_0)_q,
$$
 with equality if and only
if $V_0$ is Cohen Macaulay at $q$ (Theorem 23, Chapter VIII, Section 10 \cite{ZS}).
Intuitively, we might expect that $(H^d\cdot V_0)_q$ is the correct  number to
use (instead of the multiplicity $e({\mathcal O}_{V_0,q})$) in the projection formulas of Theorem \ref{Theorem47}. However, we see that this naive approach fails in general.

Another subtle point which comes up  is the necessity
of using a sophisticated B\'ezout theorem in the proof to show that we do not
eventually reach a situation where the dimension of the image drops
(in Theorem \ref{Theorem48}). This issue only appears when the dimension $d$
of $V_0$ is $\ge 3$.  Abhyankar proves the theorem we need in Theorem 12.3.1
of \cite{Ab7}. A slightly more general version of this theorem appears in Fulton's 
book \cite{F} (Example 12.3.1 of the Refined B\'ezout Theorem).  The only instance where characteristic $p$ appears as an
issue in the proof of Theorem \ref{Theorem19}, is in the use of Bertini's theorem for general hyperplane sections, again in the proof that the dimension of the image does not drop
(Theorem \ref{Theorem48}).  
\vskip .2truein

\noindent{\bf Valuation rings and local uniformization.}
Suppose that 
$K$ is an algebraic function field 
over ${\bf k}$,  and $V$ is a valuation ring of $K$ (containing ${\bf k}$),
whose residue field is ${\bf k}$.
Suppose that $X$ is a projective (or even proper) variety with function field $K$. Then there is a unique closed point $p\in X$ such that $V$ dominates (Section \ref{SecPre}) the local ring ${\mathcal O}_{X,p}$. 

The problem of local uniformization is to show that for every valuation ring $V$
of $K$, there exists a regular local ring $R$, which is essentially of finite
type over ${\bf k}$ and has $K$ as its quotient field, such that $V$ dominates $R$.
 
Zariski proved local uniformization over characteristic zero function fields in \cite{Z2}. 
The most difficult problems in local uniformization occur for valuations whose value groups are not finitely generated.
\vskip .2truein

\noindent{\bf The Zariski Riemann manifold.}
Zariski introduced the idea of local uniformization along a valuation
to give a local formulation of resolution of singularities.

Suppose that 
$K$ is an algebraic function field. 
Zariski constructed a ringed space $\mbox{ZR}(K)$ which we call the Zariski Riemann manifold
of $K$. A nice survey of this concept is given by Lipman in his appendix to
Chapter II of \cite{Z6}. The points of $\mbox{ZR}(K)$ are the
valuation rings of the function field $K$, which contain ${\bf k}$ and whose residue field is ${\bf k}$  \cite{Z5}. Such valuations are called zero dimensional. The local ring of the point is
the valuation ring. Zariski proved that this space is quasi compact in \cite{Z5} and Section 17, Chapter VI of \cite{ZS}.

There is a natural continuous morphism of ringed spaces $\Phi_X:\mbox{ZR}(K)\rightarrow X$ for any projective variety $X$ whose
function field is $K$. A valuation ring $A$ of $K$ maps to the point of $X$
whose local ring is dominated by $A$. The point is called the center of $A$ on $X$. The mappings $\Phi_X$ are natural, so that if $\Psi:X_1\rightarrow X$ is a 
birational map, then $\Psi\circ\Phi_{X_1}=\Psi_X$.

Because of the quasi compactness of
$\mbox{ZR}(K)$, the construction of a resolution of singularities  will follow
from patching together a finite number of local resolutions.  Zariski
accomplished this for (characteristic zero) surfaces \cite{Z3} and 3-folds \cite{Z4}.
\vskip .2truein

\noindent{\bf Step 3.}
In Section \ref{SecSmall}, we generalize the method of Jung \cite{J}, for locally resolving
singularities.  This method is discussed 
(for surfaces in characteristic zero) by Zariski on pages 16 and 17 of \cite{Z6}.
This method is also discussed by Lipman in  Lecture 2 of \cite{L1}, and by
Giraud in \cite{G2}. Abhyankar has generalized the method to the case of
tame ramification in positive characteristic, and has shown that the 
situation is extremely complicated when the ramification is wild \cite{Ab2}.

Our proof is a simplification of the proof by Abhyankar
given in ``Uniformization of points of small multiplicity'', in
Section 10 of Chapter 3 of \cite{Ab7}. Our final statement is a proof in Theorem \ref{Theorem18} of local uniformization for the function field of $V$.

In characteristic zero the method is very simple and elegant. The fundamental result is the following. Suppose that 
$f:X\rightarrow Y$ is a finite map from a normal variety $X$ (of arbitrary dimension) to a nonsingular variety $Y$. Then by the purity of the branch locus,
the branch locus of $f$ is a divisor $D_f$  on $Y$. If $D_f$ is a SNC divisor,
then $X$ must have abelian quotient singularities, which have a very simple
form, and can be easily resolved.  In fact, with respect to the
toroidal structures determined by $D_f$ and $f^{-1}(D_f)$, $f:X\rightarrow Y$
is a morphism of toroidal varieties. We can thus use one of the combinatorial
algorithms (such as in \cite{KKMS} or \cite{BrM2}) to resolve the singularities of $Y$.

This gives us a way to deduce local uniformization (assuming characteristic zero) in
dimension $d$, assuming embedded resolution in dimension $d-1$.  We take a 
singular variety $V$, and a
(zero dimensional) valuation ring $A$ of the function field of $V$, whose center is a (closed) point $p$ on $V$.
We choose a  finite map $f$ from an open neighborhood $X$ of
$p$ in $V$ to an open subset $Y$ of the affine space 
${\bold A}_{\bold k}^d$. Then we apply embedded resolution of singularities
to construct a birational, projective morphism  $\pi_1:Y_1\rightarrow Y$
so that $\pi_1^{-1}(D_f)$ is a SNC divisor on the nonsingular variety $Y_1$.
Finally, we take $X_1$ to be the normalization of the  variety $(X\times_YY_1)_{\mbox{red}}$. $X_1$ is birational to $X$. Let $p_1$ be the center of $A$ on $X_1$. Then by the above, $X_1$ has a toric singularity at
$p_1$, which we can now easily resolve.

In Section \ref{SecSmall}, we extend this argument to the case when
${\bf k}$ has characteristic $\ge 0$, and the local ring ${\mathcal O}_{V,p}$
has multiplicity less than $p$.

Our main local resolution statement is Theorem \ref{Theorem18}, which proves local uniformization for the
function field of $V$.

A slightly different  formulation of the theory of ramification from the one we present in Sections \ref{SecRam} and \ref{SecSmall},  written in  the language of etale mappings, is given by Grothendieck and Murre in \cite{GM}, especially in Section 2.3, ``Tame Ramification and Abhyankar's Theorem''.

\vskip .2truein

\noindent{\bf Hypersurfaces of maximal contact for singularities of low multiplicity.} Here we outline an alternate proof of Step 3 of Theorem \ref{Theorem82}, which was pointed out to us by one of the reviewers.
We start with a variety  $V'$ such that all points of $V'$ have multiplicity $<p$.
We now make a local study of  $V'$. 

Cossart has shown
(Theorem II.3 \cite{Co3}) that 
a hypersurface of maximal contact always exists for singularities at points
of multiplicity $<p$. For the sake of motivation, we point out that the existence of a hypersurface of maximal contact is easy to see in an important special case, that of a hypersurface singularity
$$
y^r+a_1(x)y^{r-1}+\cdots+a_r(x)=0
$$
where $x$ represents a system of variables, and $r<p$.
In this case we may obtain a hypersurface of maximal contact 
$\overline y=0$, just as in
characteristic zero, by making a Tschirnhaus transformation
$\overline y=y+\frac{a_1(x)}{r}$, to eliminate the $\overline y^{r-1}$ coefficient.

Cossart's proof  makes a detailed local analysis, using information obtained
from the Hilbert Samuel function of the singularity and a theorem of Giraud
\cite{G} to conclude that there is  a hypersurface of maximal contact 
for the singularity. Specifically, Cossart shows that
if $(f_1,\dots,f_m)$ is a local standard
 basis of the singularity, with $\nu_i=$ord$_x(f_i)$, for $1\leq i \leq m$,
 then  $\nu_i<p$, for all $i$, and
 there exists a regular system of parameters $(y_1,\dots,y_{\tau},u_1,\dots,u_d)$
 with $d\leq 3$ such that $W=V(y_1,\dots,y_{\tau})$ has maximal contact,
by Giraud's Theorem \cite{G}.

 Now using techniques of  Hironaka on 
 the Hilbert-Samuel function, the resolution of the singularity
 reduces to the resolution of the idealistic exponent
 ${\mathcal E}=\cap_{1\leq i \leq m}((f_i,\nu_i))$,
 which is supported by $W$ and is equivalent to an idealistic exponent
 $((I,b))$ with $I$ an ideal of ${\mathcal O}_{W,x}$, $b\in {\bf N}$.
 ${\mathcal E}$ is solved by principalizing $I$, which is done
 by a straight forward generalization of Theorem \ref{Theorem6} of this paper to idealistic exponents.

\vskip .2truein
\noindent{\bf Step 4.}
In Step 3, we proved local uniformization for the function field of $V$, so we have reduced, by quasi compactness of the Zariski Riemann manifold as we commented in ``The Zariski Riemann manifold'', to patching together a finite number of local resolutions.   The problem is to take several projective varieties which are birationally equivalent to $V$, and such that every valuation has a nonsingular center on at least one of the varieties.
We must construct a variety $V'$ which birationally dominates all
of these varieties and is nonsingular.
 Zariski's original characteristic zero proof \cite{Z4} made use of
general Bertini theorems which are not valid in positive
characteristic. His method was to study a birational morphism through an associated linear system.  We eliminate the need for Bertini theorems by
using the fact that a birational morphism is the blow up of a sheaf of ideals,
and  using the universal property of blow ups of ideal sheaves 
(Proposition II.7.14 \cite{Ha}, Theorem 4.2 \cite{C2}).
This is essentially  the approach of  Abhyankar's proof of this step, although
he states and proves it in his own language (Section 8 of Chapter 2 of \cite{Ab7}).

This step is where the theorem on principalization of ideals (Theorem
\ref{Theorem23}) is required. It is essential in this application that
we blow up only above the locus where the ideal sheaf is not principal.
Another essential ingredient of the proof
is that a birational morphism of nonsingular  3-folds factors in
a canonical way over the general point of a curve. It is just a sequence of
blow ups of curves which dominate the original curve. This  fails
completely in higher dimensions.

The fourth step is proven in Sections \ref{SecPat} and \ref{SecRes}.

\vskip .2truein
\noindent{\bf Step 5.}
The final fifth  step is accomplished in Section \ref{SecFin}. The original proof of
Abhyankar does not produce a resolution which is an isomorphism over the nonsingular locus. The fact that a resolution can be birationally modified
(in dimension 3) to produce a resolution which is an isomorphism over the nonsingular locus is due to Cossart \cite{Co4}. Our proof of this section is
based on his extremely simple argument. Again, the main point is that
a birational morphism  of nonsingular 3-folds factors in
a canonical way over the general point of a curve.

\section{Preliminaries}\label{SecPre}
In this section we establish notation, and review a few  results that we will use in the course of the proof.

If $R$ is a local ring, we will denote the maximal ideal of $R$ by $M(R)$. Suppose that $S$ is another local ring, and $R\rightarrow S$ is an inclusion of rings.
 We will say that $S$ dominates $R$ if 
 $M(S)\cap R=M(R)$.

Suppose that $\Lambda\subset R$ is a subset. We define the closed subset
$$
V(\Lambda)=\{P\in\text{Spec}(R)\mid \Lambda\subset P\}
$$
of $\mbox{Spec}(R)$.
We will  consider the order $\nu_R(I)$ of an ideal $I$ in $R$.
$\nu_R(I)$ is defined to be the largest integer $n$ such that
$I\subset M(R)^n$.  If $q$ is a point on a variety $W$, and
${\mathcal J}$ is an ideal sheaf on $W$, then we denote
$\nu_q({\mathcal J})=\nu_{{\mathcal O}_{W,q}}({\mathcal J}_{q})$.

\vskip .2truein \noindent{\bf Multiplicity.} Suppose that $R$ is a
$d$ dimensional Noetherian local ring, and that $q$ is an $M(R)$
primary ideal of $R$.  There exists a degree $d$ polynomial
$P_q(n)$ which has the property that the length
$$
\ell(R/q^n)=P_q(n)
$$
for $n\gg 0$. The multiplicity $e(q)$ is defined to be $d!$ times the leading coefficient of $P_q(n)$. The multiplicity of $R$ is defined to be $e(R)=e(M(R))$.
From the definition we infer that if $q'\subset q$ are $M(R)$-primary, then $e(q)\le e(q')$.
Suppose that $R^*$ is the $M(R)$-adic completion of $R$. We see that $e(qR^*)=e(q)$.

The multiplicity has the important property that $e(R)=1$ if and
only if $R$ is regular (Theorem 23, Section 10, Chapter VIII
\cite{ZS}).

\begin{Theorem}\label{Theorem70} Suppose that  $R/M(R)$ is infinite. Then there exists an ideal $q'\subset q$ generated by a system of parameters such that $e(q')=e(q)$.
\end{Theorem}
This is proven in Theorem 22, Section 10, Chapter VIII \cite{ZS}.

\begin{Theorem}\label{Theorem71} Suppose that $R$ is a local domain, and $T$ is an overring of $R$
 which is a domain and a finite $R$ module. Then $T$ is a semilocal ring. Let $p_1,\ldots,p_s$ be the maximal ideals of $T$.
  Let $K$ be the quotient field of $R$ and $M$ be the quotient field of $T$. Then
$$
[M:K]e(q)=\sum_{i=1}^s[T/p_i:R/M(R)]e(qT_{p_i}).
$$
\end{Theorem}
This is proven in Corollary 1 to Theorem 24, Section 10, Chapter VIII \cite{ZS}.
\vskip .2truein

\vskip .2truein \noindent {\bf Intersection multiplicity and
degree.} We use the notation of Fulton \cite{F} for intersection
theory. The intersection class of a $k$-cycle $\alpha$ with a
Cartier divisor $D$ on a scheme $X$ is denoted by $D\cdot \alpha$.
If $D_1,\ldots, D_k$ are Cartier divisors on $X$,    the
intersection number of $D_1,\ldots, D_k$ and $\alpha$ is denoted
by $\int_XD_1\cdot\ldots\cdot D_k\cdot \alpha$.

The degree $\text{deg }Y$ of a subscheme $Y$ of ${\bf P}^m$ of dimension $d$ is $d!$ times the leading coefficient of the Hilbert Polynomial $P_Y(n)$ of the homogeneous coordinate ring of $Y$.
It can also be computed as
$\text{deg Y}=\int_{{\bf P}^m}H^d\cdot Y$
where $H$ is a hyperplane of ${\bf P}^m$.
Proofs are given in Section 7 of Chapter 1 \cite{Ha} and the Corollary of Section 8.5 \cite{I}, or in Examples 2.5.1 and 2.5.2 \cite{F}.

\begin{Lemma}\label{Lemma44}
Suppose that $R$ is a local ring of dimension $d$, and let
$$
Y={\rm Proj}(\oplus_{n\ge 0}M(R)^n/M(R)^{n+1}).
$$
 Then $\text{deg }Y=e(R)$.
\end{Lemma}
 This follows from comparing the Hilbert polynomial of $Y$, which is equal to
 $$
\ell(M(R)^n/M(R)^{n+1})
$$
 for all large $n$, and the  polynomial $P_{M(R)}(n)$.

 \vskip .2truein
\noindent{\bf Completion.} Suppose that $\bold k$ is an
algebraically closed field, and $R$ is a normal local ring, which
is essentially of finite type over ${\bf k}$ (a localization of a
ring which is of finite type over ${\bf k}$). The $M(R)$-adic
completion $R^*$ of $R$ is a normal local ring (Theorem 32,
Section 13, Chapter VIII \cite{ZS}). Suppose that $L$ is a finite
field extension of the quotient field of $R^*$. Then the integral
closure of $R^*$ in $L$ is a local domain (since $R^*$ is
Henselian, Corollary 2 to Theorem 18, Section 7, Chapter VIII
\cite{ZS}).

We also state here the Zariski subspace theorem.  For a proof, see (10.10) \cite{Ab7}. A ring $R$ is essentially of finite type over a ring $T$
if $R$ is a localization of a finitely generated $T$ algebra.

\begin{Theorem}\label{Theorem12} Suppose that $R$ and $S$ are normal local rings which are locally of finite type over ${\bf k}$. Suppose that $S$ dominates $R$.  Let
$R^*$ be the $M(R)$-adic completion of $R$ and $S^*$ be the $M(S)$-adic completion of $S$.
Then the natural morphism $R^*\rightarrow S^*$ is an inclusion.
\end{Theorem}
\vskip .2truein \noindent{\bf Birational geometry.} In this paper,
a variety is an open subset of an integral, closed subscheme of
${\bf P}^n$. A curve, surface or 3-fold is a variety of the
corresponding dimension.

The definition of variety in \cite{C2} is a little more general. In
\cite{C2}, a variety is an open subset of a reduced,
equidimensional subscheme of ${\bf P}^n$. The resolution theorems
for varieties in \cite{C2} are thus valid for reduced
equidimensional  subschemes of  ${\bf P}^n$.

Suppose that $X$ is a projective variety over an algebraically
closed field ${\bf k}$. Let $K$ be the function field of $X$.
Suppose that $V$ is a valuation ring of $K$ which contains ${\bf
k}$. Then there exists a unique (not necessarily closed) point
$p\in X$ such that $V$ dominates the local ring ${\mathcal
O}_{X,p}$. $p$ is called the center of $V$ on $X$.

Suppose that $\Phi:X\rightarrow Y$ is a birational morphism of normal projective varieties.
 The exceptional locus of $\Phi$ is the closed subscheme of $X$ on which $\Phi$ is not an isomorphism. Irreducible components $T$ of the exceptional locus may have any dimension $1\le \text{dim }T\le\text{dim }X-1$. If $Y$ is nonsingular all components have dimension $\text{dim }X-1$.

Now suppose that $\Phi:X\rightarrow Y$ is a birational map of normal projective varieties.
There is a largest open set $U$ of $X$ such that $\Phi\mid U$ is a morphism. The complement, $F=X-U$ is called the fundamental locus of $\Phi$.
 By Zariski's main theorem (Lemma V.5.1 and Theorem V.5.2 \cite{Ha}),
  the fundamental locus $F$ has codimension $\ge 2$ in $X$.

Suppose that $X$ is a  variety and $Y$ is a  subscheme.
The ideal sheaf of $Y$ in $X$ will be denoted by ${\mathcal I}_Y$.

Suppose that $X$ is a nonsingular variety, and $Y$ is a nonsingular subvariety.
 The blow up of $X$ with center $Y$ is the blow up $X_1\rightarrow X$ of
${\mathcal I}_Y$. $X_1$ is a nonsingular variety.

Suppose that $X_1\rightarrow X$ is the blow up of the nonsingular center $Y$,
$V$ is a valuation ring of the function field $K$ of $X$,  the
center of $V$ on $X$ is $p$ and the center of $V$ on $X_1$ is $p_1$.
Then ${\mathcal O}_{X,p}\rightarrow {\mathcal O}_{X_1,p_1}$ will be called a local
blow up of  ${\mathcal O}_{X,p}$ along $V$. The induced homomorphism
$\hat {\mathcal O}_{X,p}\rightarrow \hat{\mathcal O}_{X_1,p_1}$
will also be called a local blow up.

We will make use of the fact that any birational morphism $X\rightarrow Y$ of projective varieties
is the blow up of an ideal sheaf (Theorem 7.17 \cite{Ha}), and the ``universal property of blow ups'' (Proposition 7.14 \cite{Ha}).
Foundational material on blow ups, strict transforms, and other matters are covered in Chapter
4 of \cite{C2}.

Suppose that $X\rightarrow Y$ is a dominant rational map of varieties. If the extension of function fields ${\bf k}(Y)\rightarrow {\bf k}(X)$ is finite, we denote
$[X:Y]=[{\bf k}(X):{\bf k}(Y)]$.

Suppose that $X$ is a nonsingular variety. An effective  divisor
$D$ on $X$ is a simple normal crossings divisor on $X$ if for each
$p\in X$ there exist regular parameters $x_1\ldots,x_n$ in
${\mathcal O}_{X,p}$ and $a_1,\ldots,a_n\in{\bf N}$ such that  $D$
is the divisor $x_1^{a_1}\cdots x_n^{a_n}=0$ in a neighborhood of
$p$. We will abbreviate a simple normal crossings divisor as a SNC
divisor.  Suppose that $X$ is a nonsingular variety, $D$ is a SNC
divisor on $X$ and $Y$ is a subscheme of $X$. We will say that $Y$
is transversal to  $D$ if $Y$ is nonsingular, and for $p\in Y\cap D$,
there exist regular parameters $x_1\ldots,x_n$ in ${\mathcal
O}_{X,p}$, $r\le n$  and $a_1,\ldots,a_n\in{\bf N}$ such that  $D$
is the divisor $x_1^{a_1}\cdots x_n^{a_n}=0$ in a neighborhood of
$p$, and $V(x_1,\ldots,x_r)\subset \text{Spec}({\mathcal
O}_{X,p})$ is
 the germ of $Y$.

If $D_1, D_2$ are two Cartier divisors on a variety $X$, we will write $D_1 \sim D_2$  if $D_1$ is linearly equivalent to $D_2$.

\section{The $\nu$ and $\tau$ Invariants}\label{Section4}

Suppose that $V$ is a nonsingular variety over an algebraically
closed field ${\bf k}$ of arbitrary characteristic, and ${\mathcal I}$
is a nonzero ideal sheaf of $V$. Suppose that $q\in V$ is a (not
necessarily closed) point. If $Z$ is a subscheme of $V$ with ideal
sheaf ${\mathcal I}_Z$, define $\nu_q(Z)=\nu_q({\mathcal I}_Z)$
(with notation as in Section \ref{SecPre}). Foundational material on singularities
can be found in \cite{C2}.

For $t\in{\bf N}$, let
$$
\text{Sing}_t({\mathcal I})=\{p\in V\mid \nu_p({\mathcal I})\ge
t\},
$$
which is a Zariski closed subset of $V$ (Theorem A.19 \cite{C2}).
Let
$$
r=\text{max}\{t\mid \text{Sing}_t({\mathcal
I})\ne\emptyset\}=\text{max}\{\nu_p({\mathcal I})\mid p\in V\}.
$$

Suppose that $Y$ is a nonsingular, integral subvariety of $V$. Let
$\pi_1:V_1\rightarrow V$ be the blow up of $Y$ with exceptional
divisor $F$. Let $t=\nu_q({\mathcal I}_q)$, where $q$ is the
generic point of $Y$. Then the weak transform ${\mathcal I}_1$ of
${\mathcal I}$ on $V_1$ is defined by (cf page 65 \cite{C2})
$$
{\mathcal I}{\mathcal
O}_{V_1}={\mathcal O}_{V_1}(-tF){\mathcal I}_1.
$$
$\nu_q({\mathcal I})=t$ implies that ${\mathcal I}_1$ is an ideal sheaf on $V_1$.
If $Z$ is the subscheme of $V$ defined by ${\mathcal I}$,
 then the weak transform $Z_1$ of $Z$ is the subscheme of $V_1$ defined by ${\mathcal I}_1$.
If $Z$ is a divisor on $V$, then the weak transform of $Z$
is the strict transform of $Z$ (page 38, page 65 \cite{C2}).

Suppose that $p\in \text{Sing}_r({\mathcal I})$ is a closed point,
$x_1,\ldots,x_n$ are regular parameters in ${\mathcal O}_{V,p}$
and $f\in{\mathcal I}_p$. There is an
expansion
$$
f=\sum_{i_1+\cdots+i_n\ge r}a_{i_1,\ldots, i_n}x_1^{i_1}x_2^{i_2}\cdots x_n^{i_n}
$$
with $a_{i_1,i_2,\ldots,i_n}\in{\bf k}$ in $\hat{\mathcal
O}_{V,p}={\bf k}[[x_1,\ldots,x_n]]$. The $r$-leading form $L$ of $f$
(with respect to $x_1,\ldots,x_n$)
is defined to be
$$
L(x_1,x_2,\ldots,x_n)=\sum_{i_1+\cdots+i_n=r}a_{i_1,i_2,\ldots, i_n}x_1^{i_1}\cdots x_n^{i_n}.
$$
Define $\tau(p)=\tau_{\mathcal I}(p)$ to be the dimension of the
smallest linear subspace $T$ of the {\bf k}-subspace spanned by
$x_1,x_2,\ldots,x_n$ in ${\bf k}[[x_1,\ldots, x_n]]$ such that
$L\in{\bf k}[T]$ for all $f\in{\mathcal I}_p$.
We will call the subvariety $M_p=V(T)$ of $\text{Spec}({\mathcal
O}_{V,p})$ an {\it approximate manifold} to ${\mathcal I}$ at $p$. $T$
depends on our choice of regular parameters at $p$. However, the tangent space to $M_p$ does not depend on $T$, and
the dimension $n-\tau(p)$ of $M_p$ does not depend on $T$.    After making a choice of
$x_1,\ldots,x_n$, $T$ is uniquely determined.
 It follows that if
$Y\subset\text{Sing}_r({\mathcal I})$ is a nonsingular subvariety
with $p\in Y$, then there exists an approximate manifold $M_p$ to
${\mathcal I}$ at $p$ such that the germ of $Y$ at $p$ is
contained in $M_p$.

\begin{Lemma}\label{Lemma107}  Suppose that $Y\subset\text{Sing}_r({\mathcal I})$
is a nonsingular subvariety of $V$, $\pi_1:V_1\rightarrow V$ is the blow up of
$Y$, ${\mathcal I}_1$ is the weak transform of ${\mathcal I}$
 on $V_1$, $p\in Y$,  $M_p$ is an approximate manifold to ${\mathcal I}$ at $p$ containing the germ of $Y$
at $p$, and $q\in\pi_1^{-1}(p)$. Then
\begin{enumerate}
\item[1.] $\nu_q({\mathcal I}_1)\le r$. \item[2.] $\nu_q({\mathcal
I}_1)=r$ implies $q$ is on the strict transform $M_p'$ of $M_p$ and
$\tau(p)\le\tau(q)$.
\item[3.] Suppose that $\nu_q({\mathcal I}_1)=r$ and $\tau(p)=\tau(q)$. Then there exists
an approximate manifold $M_q$ to ${\mathcal I}$ at $q$ such that
$M_q\cap \pi_1^{-1}(p)=M_p'\cap\pi_1^{-1}(p)$ where $M_p'$ is the strict transform of $M_p$ on $V_1$.
\end{enumerate}
\end{Lemma}
\begin{proof} This follows from a local calculation, as in  Lemma 6.4 and Lemma 7.5 \cite{C2}. \end{proof}

With the hypotheses of Lemma \ref{Lemma107}, we immediately
conclude that $\tau(p)\le n-\mbox{dim }Y$ and $\tau(p)= n-\mbox{dim }Y$ implies
$\nu_q({\mathcal I}_1)<r$.

We deduce from Lemma \ref{Lemma107} the following important corollary.

\begin{Lemma}\label{LemmaZL} Suppose that $V$ is a nonsingular 3-fold, $Y\subset\mbox{Sing}_r({\mathcal I})$ is a nonsingular subvariety of $V$, $\pi_1:V_1\rightarrow V$ is the blow up of $Y$, and $\mathcal{I}_1$ is the
weak transform of $\mathcal {I}$ on $V_1$. Let $F$ be the reduced exceptional divisor of $\pi_1$.
Suppose that $F\cap\mbox{Sing}_r({\mathcal I}_1)\ne\emptyset$. Then
\begin{enumerate}
\item[1.] If $Y$ is a point, then $F\cap \mbox{Sing}_r({\mathcal I}_1)$
is either an isolated point, or an irreducible nonsingular curve.
\item[2.] If $Y$ is a curve, then $F\cap\mbox{Sing}_r({\mathcal I}_1)$ 
is either a finite union of points, each in a distinct fiber over a point of $Y$, or $\mbox{Sing}_r({\mathcal I}_1)$ is an irreducible nonsingular curve 
which is a section over $Y$.
\end{enumerate}
\end{Lemma}

A nonsingular hypersurface $D_p\subset \mbox{Spec}({\mathcal O}_{V,p})$ is called an {\it approximate hypersurface} to ${\mathcal I}$ at $p$ if $D_p$ contains an
approximate manifold $M_p$ to ${\mathcal I}$ at $p$.
If $D_p$ is an approximate hypersurface to ${\mathcal I}$ at $p$, and $Y\subset \mbox{Sing}_r({\mathcal I})$ is a nonsingular subvariety containing $p$, then there exists an approximate hypersurface $D_p'$
to ${\mathcal I}$ at $p$ such that $Y\subset D_p'$ and $D_p$ and $D_p'$ have the same tangent space.
We deduce the following lemma from Lemma \ref{Lemma107}.

\begin{Lemma}\label{Lemma*107}  Suppose that $Y\subset\text{Sing}_r({\mathcal I})$
is a nonsingular subvariety of $V$, $\pi_1:V_1\rightarrow V$ is the blow up of
$Y$, ${\mathcal I}_1$ is the weak transform of ${\mathcal I}$
 on $V_1$, $p\in Y$,  $D_p$ is an approximate hypersurface  to ${\mathcal I}$ at $p$ containing the germ of $Y$
at $p$, and $q\in\pi_1^{-1}(p)$. Then
\begin{enumerate}
\item[1.] $\nu_q({\mathcal I}_1)\le r$ \item[2.] $\nu_q({\mathcal
I}_1)=r$ implies $q$ is on the strict transform $D_p'$ of $D_p$.
\item[3.] Suppose that $\nu_q({\mathcal I}_1)=r$. Then there exists
an approximate hypersurface $D_q$ to ${\mathcal I}$ at $q$ such that
$D_q\cap \pi_1^{-1}(p)=D_p'\cap\pi_1^{-1}(p)$ where $D_p'$ is the strict transform of $D_p$ on $V_1$.
\end{enumerate}
\end{Lemma}

\begin{Definition}\label{Def100} A resolution datum ${\mathcal R}=(E^+,E^-,{\mathcal I},V)$ is a 4-tuple where $E^+$ and $E^-$ are effective divisors on a nonsingular variety $V$
over an algebraically closed field $\bf k$, such that $E=E^++ E^-$ is a SNC divisor, $E^+$ and $E^-$ have no common components, and $\mathcal I$ is an
ideal sheaf on $V$. 
\end{Definition}

 For $t\ge 1$, let
$$
\text{Sing}_t({\mathcal R})=\text{Sing}_t({\mathcal I})=\{p\in V\mid \nu_p({\mathcal I})\ge t\}.
$$
$\text{Sing}_t({\mathcal R})$ is a proper Zariski closed subset of
$V$ for $t\ge 1$.

Let
$$
r=\nu({\mathcal R})=\nu({\mathcal I})=\text{max}\{\nu_p({\mathcal I})\mid p\in
V \}.
$$

For $p\in V$, let
$$
\eta(p)=\text{ the number of components of $E^-$ containing $p$}.
$$
We have $0\le \eta(p)\le \mbox{dim }V$.

Let $\tau(p)=\tau_{{\mathcal I}}(p)$ for $p\in V$.

Transversality is defined in the discussion of birational geometry in Section
\ref{SecPre} on preliminaries.

\begin{Definition}\label{Def101} A permissible transform of ${\mathcal R}$ is the blow up $\pi_1:V_1\rightarrow V$
of a   nonsingular subvariety $Y\subset \text{Sing}_{r}({\mathcal R})$
such that $Y$ is transversal to  $E=E^++E^-$.
\end{Definition}

 Let $F$ be the exceptional divisor of a permissible transform
$\pi:V_1\rightarrow V$ of ${\mathcal R}$. Then
$\pi_1^*(E^++E^-)+ F$ is a SNC divisor. We define the
transform ${\mathcal R}_1$ of ${\mathcal R}$ on $V_1$ to be ${\mathcal
R}_1=(E_1^+,E_1^-,{\mathcal I}_1,V_1)$ where ${\mathcal I}_1$ is the weak transform of
${\mathcal I}$, and
$$
E^+_1=\pi_1^*(E^+)+F, E^-_1=\text{ the strict transform of $E^-$
if }\nu({\mathcal I}_1)=\nu({\mathcal
R}),
$$
$$
E^+_1=\emptyset, E_1^-=\pi_1^*(E^++E^-)+F\mbox{ if }\nu({\mathcal I}_1)<\nu({\mathcal
R}).
$$
We remark that if $\pi_1(p)=q$ and 
$\nu_q({\mathcal I}_1)\le \nu_p({\mathcal I})=r$, then 
$\eta(q)\le\eta(p)$.

In the next several sections we prove the following theorem, which is our
main resolution theorem for ideals on a nonsingular 3-fold.

\begin{Theorem}\label{Theorem6} Suppose that $V$ is a nonsingular 3-fold over an algebraically closed field ${\bf k}$, and ${\mathcal R}=(\emptyset,E,{\mathcal I},V)$ is a resolution datum. Let $r=\nu({\mathcal R})\ge 1$. Then there exists a sequence of permissible transforms
$\pi:V_1\rightarrow V$ of ${\mathcal R}$ such that $\text{Sing}_r({\mathcal
R}_1)=\emptyset$, where ${\mathcal R}_1$ is the transform of
${\mathcal R}$ by $\pi$.
\end{Theorem}

Theorem \ref{Theorem6} implies the embedded resolution theorems
Theorem \ref{Theorem1} and Theorem \ref{Theorem23} stated in the introduction.
Using classical techniques of resolution, and  the resolution algorithm of Levi \cite{Le} and Zariski \cite{Z4}, we reduce in Sections \ref{SecERII} and \ref{SecERIII}  to a local analysis
of a sequence of points under which the order of the weak transform
of the ideal does not drop. We show that this sequence must be finite in
Sections \ref{Sect2} and \ref{Sect1}.
The essential ingredient in this local analysis is an extension of  Hironaka's termination invariant  for surfaces \cite{H2}. In Section \ref{SecERSI}, we deduce Theorem \ref{Theorem1} and \ref{Theorem23} from Theorem \ref{Theorem6}.

In the course of the proof of Theorem \ref{Theorem6}, we will construct sequences of
permissible transforms
$$
V_n\rightarrow V_{n-1}\rightarrow\cdots\rightarrow V_1\rightarrow V.
$$
To simplify notation, we will use the convention that ${\mathcal R}_n=(E_n^+,E_n^-,{\mathcal I}_n,V_n)$ where
$E_n=E_n^++E_n^-$ is the transform of ${\mathcal R}$ on $V_n$.

We will also construct  sequences of permissible transforms involving $W_i$,
$$
W_m\rightarrow W_{m-1}\rightarrow \cdots\rightarrow W_1\rightarrow V.
$$
We will use the notational convention that 
${\mathcal R}_m'=((E_m')^+,(E_m')^-,{\mathcal I}_m',W_m)$ where
$E_m'=(E_m')^++(E_m')^-$ is the transform of ${\mathcal R}$ on $W_m$.

\section{Embedded Resolution on a Non Singular Surface}\label{SecERS}

In this section, we prove resolution for an ideal sheaf on a nonsingular surface.
We define a resolution datum on a nonsingular variety over an  arbitrary field
in the same way as for algebraically closed fields.

\begin{Theorem}\label{Lemma26}
Suppose that $T$ is a nonsingular surface over a field $L$, and
${\mathcal R}=(\emptyset,E,{\mathcal I},T)$ is a resolution datum. Let $r=\nu({\mathcal R})\ge 1$. Then there exists a
sequence of permissible transforms $\pi:T_1\rightarrow T$ of ${\mathcal R}$ such that $\mbox{Sing}_r({\mathcal R}_1)=\emptyset$, where ${\mathcal R}_1$ is the transform of ${\mathcal R}$ by $\pi$.
\end{Theorem}

\begin{proof} After a few blow ups of points on the strict transforms of components of $E$, we construct $T_1\rightarrow T$ such that $\eta(p)=0$ for all $p\in\mbox{Sing}_r({\mathcal R}_1)$. We may thus assume that this is true on $T$.

Suppose that $C\subset \text{Sing}_r({\mathcal I})$ is a curve.
 Then $C$ is nonsingular, and there exists a neighborhood $U$ of $C$ in $T$ such that
  ${\mathcal I}\mid{\mathcal O}_U={\mathcal I}_C^r\mid U$. Let $\pi_1:T_1\rightarrow T$ be the blow up of $C$,
   ${\mathcal I}_1$ be the weak transform of ${\mathcal I}$ on $T_1$. Then ${\mathcal I}_1\mid\pi_1^{-1}(U)={\mathcal O}_{\pi_1^{-1}(U)}$.
    Thus $\text{Sing}_r({\mathcal I}_1)\cap \pi^{-1}(U)=\emptyset$.

We  thus reduce to the case where $\text{Sing}_r({\mathcal I})$ is
a finite set of points.

Now the proof is an extension of the proof of Theorem 3.15 \cite{C2}.
We outline the proof here, and refer to \cite{C2} for   technical details. 
  We construct a sequence of projective morphisms
\begin{equation}\label{eq201}
\cdots\rightarrow
T_n\stackrel{\pi_n}{\rightarrow}\cdots\rightarrow
T_1\stackrel{\pi_1}{\rightarrow} T_0=T
\end{equation}
where each $\pi_{n+1}:T_{n+1}\rightarrow T_n$ is the blow up of
all points in $\text{Sing}_r({\mathcal I}_n)$, where ${\mathcal
I}_n$ is the weak transform of ${\mathcal I}_{n-1}$ on $T_n$. We
must show that this sequence is finite.

Suppose that (\ref{eq201}) has infinite length. We can then find
closed points $p_n\in T_n$ such that $\pi_n(p_n)=p_{n-1}$, and
$p_n\in\text{Sing}_r({\mathcal I}_n)$ for all $n$. We then have an
infinite sequence of local rings
$$
R=R_0\rightarrow R_1\rightarrow\cdots\rightarrow
R_n\rightarrow\cdots
$$
where $R_n=\hat{\mathcal O}_{T_n,p_n}$. Let $L_1$ be the
residue field of $R$. By Lemma 3.14 \cite{C2}, the residue field of
each of these local rings is equal to $L_1$.

We will define $\delta_{p_i}\in\frac{1}{r!}{\bf N}$ so that
$\delta_{p_i}=\delta_{p_{i-1}}-1$ for all $i$, leading to a
contradiction to the assumed infinite length of (\ref{eq201}).

Suppose that $x,y$ are regular parameters in ${\mathcal O}_{T_0,p_0}$.
Identigy $L_1$ with a coefficient field of $R_0$.  For $f\in
I={\mathcal I}_{p_0}R_0$, we have an expansion
$$
f=\sum_{i+j\ge r} a_{ij}x^iy^j
$$
with $a_{ij}\in L_1$. We will call $x,y$ good parameters
for $I$ if there exists $f\in I$ such that $a_{0r}\ne 0$
($\text{ord}(f(0,y))=r$). Suppose that $x,y$ are good parameters
for $I$. For $f\in I$, we  define
$$
\delta(f;x,y)=\left\{\begin{array}{ll}
\text{min}\left\{\frac{i}{r-j}\mid j<r\text{ and
}a_{ij}\ne 0\right\}&\mbox{ if $y^r\not\,\mid f$},\\
\infty&\mbox{ if $y^r\mid f$}
\end{array}\right.
$$
and
$$
\delta(I;x,y)=\text{min}\{\delta(f;x,y)\mid f\in I\}.
$$
Now define
$$
\delta_{p_0}=\text{sup}\{\delta(I;x,y_0)\mid
y_0=y-\sum_{i=1}^nb_ix^i \text{ with $n\in{\bf N}$ and
$b_i\in L_1$}\}.
$$
If $\delta_{p_0}=\infty$, then there exists a series
$g=y-\sum_{i=1}^{\infty}b_ix^i$ such that $(g^r)=I$, a contradiction to the assumption that $p$ is isolated in $\mbox{Sing}_r({\mathcal I})$. Thus $\delta_{p_0}\in\frac{1}{r!}{\bold N}$. Let $y_0\in{\mathcal O}_{T_0,p_0}$
be such that $\delta_{p_0}=\delta(I;x,y_0)$. Then we check that $x,y_1$ are regular parameters in $R_1$, where $y_0=xy_1$, that $x,y_1$ are good parameters for $I_1=\widehat{({\mathcal I}_1)_{p_1}}$, and $\delta_{p_1}=\delta_{p_0}-1$.
We thus see that $\delta_{p_i}<0$ for $i\gg 0$, a contradiction to our assumption that (\ref{eq201}) has infinite length.

\end{proof}

\begin{Remark}\label{RemarkL5} Theorem \ref{Lemma26} and its proof are applicable when $T=\mbox{Spec}(A)$, where $A$ is an excellent 2 dimensional
regular local ring containing a field.
\end{Remark}

\section{Reduction of Embedded Resolution to Resolution}\label{SecERII}

In this section, we suppose that  $V$ is a nonsingular 3-fold (a 3-dimensional variety over
 an algebraically closed field ${\bf k}$) and  ${\mathcal R}$ is a resolution datum on  $V$
with $r=\nu({\mathcal R})\ge 1$. We will transform
${\mathcal R}$ into a nice stable form (Theorem \ref{Theorem4}), from which we will obtain our main embedded resolution
theorem (Theorem \ref{Theorem6}) in  Sections \ref{SecERIII} - \ref{Sect1}.

\begin{Lemma}\label{Theorem2} Suppose that $p\in \text{Sing}_{r}({\mathcal R})$,  and an approximate hypersurface $D_p$ of ${\mathcal I}$ at $p$  is transversal to $E^+$ and is not a component of $E^+$.
Suppose that $\pi_1:V_1\rightarrow V$ is a permissible transform of a nonsingular subvariety $Y$ of $V$ containing $p$.
Let ${\mathcal R}_1$ be the  transform of ${\mathcal R}$ on $V_1$.
  Suppose that $p_1\in\pi_1^{-1}(p)\cap\text{Sing}_{r}({\mathcal R}_1)$. Let $D_{p_1}$
   be an approximate hypersurface of ${\mathcal I}_1$ at $p_1$. Then $D_{p_1}$ is transversal to  $E_1^+$, and is not a component of $E_1^+$.
\end{Lemma}

\begin{proof} The proof is immediate from 3 of Lemma \ref{Lemma*107}.
\end{proof}

\begin{Theorem}\label{Theorem4} Let ${\mathcal R}=(\emptyset,E,{\mathcal I},V)$. Then there exists a sequence of permissible transforms
$\pi:V_1\rightarrow V$   such that if ${\mathcal R}_1$ is the
transform of ${\mathcal R}$ on $V_1$, then
\begin{enumerate}
\item[1.] $\text{Sing}_{r}({\mathcal R}_1)\cap E_1^-=\emptyset$, so that $\eta(p)=0$ for $p\in\text{Sing}_{r}({\mathcal R}_1)$.
\item[2.] $\text{Sing}_{r}({\mathcal R}_1)\subset E_1^+$ and $\mbox{dim Sing}_r({\mathcal R}_1)\le 1$.
\item[3.] All irreducible curves $C\subset \mbox{Sing}_r({\mathcal R}_1)$ 
are nonsingular.
\item[4.]
If $p\in \text{Sing}_{r}({\mathcal R}_1)$,
then there exists an approximate hypersurface $D_p$ of ${\mathcal I}_1$ at $p$ which is transversal to 
$E_1^+$ and is not a component of $E_1^+$.
\end{enumerate}
\end{Theorem}

\begin{proof} We first establish 1. We begin by blowing up all triple points of $E$ (points $p$ with $\eta(p)=3$)
in $\mbox{Sing}_r({\mathcal R})$
to construct $V_1\rightarrow V$ such that   $\eta(p)\le 2$ for $p\in \mbox{Sing}_r({\mathcal R}_1)$. We then blow up all double curves of $E$
(curves $C$ such that  $\eta(p)=2$ for all $p\in C$)
in $\mbox{Sing}_r({\mathcal R}_1)$
to construct $V_2\rightarrow V_1$ such that   $\eta(p)\le 2$ for $p\in \mbox{Sing}_r({\mathcal R}_2)$, and $\eta(p)=2$
at only finitely many points.

Suppose that $p\in \mbox{Sing}_r({\mathcal R}_2)$ is an (isolated) point with $\eta(p)=2$. Then after a finite number of blow ups of points $W_1\rightarrow V_2$
 centered at the point on the strict transform of the double curve of $E_2$ through $p$
which dominates $p$, the transform ${\mathcal R}_1'$ of ${\mathcal R}$ on $W_1$ satisfies $\eta(p_1)<2$ at all points 
$p_1\in\mbox{Sing}_r({\mathcal R}_1')$ which dominate $p$.  To see this, 
let $x,y,z$ be regular parameters in
${\mathcal O}_{V_2,p}$ such that $xy=0$ is a local equation of $E_2^-$ at $p$. Then $x=y=0$ is a local equation of
the double curve $C$ through $p$, and after blowing up $n$ times $\psi:W_n\rightarrow V$ along the strict transform of $C$, the local ring 
on $W_n$ of the point
$p_n$ on the strict transform of $C$ which dominates $p$ has regular parameters 
$x=x_1z^n, y=y_1z^n, z$. Substituting into generators of $({\mathcal I}_2)_p$, we see that $\nu_{p_n}({\mathcal I}_n')<r$
for $n\gg 0$, since $C\not\subset \mbox{Sing}_r({\mathcal R}_2)$.

Repeating this for all $p\in \mbox{Sing}_r({\mathcal R}_2)$ with $\eta(p)=2$, we 
construct $V_3\rightarrow V_2$ such that  $\eta(p)\le 1$ for $p\in \mbox{Sing}_r({\mathcal R}_3)$.

If there is an irreducible  component $F$ of $E_3^{-}$ contained in $\mbox{Sing}_r({\mathcal R}_3)$, then we can blow up
$F$ by $V_4\rightarrow V_3$ to eliminate $F$ from $\mbox{Sing}_r({\mathcal R}_4)$. In this way we reduce to the assumption that 
no irreducible  components $F$ of $E_4^-$ are contained in $\mbox{Sing}_r({\mathcal R}_4)$.

Now suppose that $C$ is a curve in $V_4$ which is contained in $\mbox{Sing}_r({\mathcal R}_4)$, and there is an
irreducible  component
$F$ of $E_4^-$ containing $C$ (so that $\eta(p)=1$ for all $p\in C$). By blowing up points on the strict transform of
$C$ we can construct  $W_1\rightarrow V_4$ such that the strict transform $C_1$ of  $C$ on $W_1$ is nonsingular and intersects $E_1'$ transversely (by Theorem \ref{Lemma26}). Now we construct the
sequence of permissible transforms 
\begin{equation}\label{eqL1}
\cdots\rightarrow W_n\rightarrow W_{n-1}\rightarrow \cdots \rightarrow W_2\rightarrow W_1
\end{equation}
where each morphism $W_{i+1}\rightarrow W_i$ is the blow up of the section $C_i$ over $C_1$ on the strict transform $F_i$ of  $F$ on $W_i$
($F_i$ is isomorphic to $F_1$ by this morphsim), and we continue as long as $C_i$ is contained in
$\mbox{Sing}_r({\mathcal R}_i')$ (where ${\mathcal R}_i'$ is the transform of ${\mathcal R}$ on $W_i$).

Let $\zeta$ be the generic point of $C$ in $V_4$, and let $\zeta_i$ be the generic points of $C_i$ in $W_i$ for all $i$. Base change of (\ref{eqL1}) by the inclusion of $\mbox{Spec}({\mathcal O}_{V_4,\zeta})\cong \mbox{Spec}({\mathcal O}_{W_1,\zeta_1})$ into $W_1$ induces a sequence of birational 
projective morphisms of regular, 2 dimensional schemes. For all $i\ge 1$,
${\mathcal I}_{C_i}\otimes_{{\mathcal O}_{V_4}}{\mathcal O}_{V_4,\zeta}$
is the ideal sheaf of a closed point $\gamma_i$ in $W_i\times_{V_4}\mbox{Spec}({\mathcal O}_{V_4,\zeta})$ which dominates 
$\gamma_{i-1}$, with $\gamma_0=\mbox{Spec}({\mathcal O}_{V_4,\zeta}/
{\mathcal I}_{C,\zeta})$. We see that by
taking base change of the sequence (\ref{eqL1}) by localization at $\zeta$, we obtain a sequence of blow ups of closed points above $\text{Spec}({\mathcal O}_{V_4,\zeta})$. 
Now by Remark \ref{RemarkL5} to Theorem \ref{Lemma26}, we have 
that $\gamma_i\not\in \mbox{Sing}_r(W_i\times_{V_4}\mbox{Spec}({\mathcal O}_{V_4,\zeta}))$ for all $i\gg0$. Since 
${\mathcal O}_{W_i,\zeta_i}\cong 
\mathcal O_{W_i\times_{V_4}\mbox{Spec}({\mathcal O}_{V_4,\zeta}),\gamma_i}$,
we see that $C_i\not\in \mbox{Sing}_r({\mathcal R}_i')$ for all $i\gg0$.

Repeating this for all curves $C$ in $V_4$ which are contained in $\mbox{Sing}_r({\mathcal R}_4)$, and such that there is an
irreducible  component
$F$ of $E_4^-$ containing $C$, we construct $V_5\rightarrow V_4$ such that $\eta(p)\le 1$ for all $p\in \mbox{Sing}_r({\mathcal R}_5)$, and 
there are at most finitely many points  $p\in \mbox{Sing}_r({\mathcal R}_5)$ such that $\eta(p)=1$.

Suppose that $p \in \mbox{Sing}_r({\mathcal R}_5)$ and $\eta(p)=1$. Let $A=\widehat{{\mathcal O}_{V_5,p}}$, $I=\widehat{({\mathcal I}_5)_p}$.
Let $F$ be the irreducible component of $E_5^-$ containing $p$. Choose regular parameters $x,y,z$ in ${\mathcal O}_{V_5,p}$ such that $z=0$ is a local equation of $F$. 

For $f\in I$, we have an expansion
$$
f=\sum_{i=0}^{r-1}a_i(x,y)z^i+z^r\Lambda
$$
in $A={\bf k}[[x,y,z]]$.  Define $C(I)$ to be the ideal in $B=\widehat{{\mathcal O}_{F,p}}={\bf k}[[x,y]]$
generated by 
$$
\left\{a_i(x,y)^{\frac{r!}{r-i}}\mid f\in I\right\}.
$$
We easily verify that for an ideal $H\subset A$, $\nu_A(H)\ge r$ if and only if $\nu_B(C(H))\ge r!$.

By  Remark \ref{RemarkL5} to Theorem \ref{Lemma26}, applied to $\mbox{Spec}(B)$ and $J=C(I)$, there exists a sequence of
permissible transforms
\begin{equation}\label{eqL2}
T_m\rightarrow T_{m-1}\rightarrow \cdots\rightarrow T_0=\mbox{Spec}(B)
\end{equation}
such that the weak transform  $J_m$ of $J$ on $T_m$ satisfies $\nu(J_m)<r!$.
Since $\mbox{Sing}_{r!}(J)$ is the maximal ideal of $B$ (as $p$ is isolated in $F\cap \mbox{Sing}_r({\mathcal I})$
by our reduction) every curve blown up in (\ref{eqL2}) is the strict transform of an exceptional divisor of some
$T_{i+1}\rightarrow T_i$.  Thus there exists a sequence of permissible transforms
\begin{equation}\label{eqL3}
W_m\rightarrow W_{m-1}\rightarrow \cdots\rightarrow W_0=V_5
\end{equation}
such that each $W_{i+1}\rightarrow W_i$ is either the blow up of a point on the strict transform of $F$ or the blow
up of a nonsingular curve on the strict transform of $F$, and (\ref{eqL2}) is obtained from (\ref{eqL3})
by taking the strict transforms of $F$ by the successive morphisms, and
taking base change by the inclusion $\mbox{Spec}(B)\rightarrow F$. 
Let ${\mathcal R}_i'$ be the transform of ${\mathcal R}_4$ on $W_i$.

For $p_1$ in the strict transform of $F$ by the map $W_1\rightarrow V_5$ 
which is also in the fiber above $p$ (where ${\mathcal R}_1'$ is the transform of ${\mathcal R}$ on $W_1$), we verify that
$C(\widehat{({\mathcal I}_1')_{p_1}})$ is the completion of the localization at $p_1$ of the weak transform ${\mathcal J}_1$ of $J$
on $T_1$, and $p_1\in \mbox{Sing}_r({\mathcal R}_1')$ if and only if $p_1\in \mbox{Sing}_{r!}({\mathcal J}_1)$.
By induction, we see that this property persists throughout the sequences (\ref{eqL2}) and (\ref{eqL3}), so that
$\mbox{Sing}_r({\mathcal R}_m')$ is disjoint from the strict transform of $F$ above $p$.

Iterating this construction for all (the finitely many) $p\in \mbox{Sing}_r({\mathcal R}_5)$ with $\eta(p)=1$,
we construct $V_6\rightarrow V_5$ such that  the conclusion 1 of the theorem holds on $V_6$.

We now establish 2.

Suppose that $S$ is a 2 dimensional irreducible component of $\mbox{Sing}_r({\mathcal R}_6)$. Then $S$ is nonsingular and isolated 
in $\mbox{Sing}_r({\mathcal R}_6)$. For $p\in S$, $S$
is an approximate hypersurface $D_p$ for ${\mathcal I}_6$ at $p$.
$E_6^+$ is transversal to  $D_p$, by Lemma \ref{Theorem2}. Thus $S$ is transversal to
$E_6=E_6^++E_6^-$ at all $p\in S$. The blow up
$W_1\rightarrow V_6$ of $S$ is thus permissible. We have
$\mbox{Sing}_r({\mathcal R}_1')\cap\pi_2^{-1}(S)=\emptyset$ (where ${\mathcal R}_1'$ is the transform
of ${\mathcal R}_6$ on $W_1$).

Repeating this construction for all 2 dimensional components of 
$\mbox{Sing}_r({\mathcal R}_6)$, we construct $V_7\rightarrow V_6$ such that 
 $\mbox{dim }\mbox{Sing}_r({\mathcal R}_7)\le 1$.

Suppose that $C_1,\ldots,C_n$ are the Zariski closures in $V_7$ of the 
irreducible curves in $\mbox{Sing}_r({\mathcal R}_7)-E_7$.
There exists (by Corollary 4.4 \cite{C2}) a sequence of blow ups of points
$\psi_1:W_1\rightarrow V_7$ with transform
${\mathcal R}_1'$ of 
${\mathcal R}_7$ on $W_1$, such that the strict transforms
$C_1',\ldots, C_n'$ of $C_1,\ldots, C_n$ on $W_1$ are nonsingular and 
disjoint.  Suppose that $p\in W_1$ is such that 
$E_1'$ is not transversal to $C_1'$ (necessarily $p\not\in (E_1')^-$).
There exist regular parameters $x,y,z$ in
${\mathcal O}_{W_1,p}$ such that $x=y=0$ are local equations of
$C_1'$ at $p$. Let
$f(x,y,z)=0$ be a local equation of $E_1'=(E_1')^++(E_1')^-$ at $p$ in
$\hat{\mathcal O}_{W_1,p}$. We necessarily have that $f\not\in
(x,y)\hat{\mathcal O}_{W_1,p}$. Let $\psi_2:W_2\rightarrow W_1$ be
the blow up of $p$, and suppose that $p_1\in\pi_2^{-1}(p)$ is on
the strict transform of $C_1'$. Then there
exist regular parameters $x_1,y_1,z_1$ in ${\mathcal O}_{W_2,p_1}$
such that
$$
x=x_1z_1, y=y_1z_1, z=z_1.
$$
Since $f\not\in(x,y)$, we see by substituting
$$
x=x_nz_n^n, y=y_nz_n^n, z=z_n
$$
into $f(x,y,z)$, that after finitely many blow ups of points
$W_{n+1}\rightarrow W_1$, the strict transform of $C_1$ is
disjoint from the strict transform of $E_1'$.

Since $C_1'$ is nonsingular, the strict transform $C_1''$ of $C_1'$ on $W_{n+1}$ is transversal to  the exceptional divisor of $W_{n+1}\rightarrow W_1$. Thus, $C_1''$ is transversal to  $E_{n+1}'$. The blow up $W_{n+2}\rightarrow W_{n+1}$ of $C_1''$ is thus permissible.

We repeat this construction for $C_2,\ldots, C_n$ to construct $V_8\rightarrow V_7$  such that $\mbox{Sing}_r({\mathcal R}_8)-
E_8$ is a finite set of points.  Let $V_9\rightarrow V_8$ be the blow up
of these points. The conclusions 1 and 2 of the theorem  hold on $V_9$. The conclusions 3 and 4
of the  theorem now hold by
Lemma \ref{LemmaZL} and   Lemma \ref{Theorem2}.

\end{proof}

Now suppose that $S$ is a reduced effective divisor on $V$ such that $S$ 
is nonsingular, and suppose that ${\mathcal R}=(E^+,E^-,{\mathcal I}_S,V)$ is
a resolution datum on $V$, such that $E^++S$ is a SNC divisor. Let
$$
\mbox{NSNC}({\mathcal R})=\{p\in S\mid E^++E^-+S\mbox{ is not a SNC divisor at $p$}\}.
$$
$\mbox{NSNC}({\mathcal R})$ is a closed subset of $\mbox{Sing}_1({\mathcal R})=S$. 

We define a NC-permissible transform of ${\mathcal R}$ to be the blow up
$\pi_1:V_1\rightarrow V$ of a nonsingular subvariety $Y$ of $\mbox{NSNC}({\mathcal R})$, such that $Y$ is transversal to $E=E^++E^-$. We define the transform ${\mathcal R}_1$ of ${\mathcal R}$ with respect to a 
NC-permissible transform $\pi_1$ as for a permissible transform in Section 
\ref{Section4}.

Observe that if $\pi_1$ is a NC-permissible transform of ${\mathcal R}$, then 
${\mathcal I}_1={\mathcal I}_{S_1}$, where $S_1$ is the strict transform of $S$
by $\pi_1$. We have that $S_1$ is nonsingular, $E_1^++E_1^-$ is a SNC divisor, and $E_1^++S_1$ is a SNC divisor.

\begin{Lemma}\label{LemmaFL} Let ${\mathcal R}=(\emptyset,E,{\mathcal I}_S,V)$
be a resolution datum on $V$, where $S$ is a reduced effective divisor on $V$ such that $S$ is nonsingular. Then there exists a sequence of NC-permissible
transforms $\pi_1:V_1\rightarrow V$ such that if ${\mathcal R}_1$ is the transform
of ${\mathcal R}$ on $V_1$, then $\mbox{NSNC}({\mathcal R}_1)=\emptyset$.
\end{Lemma}

\begin{proof} A simplification of the proof of 1 of Theorem \ref{Theorem4},
with $\mbox{Sing}_r$ replaced with $\mbox{NSNC}$, shows that there exists a
sequence of NC-permissible transforms $\pi:V_1\rightarrow V$, such that
$\eta(p)=0$ for $p\in\mbox{NSNC}({\mathcal R}_1)$. Since $E_1^++S_1$ is a SNC
divisor, we thus have that $E_1+S_1$ is a SNC divisor, so that $\mbox{NSNC}({\mathcal R}_1)=\emptyset$.
\end{proof}

\section{Reduction to Local Resolution}\label{SecERIII}
In this section we reduce the proof of Theorem \ref{Theorem6}
to a local problem, which is proven in Sections \ref{Sect2} and \ref{Sect1}.

We may assume that
the conclusions  of Theorem \ref{Theorem4} hold (but we now have
${\mathcal R}=(E^+,E^-,{\mathcal I},V)$).

 Suppose that $Y\subset
\text{Sing}_r({\mathcal R})$ is an irreducible component. 
 By the conclusions of Theorem
\ref{Theorem4},  $Y$ is nonsingular of dimension $\le 1$, $Y\cap E^{-}=\emptyset$, $Y\subset E^+$
and there exists an approximate hypersurface $D_p$ to ${\mathcal I}$ at $p$ which is transversal to  $E^+$ such that $Y\subset D_p$. Thus $Y$ is transversal to $E$, and the blow up  $\pi_1:V_1\rightarrow V$ of $Y$ is a permissible transform of ${\mathcal R}$.
Furthermore, it follows from Lemmas \ref{LemmaZL} and  \ref{Theorem2} that the conclusions of
Theorem \ref{Theorem4} hold for the transform ${\mathcal R}_1$ of ${\mathcal R}$ on $V_1$. 

We may thus  construct a sequence of  permissible
transforms
\begin{equation}\label{In5}
\cdots\rightarrow V_n\rightarrow V_{n-1}\rightarrow \cdots\rightarrow V_1\rightarrow V
\end{equation}
by applying the following algorithm:
\vskip .2truein

\begin{enumerate}
\item[1.] If there exists an irreducible curve $C\subset \text{Sing}_r({\mathcal
R}_n)$  then perform the permissible transform $V_{n+1}\rightarrow V_n$ obtained by  blowing up  $C$.
\item[2.] If $\text{Sing}_r({\mathcal
R}_n)$ is a finite set of points, then blow up a point in $\text{Sing}_r({\mathcal
R}_n)$. 
\end{enumerate}
\vskip .2truein
As remarked in the paragraph above (\ref{In5}), the above blow ups are permissible, and the conclusions of Theorem \ref{Theorem4} 
hold for all transforms ${\mathcal R}_n$ of ${\mathcal R}$ along the sequence
(\ref{In5}). 

We must show that the sequence terminates after a finite number of steps
with $\text{Sing}_r({\mathcal
R}_n)=\emptyset$.

We first observe that there exists a $V_n$ in (\ref{In5}) such that $\text{Sing}_r({\mathcal
R}_n)$ is a finite set. Suppose otherwise. Then (\ref{In5}) is an infinite sequence,
constructed by only performing Step 1 of the algorithm.  By Lemma \ref{LemmaZL}, there exists
a curve $C$ in $\text{Sing}_r({\mathcal
R})$ such that an infinite number of sections over $C$ are blown up in (\ref{In5}). let $\zeta$ be the generic point of $C$.
As in the argument following (\ref{eqL1}) in the proof of Theorem \ref{Theorem4},
the sequence of blow ups of closed points over $\mbox{Spec}({\mathcal O}_{V,\zeta})$ obtained  by making the base change of (\ref{In5}) with the inclusion of $\mbox{Spec}({\mathcal O}_{V,\zeta})$ into $V$,
does not reduce the order of the weak transform of ${\mathcal I}_n$, a contradiction to
Remark \ref{RemarkL5} to Theorem \ref{Lemma26}. We may thus assume that
$\text{Sing}_r({\mathcal R})$ is a finite set. 

With the assumptions that $\mbox{Sing}_r({\mathcal I})$ is a finite set
of points, and that the conclusions of Theorem \ref{Theorem4} hold, we now
construct a sequence of permissible transforms

\begin{equation}\label{eqL4}
\cdots\rightarrow V_n\rightarrow V_{n-1}\rightarrow \cdots\rightarrow V_1\rightarrow V
\end{equation}
by applying the following algorithm:
\vskip .2truein
\begin{enumerate}
\item[1.] If there is a curve in $\mbox{Sing}_r({\mathcal I}_n)$, then blow
up the union of one dimensional components of $\mbox{Sing}_r({\mathcal I}_n)$.
\item[2.] Otherwise, blow up $\mbox{Sing}_r({\mathcal I}_n)$,  which is a finite union of points.
\end{enumerate}
 
By Lemma \ref{Lemma*107}, Lemma \ref{LemmaZL}, 4 of Theorem \ref{Theorem4},
and since we started with $\mbox{Sing}_r({\mathcal I})$ on $V$ being a finite set of points, under this algorithm $\mbox{Sing}_r({\mathcal I}_n)$ is always a disjoint union of points and nonsingular curves, which are transversal to  $E_n$.
The conclusions of Theorem \ref{Theorem4} hold on all $V_n$.

If the sequence (\ref{eqL4}) is of finite length, then we have achieved the conclusions of Theorem \ref{Theorem6}.

Suppose that (\ref{eqL4}) has infinite length, and  there does not
exist a sequence of points $q_n\in\mbox{Sing}_r({\mathcal I}_n)$ such that
$q_{n+1}$ maps to $q_n$ for all $n$. We will derive a contradiction.

Let $K={\bf k}(V)$ be the function field of $V$, and let $\mbox{ZR}(K)$ be the 
Zariski Riemann manifold of $K$ (see Section \ref{Intro3}). Let $\Phi_n:\mbox{ZR}(K)\rightarrow V_n$ be the natural projection. Let $U_n=V_n-\mbox{Sing}_r({\mathcal I}_n)$, an open subset of $V_n$.

Suppose that $A\in\mbox{ZR}(K)$. Let $q_n$ be the center of $A$ on $V_n$. By assumption, $q_n\not\in \mbox{Sing}_r({\mathcal I}_n)$ for large $n$. Thus $A\in\Phi_n^{-1}(U_n)$. It follows that 
$$
\cup_{n=0}^{\infty}\Phi_n^{-1}(U_n)=\mbox{ZR}(K).
$$
Since $\mbox{ZR}(K)$ is quasi compact (\cite{Z5} or \cite{ZS}), there exists an $n$ such that $\Phi_n^{-1}(U_n)=\mbox{ZR}(K)$, so that $\mbox{Sing}_r({\mathcal I}_n)=\emptyset$,
and (\ref{eqL4}) has finite length, a contradiction.

Assume that (\ref{eqL4}) is not finite.  We must then have an infinite sequence of points $q_n\in V_n$ such that 
$q_n\in\mbox{Sing}_r({\mathcal I}_n)$, and $q_n$ maps to $q_{n-1}$  for all $n$.

By Lemma \ref{Lemma107}, $\tau(q_n)\ge \tau(q_{n-1})$ for all $n$.
Thus, by induction on $\tau$, we may assume that $\tau(q_n)=\tau(q)$ for all $n$.

We will show that this is not possible, by associating to each point $q_n$ an
element $\Omega(q_n)$ in an ordered set such that $\Omega(q_{n+1})<\Omega(q_n)$ for all $n$. We will further show that $\Omega$ cannot decrease indefinitely.
It will then follow that there exists an $n$ such that $\mbox{Sing}_r(V_n)=\emptyset$, proving Theorem \ref{Theorem6}.

We have a
sequence
\begin{equation}\label{eq35}
R_0\rightarrow R_1\rightarrow \cdots
\end{equation}
of infinite length, where $R_i=\hat{\mathcal O}_{V_i,q_i}$ is the
completion of the local ring of $V_i$ at $q_i$, and
$\nu_{q_i}({\mathcal I}_i)=r$  for all $i$.
 Let $I_i=({\mathcal I}_i)_{q_i}R_i$.

Since ${\mathcal O}_{V_n,q_n}$ is excellent, we have $\mbox{Sing}_r({\mathcal I}_n)\cap\mbox{Spec}(R_n)=\mbox{Sing}_r(I_n)$ for all $n$.  

We may extend our definition of $\tau$ to an ideal in a power series ring.
We have $\tau(q_n)=\tau(I_n)$ for all $n$.

We will consider separately the cases with different values of $\tau$
in the following Sections \ref{Sect2} and \ref{Sect1}, to derive a contradiction, showing that (\ref{eq35}) must have finite length
in all cases.

The case where $\tau(q)=3$ is immediate from  Lemma
\ref{Lemma107}.

\section{Reduction When $\tau(q)=2$}\label{Sect2}

The proof in this case is an extension of the case of an ideal sheaf on a nonsingular surface,  given in Section \ref{SecERS}. The proof is based on Hironaka's termination argument in \cite{H2}.

Suppose that $T$ is a power series ring in 3 variables over an algebraically closed field ${\bf k}$, and
$J\subset T$ is an ideal. Suppose that $x,y,z$ are regular parameters in $T$, and $r$ is a positive integer.

For $g\in T$, we have an expansion
\begin{equation}\label{La1}
g=\sum b_{ijk}x^iy^jz^k
\end{equation}
with $b_{ijk}\in{\bf k}$.
Define 
$$
\gamma(g;x,y,z)=\text{min}\{\frac{k}{r-(i+j)}\mid b_{ijk}\ne 0\text{ and }i+j<r\}
\in \frac{1}{r!}{\bold N}\cup \{\infty\}.
$$

Define
$$
\gamma(J;x,y,z)=\text{min}\{\gamma(g;x,y,z)\mid g\in J\}.
$$

We have  $\gamma(J;x,y,z)=\infty$ if and only if $J\subset (x,y)^r$.
Further,
\begin{enumerate}
\item[1.] $\nu_T(J)=r$ if and only if $\gamma(J;x,y,z)\ge 1$ and
\item[2.] $V(x,y)$ is contained in an approximate manifold of $J$ if and only if
$$
\gamma(J;x,y,z)>1.
$$
\end{enumerate}

Let $\gamma=\gamma(J;x,y,z)$.
For $g\in J$, define

\begin{equation}\label{eqXS1}
[g]_{xyz}=\sum_{(i+j)\gamma+k=r\gamma}b_{ijk}x^iy^jz^k.
\end{equation}
There is an expansion
\begin{equation}\label{eqXS2}
g=[g]_{xyz}+\sum_{(i+j)\gamma+k>r\gamma}b_{ijk}x^iy^jz^k.
\end{equation}

Define
$$
A_{\gamma}=\left\{(i,j)\mid \frac{k}{r-(i+j)}=\gamma, i+j<r,\mbox{ and }  b_{ijk}\ne 0\mbox{ for some }
g\in J\right\}.
$$
Then for $g\in J$,
\begin{equation}\label{eqXS3}
[g]_{xyz}=L_g(x,y)+\sum_{(i,j)\in A_{\gamma}}b_{i,j,\gamma(r-i-j)}x^iy^jz^{\gamma(r-i-j)},
\end{equation}
where $L_g(x,y)=\sum_{i+j=r}b_{ij0}x^iy^j$.

Regular parameters $x,y,z$ in $T$ will be called
good parameters for $J$ if $\nu_{\overline T}(J\overline T)=r$ and $\tau(J\overline T)= 2$,
where $\overline T=T/zT$.

If $\nu_T(J)=r$ and $\tau(J)=2$, then regular parameters $x,y,z$ in $T$ such 
that $V(x,y)$ is an approximate manifold of $J$ are good parameters for
$J$.

\begin{Definition} 
Suppose that $x,y,z$ are good parameters for $J$ and  
$\nu_T(J)=r$.
Then $J$ is solvable with respect to $x,y,z$ if $\gamma=\gamma(J;x,y,z)\in{\bold N}$ and
there exist $\alpha,\beta\in {\bf k}$
such that
$$
[g]_{xyz}=L_g(x-\alpha z^{\gamma},y-\beta z^{\gamma}).
$$
for all $g\in J$ with $\nu_{T}(g)=r$.
\end{Definition}

\begin{Lemma}\label{LemmaXS11} Suppose that 
$x,y,z$ are good parameters for $J$, $\nu_T(J)=r$ and $\mbox{Sing}_r(J)$ is the maximal
ideal of $T$. Then
there exists a change of variables 
$$
x_1= x-\sum_{i=1}^n\alpha_i z^i,
y_1=y-\sum_{i=1}^n\beta_iz^i
$$
such that $J$ is not solvable with respect to $x_1,y_1,z$.
\end{Lemma}

\begin{proof}If there doesn't exist such a change of variables then 
we construct series
$$
x_{\infty}=x-\sum_{i=1}^{\infty}\alpha_iz^{i},
y_{\infty}=y-\sum_{i=1}^{\infty}\beta_iz^{i}
$$
such that 
$\gamma(J;x_{\infty},y_{\infty},z)=\infty$,
and thus
$$
J\subset (x-\sum_{i=1}^{\infty}\alpha_iz^i,y-\sum_{i-1}^{\infty}\beta_iz^i)^r\subset T,
$$
a contradiction to our assumption that $\mbox{Sing}_r(J)$ is the maximal ideal of $T$.
\end{proof}

Observe that if $\nu_T(J)=r$, $\tau(J)=2$, 
$x,y,z$ are good parameters for $J$ and $J$ is not solvable for $x,y,z$, then
$\gamma(J;x,y,z)>1$.

\begin{Lemma}\label{LemmaRSp9} Suppose that 
$\nu_T(J)=r$, $\tau(J)=2$, $\mbox{Sing}_r(J)$ is the maximal ideal of $T$,
$x,y,z$ are good parameters for $J$ and $J$ is not solvable with respect to  $x,y,z$.

Suppose that $T_1$ is the completion of a local ring of the blow up of the
maximal ideal of $T$ such that the weak transform $J_1$ of $J$ in $T_1$ satisfies $\nu_{T_1}(J_1)=r$ and $\tau(J_1)=2$. Then
\begin{enumerate}
\item[1.] $T_1$ has regular parameters $x_1,y_1,z_1$ defined by 
\begin{equation}\label{eqXS4}
x=x_1z_1, y=y_1z_1, z=z_1.
\end{equation}
\item[2.] $x_1,y_1,z_1$ are good parameters for $J_1$ and $V(x_1,y_1)$ is an 
approximate manifold of $J_1$.
\item[3.] $\gamma(J_1;x_1,y_1,z_1)=\gamma(J;x,y,z)-1$.
\item[4.] $J_1$ is not solvable with respect to $x_1,y_1,z_1$.
\end{enumerate}
\end{Lemma}

\begin{proof}
Assertion 1 follows from Lemma \ref{Lemma107}, as $V(x,y)$ is an approximate manifold for $J$. Assertions 3 and 4 follow from substitution of (\ref{eqXS4}) in (\ref{La1}), (\ref{eqXS1}), (\ref{eqXS2}) and (\ref{eqXS3}).

Since we are assuming that $\tau(J_1)=2$,
we  have that there exist $\alpha,\beta\in{\bf k}$ such that an
approximate manifold of $J_1$ has the form
\begin{equation}\label{eqL33}
V(x_1+\alpha z_1,y_1+\beta z_1),
\end{equation}
and $x_1,y_1,z_1$ are good parameters for $J_1$.

Since we are assuming that $\nu_{T_1}(J_1)=r$, we have $\gamma(J_1;x_1,y_1,z_1)\ge 1$. If $\gamma(J_1;x_1,y_1,z_1)>1$, then $V(x_1,y_1)$ is an approximate manifold of $I_1$.

Suppose that $\gamma(J_1;x_1,y_1,z_1)=1$. For $g\in J$ with $\nu_T(g)=r$
and  $g_1=\frac{g}{x_1^r}$, we have $\nu_{T_1}(g_1)=r$, and 
$$
[g_1]_{x_1y_1z_1}=\frac{1}{x_1^r}[g]_{xyz}.
$$
Since $\gamma(J_1;x_1,y_1,z_1)=1$, $[g_1]_{x_1y_1z_1}$ is the $r$-leading form
of $g_1$ with respect to $x_1, y_1, z_1$. By (\ref{eqL33}), there exists a form $\Psi_g(u,v)$, which depends on $g$, such that
$$
[g_1]_{x_1y_1z_1}=\Psi_g(x_1+\alpha z_1,y_1+\beta z_1)=L_g(x_1,y_1)+z_1\Omega
$$
where $\Omega\in T_1$ and $L_g(x,y)$ is the $r$-leading form of $g$ with respect to $x, y, z$.  Setting $z_1=0$, we have $L_g(x_1,y_1)=\Psi_g(x_1,y_1)$, so that
 $$
[g_1]_{x_1y_1z_1}=L_g(x_1+\alpha z_1,y_1+\beta z_1).
$$
We conclude that $J_1$ is solvable with respect to $x_1,y_1,z_1$, a contradiction to the
assumption that $\gamma(J_1;x_1,y_1,z_1)=1$.
\end{proof}

We now prove that (\ref{eq35}) cannot have infinite length.
 
Since $\tau(q_i)=2$ for all $i$, 
the assumption that (\ref{eq35}) is infinite and Lemma \ref{Lemma107} 
imply that for all $i$, $\mbox{Sing}_r(I_i)$ is the maximal ideal of $R_i$ and $R_{i+1}$ is the completion of a local ring of a closed point of the
blow up of the maximal ideal of $\text{Spec}(R_i)$.

Moreover,  $\tau(q)=2$ implies  that there exist regular parameters $(x,y,z)$ in $R_0$ such that $V(y,z)$ is an approximate manifold of $I_0$, and thus $x,y,z$ are good 
parameters for $I_0$.

By Lemma \ref{LemmaXS11}, there exists a change of variables in $R_0$ so
that we find good parameters $x,y,z$ for $I_0$ such that $I_0$ is not solvable with respect to $x,y,z$. Set $\Omega(q_0)=\gamma(I_0;x,y,z)$.

By Lemma \ref{LemmaRSp9}, we can inductively define positive rational numbers $\Omega(q_n)$  such that $\Omega(q_{n+1})=\Omega(q_n)-1$ for all $n$, giving
a contradiction if (\ref{eq35}) has infinite length.

\section{Reduction When $\bf \tau(q)=1$.}\label{Sect1}

This is the most interesting and difficult case. The proof is based on
Hironaka's termination argument in \cite{H2}.

\subsection{Definition of $\Omega$}

In this subsection we define $\Omega$ for an ideal $I$ in a power series ring
$T$ of three variables, over an algebraically closed field ${\bf k}$.

Suppose that $x,y,z$ are regular parameters in $T$, and $r$ is a positive integer.

For $g\in T$, we have an expansion
  $$
g=\sum b_{ijk}x^iy^jz^k
 $$
with $b_{ijk}=b_{ijk}(g)\in{\bf k}$. We define
$$
\Delta=\Delta(I;x,y,z)=\left\{\left(\frac{i}{r-k},\frac{j}{r-k}\right)\in{\bf
Q}^2\mid  k<r \text{ and }b_{ijk}\ne 0 \text{ for some }g\in I\right\}.
$$

Let $|\Delta(I;x,y,z)|$ be the smallest convex set in
${\bf R}^2$ such that $\Delta\subset |\Delta|$, and $(a,b)\in
|\Delta|$ implies $(a+c,b+d)\in |\Delta|$ for all $c,d\ge 0$.
For $\gamma\in {\bf R}$, let $S(\gamma)$ be the line through $(\gamma,0)$ with slope -1.

Suppose that $|\Delta|\ne\emptyset$.
We define $\alpha_{xyz}(I)$ to be the smallest $a$ appearing in any $(a,b)\in|\Delta|$, $\beta_{xyz}(I)$ to be the smallest $b$ such that $(\alpha_{xyz}(I),b)\in|\Delta|$.
Let $\gamma_{xyz}(I)$ be the smallest number $\gamma$ such that $S(\gamma)\cap|\Delta|\ne\emptyset$
and let $\delta_{xyz}(I)$ be such that $(\gamma_{xyz}(I)-\delta_{xyz}(I),\delta_{xyz}(I))$
is the lowest point on $S(\gamma_{xyz}(I))\cap|\Delta|$. $(\alpha_{xyz}(I),\beta_{xyz}(I))$
and $(\gamma_{xyz}(I)-\delta_{xyz}(I),\delta_{xyz}(I))$ are vertices of $|\Delta|$.
Define $\epsilon_{xyz}(I)$ to be the absolute value of the  largest slope of a line  through
$(\alpha_{xyz}(I),\beta_{xyz}(I))$ such that no points of $|\Delta|$ lie below it.

We now define
$$
\Omega(I;x,y,z) = (\beta_{xyz}(I),\frac{1}{e_{xyz}(I)},\alpha_{xyz}(I)),
$$
which is in the ordered set (by the Lex order)
$$
\left(\frac{1}{r!}{\bf N}\right)
\times \left({\bf Q}\cup{\infty}\right)\times \left(\frac{1}{r!}{\bf N}\right)
$$

We observe that:
\begin{enumerate}
\item[1.] The vertices of $|\Delta|$ are points of $\Delta$, which lie in the
lattice $\frac{1}{r!}{\bold Z}\times \frac{1}{r!}{\bold Z}$.
\item[2.] $\nu_R(I)<r$ holds if and only if $|\Delta|$ contains a point on $S(c)$
with $c<1$ which holds if and only if there is a vertex $(a,b)$ with  $a+b<1$.
\item[3.] $\alpha_{xyz}(I)<1$ if and only if $I\not\subset (x,z)^r$.
\item[4.] A vertex of $|\Delta|$ lies below the line $b=1$ if and only if 
$I\not\subset (y,z)^r$.
\end{enumerate}

\subsection{Well preparation}\label{Wp}

 $(x,y,z)$ will be called good parameters for $I$ if  
$\nu_T(I)=r$ and 
 $b_{00r}\ne 0$ for all $g\in I$  such that $\nu_T(g)=r$. Good parameters exist if $\tau(I)=1$.

Suppose  that $\nu_T(I)=r$ and 
 and $(x,y,z)$ are good parameters for $I$. 
Let $\Delta=\Delta(I;x,y,z)$ and
suppose that $(a,b)$ is a vertex of $|\Delta|$.
Define
$$
S_{(a,b)}=\left\{k\mid k<r, \mbox{ and }(\frac{i}{r-k},\frac{j}{r-k})=(a,b)\right\}.
$$
For $g\in I$,  define
$$
\{g\}^{ab}_{xyz}
=b_{00r}z^r+\sum_{k\in S_{(a,b)}}b_{a(r-k),b(r-k),k}x^{a(r-k)}y^{b(r-k)}z^k.
$$

We will say that  the vertex $(a,b)$ is not prepared
on $|\Delta|$ if
$a,b$ are integers and there exists $\eta\in{\bf k}$ such that 
$$
\{g\}^{ab}_{xyz}=b_{00r}(g)(z-\eta x^ay^b)^r
$$
for every $g\in I$ with $\nu_T(g)=r$.
 We will say that the vertex $(a,b)$ is prepared on 
$|\Delta|$ if such an $\eta$ does not exist.

If  the vertex $(a,b)$ is not prepared, then we can make an $(a,b)$ preparation,
which is the change of parameters
$$
z_1=z-\eta x^ay^b.
$$
If all  vertices $(a,b)$ of $|\Delta|$ are prepared, then
we say that $(I;x,y,z)$ is well prepared.

\begin{Lemma}\label{LemmaRSp14}
Suppose that $\nu_T(I)=r$,
$\tau(I)=1$, and
 $(x,y,z)$ are good parameters of $I$.
Then $z=0$ is an
approximate manifold of $I$ if and only if the vertices $(0,1)$
and $(1,0)$ are prepared on $|\Delta(I;x,y,z)|$.
\end{Lemma}

\begin{proof} Suppose that $(0,1)$ and $(1,0)$ are prepared on $|\Delta(I;x,y,z)|$.
Since $\tau(I)=1$, there is a linear form
$a x+b y+ c z$ with $ a, b, c\in{\bf k}$ such that $V( a x+b y+ cz)$
is an approximate manifold of $I$. For each $g\in I$ with $\nu_T(g)=r$,
there exists $0\ne d_g\in{\bf k}$ such that the $r$-leading form of $g$ is $d_g( a x+ b y+ c z)^r$.

As $(x,y,z)$ are good parameters for $I$,
$b_{00r}(g)\ne0$, so we must have $ c\ne 0$. If $ a\ne 0$, then $(1,0)$
is a vertex of $|\Delta(I;x,y,z)|$, which is not prepared since for $g\in I$
with $\nu_T(g)=r$,
$$
\{g\}_{xyz}^{(1,0)}= d_g c^r(z+\frac{ a}{ c} x)^r.
$$
Thus $ a=0$. If $ b\ne 0$, then $(0,1)$ is a vertex of $|\Delta(I;x,y,z)|$
which is not prepared since for $g\in I$ with $\nu_T(g)=r$,
$$
\{g\}_{xyz}^{(0,1)}= d_g c^r(z+\frac{ b}{ c} y)^r.
$$
Thus  $ a= b=0$ and $z=0$ is an approximate manifold of $I$.

The proof of the converse also follows from the above arguments.
\end{proof}

\begin{Lemma}\label{RSp12} Consider the terms in the expansion 
\begin{equation}\label{eqRSp13}
h=\sum_{\lambda=0}^k\eta^{k-\lambda}\binom{k}{\lambda}x^{i+(k-\lambda)a}
y^{j+(k-\lambda)b}z_1^{\lambda}
\end{equation}
obtained by substituting  $z_1=z-\eta x^ay^b$ into the  monomial $x^iy^jz^k$.
Define a projection for $(a,b,c)\in{\bold N}^3$ such that $c<r$,
$$
\pi(a,b,c) =(\frac{a}{r-c},\frac{b}{r-c}).
$$
\begin{enumerate}
\item[1.] Suppose that $k<r$. Then the exponents of monomials in (\ref{eqRSp13})
with nonzero coefficients  project into the line segment joining $(a,b)$ to 
$(\frac{i}{r-k},\frac{j}{r-k})$. 
\begin{enumerate}
\item[a.] If $(a,b)=(\frac{i}{r-k},\frac{j}{r-k})$, then all these monomials
project to $(a,b)$.
\item[b.] If $(a,b)\ne (\frac{i}{r-k},\frac{j}{r-k})$, then 
$x^{i}y^{j}z_1^k$ is the unique monomial in 
  (\ref{eqRSp13}) which projects onto $(\frac{i}{r-k},\frac{j}{r-k})$.
No monomial in (\ref{eqRSp13}) projects to $(a,b)$.
\end{enumerate}
\item[2.] Suppose that $r\le k$, and $(i,j,k)\ne (0,0,r)$. Then all exponents in (\ref{eqRSp13})
with nonzero coefficients and  $z_1$ exponent less than $r$ project into 
$$
\left((a,b)+{\bold Q}^2_{\ge 0}\right)-\{(a,b)\}.
$$
\item[3.] Suppose that $(i,j,k)=(0,0,r)$. Then all exponents in (\ref{eqRSp13})
with nonzero coefficients and $z_1$ exponent less than $r$ project to $(a,b)$.
\end{enumerate}
\end{Lemma}

Lemma \ref{RSp12} is proved by a straight forward calculation.
We deduce Lemma \ref{LemmaRSp15} from Lemma \ref{RSp12}.

\begin{Lemma}\label{LemmaRSp15}
Suppose that the vertex $(a,b)$ is not perpared on $|\Delta|$, and $z_1=z-\eta x^ay^b$ is an $(a,b)$ preparation. Then
\begin{enumerate}
\item[1.]  $|\Delta(I;x,y,z_1)|\subset |\Delta(I;x,y,z)|-\{(a,b)\}$.
\item[2.] If $(a',b')$ is another vertex of $|\Delta(I;x,y,z)|$,
then $(a',b')$ is a vertex of 
\linebreak $|\Delta(I;x,y,z_1)|$,
 and for $g\in I$ with $\nu_T(g)=r$,
$\{g\}^{a',b'}_{x,y,z_1}$ is obtained from  $\{g\}^{a',b'}_{x,y,z}$
by substituting $z_1$ for $z$. 
\end{enumerate}
\end{Lemma}

\begin{Lemma}\label{LemmaRSp16}  Suppose that $\mbox{Sing}_r(I)$ has dimension $<2$, $\nu_T(I)=r$, $\tau(I)=1$
and $(x,y,z)$ are good parameters for $I$. Then there is a formal series
 $\Psi(x,y)\in {\bf k}[[x,y]]$
such that under the substitution $z=z_1-\Psi(x,y)$, $(x,y,z_1)$ are good parameters
for $I$ and 
$(I;x,y,z_1)$ is well prepared.
\end{Lemma}

Lemma \ref{LemmaRSp16} is proved by successive removal of  vertices which are not
 prepared, using Lemma \ref{LemmaRSp15}. 

\begin{Remark}\label{RemarkL26}
It is not always possible to well prepare
after a finite number of  preparations, so that $z_1$ could be a
formal series in $z,x,y$. An example is given in Example 7.19 \cite{C2}.
\end{Remark}

\subsection{Very Well Preparation}\label{VWp}

We will also consider change of variables of the form
$$
y_1=y-\eta x^n
$$
for $\eta\in {\bf k}$, $n$ a positive integer, which we will call translations.

\begin{Lemma}\label{LemmaRSp207} Consider  the expansion 
\begin{equation}\label{eqRSp16}
h=\sum_{\lambda=0}^j\eta^{j-\lambda}\binom{j}{\lambda}x^{i+(j-\lambda)n}
y_1^{\lambda}z^k
\end{equation}
obtained by substituting  $y_1=y-\eta x^n$ into the  monomial $x^iy^jz^k$.
Consider the projection for $(a,b,c)\in{\bold N}^3$ such that $c<r$ defined by
$$
\pi(a,b,c) =(\frac{a}{r-c},\frac{b}{r-c}).
$$
Suppose that $k<r$. Set $(a,b)=(\frac{i}{r-k},\frac{j}{r-k})$.
Then $x^{i}y_1^{j}z^k$ is the unique monomial  
 in (\ref{eqRSp16}) whose coefficients project onto $(a,b)$.
All other monomials in (\ref{eqRSp16}) with non-zero coefficient project to points
below $(a,b)$ on the line through $(a,b)$ with slope $-\frac{1}{n}$.
\end{Lemma}

Lemma \ref{LemmaRSp207} is proved by a straight forward calculation. 
We deduce Lemma \ref{LemmaL30} from Lemma \ref{LemmaRSp207} and the definition of well preparedness.

\begin{Lemma}\label{LemmaL30}  Suppose that $(I;x,y,z)$ is well prepared, 
$y_1=y-\eta x^n$ is a translation, and $z_1=z-\psi(x,y_1)$ is a subsequent
well preparation. Then
$$
\alpha_{x,y_1,z_1}(I)=\alpha_{xyz}(I),
\beta_{x,y_1,z_1}(I)=\beta_{xyz}(I)\text{ and }
\gamma_{x,y_1,z_1}(I)=\gamma_{xyz}(I).
$$
\end{Lemma}

\begin{Definition}\label{DefinitionRSp22}
Suppose that $I\subset R$ is an ideal such that $\nu_T(I)=r$, $\tau(I)=1$, $\mbox{Sing}_r(I)$ has dimension $<2$ and $(x,y,z)$ are good parameters for $I$.
Let $\alpha=\alpha_{xyz}(I)$, $\beta=\beta_{xyz}(I)$, $\gamma=\gamma_{xyz}(I)$,
$\delta=\delta_{xyz}(I)$, $\epsilon=\epsilon_{xyz}(I)$.
$(I;x,y,z)$ will be said to be very well prepared if
it is well prepared and one of  
the following conditions holds.
\begin{enumerate}
\item[1.] $(\gamma-\delta,\delta)\ne (\alpha,\beta)$ and if we make a translation
$y_1=y-\eta x$, with subsequent well preparation $z_1=z-\Psi(x,y_1)$, then
$$
\alpha_{xy_1z_1}(I)=\alpha, \beta_{xy_1z_1}(I)=\beta,
\gamma_{xy_1z_1}(I)=\gamma
$$
and
$$
\delta_{xy_1z_1}(I)\le\delta.
$$
\item[2.] $(\gamma-\delta,\delta)=(\alpha,\beta)$ and one of the following cases hold:
\begin{enumerate}
\item[a.] $\epsilon=0$
\item[b.] $\epsilon\ne 0$ and $\frac{1}{\epsilon}$ is not an integer 
\item[c.] $\epsilon\ne 0$ and $n=\frac{1}{\epsilon}$ is a (positive) integer and for any $\eta\in {\bf k}$, if
$y_1=y-\eta x^n$ is a translation, with subsequent well preparation 
$z_1=z-\Psi(x,y_1)$, then $\epsilon_{xy_1z_1}(I)=\epsilon$.
Further, if $(c,d)$ is the lowest point on the line  through $(\alpha,\beta)$ with
slope $-\epsilon$ in $|\Delta(I;x,y,z)|$
and $(c_1,d_1)$ is the lowest point on this line in $|\Delta(I;x,y_1,z_1)|$,
then $d_1\le d$.
\end{enumerate}
\end{enumerate}
\end{Definition}

\begin{Lemma}\label{LemmaRSp18} 
Suppose that $I\subset T$ is an ideal such that $\nu_T(I)=r$, $\tau(I)=1$, $\mbox{Sing}_r(I)$ has dimension $<2$ and $(x,y,z)$ are good parameters for $I$.
Then there are formal substitutions
$$
z_1=z-\Psi(x,y), y_1=y-\phi(x)
$$
where $\Psi(x,y),\phi(x)$ are series  such that $(I;x,y_1,z_1)$ is very well prepared.
\end{Lemma}

\begin{proof} 
By Lemmas \ref{LemmaRSp16} and  \ref{LemmaL30}, we can find good parameters $x,y,z$ for $I$ such that $(I;x,y,z)$ is well prepared, $\alpha,\beta,\gamma$ do not change under
translation $y_1=y-\eta x$ followed by subsequent well preparation, and $\delta$ is maximal. If $(\gamma-\delta,\delta)\ne (\alpha,\beta)$,
then we have
achieved  case 1 of Definition \ref{DefinitionRSp22}, so that $(I;x,y,z)$ is very well prepared.

We now assume that $(I;x,y,z)$ is well prepared, and $(\alpha,\beta)=(\gamma-\delta,\delta)$.
If $\epsilon=0$ or $\frac{1}{\epsilon}$ is not an integer then $(I;x,y,z)$ is very well prepared.

Suppose that $n=\frac{1}{\epsilon}$ is an integer. 
We then choose
$\eta\in {\bf k}$ such that with the translation $y_1=y-\eta x^n$, and subsequent
well preparation, we maximize $d$ for points $(c,d)$ of the line through $(\alpha,\beta)$ with slope
$-\epsilon$ 
on the boundary of $|\Delta(I;x,y_1,z_1)|$.
  By Lemma \ref{LemmaL30},
$\alpha,\beta$ and $\gamma$ are not changed. If we now have that $d\ne \beta$, we are very well prepared.

 If the process does not end after a finite number of iterations, then 
we construct formal series $y'=y-\phi(x)$ and $z'=z-\psi(x,y)$ such that
$|\Delta(I;x,y',z')|$  has the single vertex $(\alpha,\beta)$,
and $(I;x,y',z')$ is thus very well prepared.
\end{proof}

\subsection{Effect of a permissible blow up on $\Delta$}\label{SubsecPB}

\begin{Definition}\label{DefinitionRSp30}
Suppose that $I\subset T$ is an ideal such that $\nu_T(I)=r$, $\tau(I)=1$,
$\mbox{Sing}_r(I)$ has dimension $<2$,  and $(x,y,z)$ are good parameters for $I$.

We consider 4 types of  transformations $T\rightarrow T_1$,
where $T_1$ is the completion of the local ring of a blow up of $T$,
and $T_1$ has regular parameters $(x_1,y_1,z_1)$ related to the regular parameters
$(x,y,z)$ of $T$ by one of the following rules. In all cases, 
if $\nu_{T_1}(I_1)=r$ and $\tau(I_1)=1$, then $(x_1,y_1,z_1)$ are good parameters for
$I_1$.
\begin{enumerate}
\item[Tr1] $\text{Sing}_r(I)=V(x,y,z)$, 
$$
x=x_1, y=x_1(y_1+\eta), z=x_1z_1,
$$
with $\eta\in {\bf k}$.
Then $I_1=\frac{1}{x_1^r}IT_1$ is the weak transform of $I$ in $T_1$.
\item[Tr2] $\text{Sing}_r(I)=V(x,y,z)$, 
$$
x=x_1y_1, y=y_1, z=y_1z_1.
$$
Then $I_1=\frac{1}{y_1^r}IT_1$ is the weak transform of $I$ in $T_1$.
\item[Tr3] 
$\text{Sing}_r(I)=V(x,z)$,
$$
x=x_1, y=y_1, z=x_1z_1.
$$
Then $I_1=\frac{1}{x_1^r}IT_1$ is the weak transform of $I$ in $T_1$.
\item[Tr4]
$\text{Sing}_r(I)=V(y,z)$,
$$
x=x_1, y=y_1, z=y_1z_1.
$$
Then $I_1=\frac{1}{y_1^r}IT_1$ is the weak transform of $I$ in $T_1$.
\end{enumerate}
\end{Definition}

\begin{Lemma}\label{LemmaRSp19}  Suppose that assumptions are as in Definition
\ref{DefinitionRSp30}. Further
suppose that 
$(x,y,z)$ and $(x_1,y_1,z_1)$ are related by a  transformation of one of the above types
Tr1 - Tr4, and $\nu_{T_1}(I_1)=r$, $\tau(I_1)=1$. 
 Then there is a 1-1 correspondence
$$
\sigma:\Delta(I,x,y,z)\rightarrow \Delta(I_1,x_1,y_1,z_1)
$$
defined by 
\begin{enumerate}
\item[1.] $\sigma(a,b)=(a+b-1,b)$ if the transformation is a Tr1 with $\eta=0$.
\item[2.] $\sigma(a,b)=(a,a+b-1)$ if the transformation is a Tr2.
\item[3.] $\sigma(a,b)=(a-1,b)$ if the transformation is a Tr3.
\item[4.] $\sigma(a,b)=(a,b-1)$ if the transformation is a Tr4.
\end{enumerate}
\end{Lemma}
The proof of Lemma \ref{LemmaRSp19} is a straight forward calculation.
We deduce Lemma \ref{LemmaRSp20} from Lemma \ref{LemmaRSp19}.

\begin{Lemma}\label{LemmaRSp20} In each of the four cases of the preceeding 
lemma, if $\sigma(a,b)=(a_1,b_1)$ is a vertex of $|\Delta(I_1;x_1,y_1,z_1)|$,
then $(a,b)$ is a vertex of $|\Delta(I;x,y,z)|$, and if
$(I;x,y,z)$ is $(a,b)$ prepared  then $(I_1;x_1,y_1,z_1)$ is $(a_1,b_1)$
prepared. In particular, $(I_1;x_1,y_1,z_1)$ is well prepared if 
$(I;x,y,z)$ is well prepared.
\end{Lemma}

\begin{Lemma}\label{LemmaRSp21}
Suppose that assumptions are as in Definition \ref{DefinitionRSp30} and
$(x,y,z)$ and $(x_1,y_1,z_1)$ are related by a Tr3 transformation. 
Suppose that $\nu_{T_1}(I_1)=r$ and $\tau(I_1)=1$.
If $(I;x,y,z)$ is very well prepared, then $(I_1;x_1,y_1,z_1)$
is very well prepared, 
$$
\beta_{x_1y_1z_1}(I_1)=\beta_{xyz}(I),
\delta_{x_1y_1z_1}(I_1)=\delta_{xyz}(I),
\epsilon_{x_1y_1z_1}(I_1)=\epsilon_{xyz}(I)
$$
and 
$$
\alpha_{x_1y_1z_1}(I_1)=\alpha_{xyz}(I)-1,
\gamma_{x_1y_1z_1}(I_1)=\gamma_{xyz}(I)-1.
$$
We further have $\text{Sing}_r(I_1)\subset V(x_1,z_1)$.
\end{Lemma}

\begin{proof} Well preparation is preserved by Lemma \ref{LemmaRSp20}. 
We  deduce from the
definition of Tr3 that $(I_1;x_1,y_1,z_1)$ is also very well prepared.

Since $\text{Sing}_r(I)
=V(x,z)$, we have  $\text{Sing}_r(I_1)\subset V(x_1)$, the exceptional divisor of our  Tr3 transformation. $\mbox{Sing}_r(I_1)$ is either the maximal ideal, or the germ of a nonsingular irreducible curve, by Lemma \ref{LemmaZL}.

Suppose that $\mbox{Sing}_r(I_1)$ is a nonsingular curve $C$. Then $C$ has formal local equations
$$
x_1=z_1-\sum_{t=1}^{\infty}a_ty_1^t=0
$$
 for some $a_t\in{\bf k}$
by Lemma \ref{Lemma107}.
Suppose that $g\in I$ with $\nu_T(g)=r$. Expand
$$
g=(\sum_{i+k=r}\sum_{j=0}^{\infty}b_{ijk}x^iy^jz^k)+h
$$
with $h\in (x,z)^{r+1}$. 
$$
\frac{g}{x_1^r}=\left(\sum_{i+k=r}\sum_{j=0}^{\infty}b_{ijk}y_1^jz_1^k)\right)
+x_1\frac{h}{x_1^{r+1}}
$$
with $\frac{h}{x_1^{r+1}}\in T_1$.

$\nu_T(\frac{g}{x_1^r})=r$ and
$\frac{g}{x_1^r}\in (x_1,z_1-\sum_{t=1}^{\infty}a_ty_1^t)^r$ implies
$$
\sum_{i+k=r}\sum_{j=0}^{\infty}b_{ijk}y_1^jz_1^k
=b_{00r}(z_1-\sum_{t=1}^{\infty}a_ty_1^t)^r.
$$
Thus
$$
g=b_{00r}(z-\sum_{t=1}^{\infty}a_txy^t)^r+h.
$$
Suppose that $a_t\ne 0$ for some $t$. Let $t_0$ be the smallest value
of $t$ for which $a_{t}\ne 0$. Then $(\alpha,\beta)=(1,t_0)$, since
$I\subset (x,z)^r$. But $(1,t_0)$ is then a vertex of 
$|\Delta(I;x,y,z)|$ which is not prepared, a contradiction. Thus
$a_t=0$ for all $t$ and $C=V(x,z_1)$.
We conclude that $\text{Sing}_r(I_1)\subset V(x_1,z_1)$

The remaining statements of the lemma follow from Lemma \ref{LemmaRSp19}.

\end{proof}

\begin{Lemma}\label{LemmaRSp23} Suppose that assumptions are as in 
Definition \ref{DefinitionRSp30}, $(x,y,z)$ and $(x_1,y_1,z_1)$ are related by a 
Tr1 transformation with $\eta=0$, $(I;x,y,z)$ is very well prepared and
$\nu_{T_1}(I_1)=r$, $\tau(I_1)=1$.
Let $(I_1;x_1,y',z')$ be a very well preparation of $(I_1;x_1,y_1,z_1)$. Let
$$
\alpha=\alpha_{xyz}(I), \beta=\beta_{xyz}(I), \gamma=\gamma_{xyz}(I),
\delta=\delta_{xyz}(I),\epsilon=\epsilon_{xyz}(I),
$$
$$
\alpha_1=\alpha_{x_1y'z'}(I_1), \beta_1=\beta_{x_1y'z'}(I_1), \gamma_1=\gamma_{x_1y'z'}(I_1),
\delta_1=\delta_{x_1y'z'}(I_1),\epsilon_1=\epsilon_{x_1y'z'}(I_1).
$$
Then 
$$
(\beta_1,\frac{1}{\epsilon}_1,\alpha_1)<(\beta,\frac{1}{\epsilon},\alpha)
$$
in the lexicographical order.
If   
$\beta_1=\beta$, and $\epsilon\ne 0$,
then 
$\frac{1}{\epsilon_1}=\frac{1}{\epsilon}-1$. 

We further have $\text{Sing}_r(I_1)\subset V(x_1,z')$.
\end{Lemma}

\begin{proof}  $(I_1;x_1,y_1,z_1)$ is well
prepared by Lemma \ref{LemmaRSp20}.

Let $\sigma:\Delta(I;x,y,z)\rightarrow\Delta(I_1;x_1,y_1,z_1)$
be the 1-1 correspondence of 1 in Lemma \ref{LemmaRSp19}, which is defined by
$\sigma(a,b)=(a+b-1,b)$.
$\sigma$ transforms lines of slope $m\ne -1$ to lines of slope $\frac{m}{m+1}$, and
transforms lines of slope -1 to vertical lines.

We deduce the formulas for transformation of $\alpha,\beta,\gamma,\delta,\epsilon$ 
from consideration of the effect of a Tr1 transformation followed by very well preparation on $(I;x,y,z)$ in the different cases of
 Definition \ref{DefinitionRSp22}. We need the assumption $\mbox{Sing}_r(I)=V(x,y,z)$ so that $\beta<1$, to conclude that $\alpha_1<\alpha$ if $\epsilon=0$.

By Lemmas \ref{Lemma107} and \ref{LemmaRSp14}, $\text{Sing}_r(I_1)\subset V(x_1,z_1)$. If $\text{Sing}_r(I_1)
=V(x_1,z_1)$, then all vertices of $|\Delta(I_1;x_1,y_1,z_1)|$ have the first 
coordinate $\ge 1$, and thus $\text{Sing}_r(I_1)=V(x_1,z')$.

\end{proof}

\begin{Lemma}\label{LemmaRSp24}Suppose that assumptions are as in 
Definition \ref{DefinitionRSp30}, $(x,y,z)$ and $(x_1,y_1,z_1)$ are related by a 
Tr2 or a Tr4 transformation   and $(I;x,y,z)$ is  well prepared.
Suppose that $\nu_{T_1}(I_1)=r$ and $\tau(I_1)=1$.
We have that $(I_1;x_1,y_1,z_1)$ is well prepared and $\text{Sing}_r(I_1)\subset V(y_1,z_1)$.
Let $(I_1;x_1,y',z')$ be a very well preparation. Then
$$
\beta_{x_1,y_1,z_1}(I_1)=\beta_{x_1,y',z'}(I_1)<\beta_{xyz}(I).
$$
\end{Lemma}

\begin{proof} First suppose that $T\rightarrow T_1$ is a Tr2 transformation.
The correspondence $\sigma:\Delta(I;x,y,z)\rightarrow\Delta(I_1;x_1,y_1,z_1)$
of 2 in Lemma \ref{LemmaRSp19} is defined by $\sigma(a,b)=(a,a+b-1)$. $\sigma$ takes lines
of slope $m$ to lines of slope $m+1$ Thus
$$
(\alpha_{x_1y_1z_1}(I_1),\beta_{x_1y_1z_1}(I_1)
=(\alpha_{xyz}(I),\alpha_{xyz}(I)+\beta_{xyz}(I)-1).
$$
Since $\text{Sing}_r(I)=V(x,y,z)$ by assumption, we have $\alpha_{xyz}(I)<1$, and thus
$$
\beta_{x_1y_1z_1}(I_1)<\beta_{xyz}(I).
$$
 $(I_1;x_1,y_1,z_1)$ is well prepared by
Lemma \ref{LemmaRSp20}. Thus the vertex $(\alpha_{x_1y_1z_1}(I_1),\beta_{x_1y_1z_1}(I_1))$
is not affected by very well preparation. We thus have $\beta_{x_1y'z'}(I_1)<\beta_{xyz}(I)$.

Now suppose that $T\rightarrow T_1$ is a Tr4 transformation. The 1-1 correspondence 
$\sigma:\Delta(I;x,y,z)\rightarrow\Delta(I_1;x_1,y_1,z_1)$ of 4 in Lemma 
\ref{LemmaRSp19} is defined by $\sigma(a,b)=(a,b-1)$. Thus
$\beta_{x_1y_1z_1}(I_1)<\beta_{xyz}(I)$.
$(I_1;x_1,y_1,z_1)$ is well prepared by Lemma \ref{LemmaRSp20}, and the vertex
$(\alpha_{x_1y_1z_1}(I_1),\beta_{x_1y_1z_1}(I_1))$ is not affected by very well
preparation. Thus
$\beta_{x_1y'z'}(I_1)=\beta_{x_1y_1z_1}(I_1)<\beta_{xyz}(I)$.

The proof that $\mbox{Sing}_r(I_1)\subset V(y_1,z_1)$ is as in the proofs of
Lemma \ref{LemmaRSp21} and Lemma \ref{LemmaRSp23}.
\end{proof}

\begin{Lemma}\label{LemmaRSp25} Suppose that $p$ is a prime,
$s$ is a nonnegative integer and
$r_0$ is a positive integer such that $p\not\,\mid r_0$. 
Let $r=r_0p^s$. Then
\begin{enumerate}
\item[1.] $\binom{r}{\lambda}\equiv 0\text{ mod }p, \text{ if }p^s\not\,\mid \lambda$, $0\le\lambda\le r$ is an integer.
\item[2.] $\binom{r}{\lambda p^s}\equiv \binom{r_0}{\lambda}
\text{ mod }p$, if $0\le\lambda\le r_0$ is an integer.
\end{enumerate}
\end{Lemma}

\begin{proof} Compare the expansions over ${\bold Z}_p$ of
$(x+y)^r=(x^{p^s}+y^{p^s})^{r_0}$.
\end{proof}

The following theorem is the most delicate part of the proof, where we really
see the difference between characteristic 0 and characteristic $p$.

\begin{Theorem}\label{TheoremRSp27}
Suppose that assumptions are as in 
Definition \ref{DefinitionRSp30}, $(x,y,z)$ and $(x_1,y_1,z_1)$ are related by a 
Tr1 transformation with $\eta\ne 0$ and
  $(I;x,y,z)$ is very well prepared. Suppose that $\nu_{T_1}(I_1)=r$,
$\tau(I_1)=1$ and $\beta_{xyz}(I)>0$.  Then there exist good parameters
 $(x_1,y_1',z_1')$ in $T_1$ such that  $(I_1;x_1,y_1',z_1')$
is very well prepared, with
  $\beta_{x_1,y_1',z_1'}(I_1)<\beta_{x,y,z}(I)$
and $\text{Sing}_r(I_1)\subset V(x_1,z_1')$.
\end{Theorem}

\begin{proof} Let $\alpha=\alpha_{xyz}(I)$, $\beta=\beta_{xyz}(I)$, $\gamma=\gamma_{xyz}(I)$, $\delta=\delta_{xyz}(I)$. Lemma \ref{LemmaRSp14} implies $z=0$ is an approximate manifold of $I$. Thus since $\tau(I)=1$,
\begin{equation}\label{eqX10}
\alpha+\beta>1.
\end{equation}
Apply the translation $y'=y-\eta x$ and well prepare by some substitution 
$z'=z-\Psi(x,y')$.
This does not change $\alpha$, $\beta$ or $\gamma$.
Set  $\delta'=\delta_{x,y',z'}(I)$.

First assume that $\delta'<\beta$. We have regular parameters $(x_1, y_1, \overline z_1)$ in $T_1$ such that  $x=x_1, y'=x_1y_1, 
z'=x_1\overline z_1$,  so
we have, by 1 in Lemma \ref{LemmaRSp19}, a 1-1 correspondence 
$$
\sigma:\Delta(I;x,y',z')\rightarrow \Delta(I_1;x_1,y_1,\overline z_1)
$$
defined by $\sigma(a,b)=(a+b-1,b)$. 
$(I_1;x_1,y_1,\overline z_1)$ is well prepared
by Lemma \ref{LemmaRSp20}. 
Since $\delta'<\beta$, the line segment in $|\Delta(I;x,y',z')|$ through 
$(\alpha,\beta)$ and $(\gamma-\delta',\delta')$ has slope $\le -1$, so it is
transformed by $\sigma$ to a segment with positive slope or a vertical line in
$|\Delta(I_1;x_1,y_1,\overline z_1)|$.
We have that for $\alpha_1=\alpha_{x_1y_1\overline z_1}(I_1)$ and
$\beta_1=\beta_{x_1y_1\overline z_1}(I_1)$, $(\alpha_1,\beta_1)=\sigma(\gamma-\delta',\delta')$ and
$\beta_1=\delta'<\beta$. Since $(I_1;x_1,y_1,\overline z_1)$ is well prepared, very well
preparation does not effect the vertex $(\alpha_1,\beta_1)$, by Lemma \ref{LemmaL30}.
 After very well preparation, we thus have regular parameters
$x_1,y_1',z_1'$ in $T_1$ such that $(I_1;x_1,y_1',z_1')$ is very well
prepared and $\beta_{x_1,y_1',z_1'}(I_1)<\beta_{xyz}(I)$.
$\mbox{Sing}_r(I_1)\subset V(x_1,z_1')$ as in the proof of Lemma \ref{LemmaRSp23}.

 We have thus reduced the proof to showing that
$\delta'<\beta$.

 If $(\alpha,\beta)\ne(\gamma-\delta,\delta)$, we have $\delta'\le\delta<\beta$
 by Lemma \ref{LemmaL30}, since $(I;x,y,z)$ is very well prepared.

Suppose that $(\alpha,\beta)=(\gamma-\delta,\delta)$. Set
$$
W=\left\{(i,j,k)\in{\bf N}^3\mid k<r\text{ and
}\left(\frac{i}{r-k},\frac{j}{r-k}\right)=(\alpha,\beta)\right\}.
$$
For $g=\sum a_{ijk}x^iy^jz^k\in I$, with $\nu_T(g)=r$, set
$$
F_g=\sum_{(i,j,k)\in W}a_{ijk}x^iy^jz^k.
$$
By assumption, $(\alpha,\beta)$ is  prepared on
$|\Delta(I;x,y,z)|$, which implies that if $\alpha,\beta$ are
integers, then there does not exist $\lambda\in{\bf k}$ such that
$$
a_{00r}z^r+F_g=a_{00r}(z-\lambda x^{\alpha}y^{\beta})^r
$$
for all $g\in I$ with $\nu_T(g)=r$. Moreover,
\begin{equation}\label{eq40}
F_g(x,y,z)=F_{g}(x,y'+\eta x,z)=\sum_{(i,j,k)\in
W}\sum_{\lambda=0}^ja_{ijk}\eta^{\lambda}\binom{j}{\lambda}
x^{i+\lambda}(y')^{j-\lambda}z^k.
\end{equation}

By   Lemma \ref{LemmaRSp207}, the terms in the expansion of
$g(x,y'+\eta x,z)$ contributing to $(\gamma,0)$ in
$|\Delta(I;x,y',z)|$, where $\gamma=\alpha+\beta$, are
$$
F_{g,(\gamma,0)}=\sum_{(i,j,k)\in W}a_{ijk}\eta^jx^{i+j}z^k.
$$
If $(I;x,y',z)$ is $(\gamma,0)$ prepared, then $\delta'=0<\beta$.
Suppose that $(\gamma,0)$ is not prepared  on
$|\Delta(I,x,y',z)|$. Then $\gamma\in{\bf N}$, and there exists
$\psi\in{\bf k}$ such that for all $g\in I$ with $\nu_T(g)=r$,
\begin{equation}\label{eq36}
a_{00r}(z-\psi x^{\gamma})^r=a_{00r}z^r+F_{g,(\gamma,0)},
\end{equation}
so that, with $\omega=\frac{-\psi}{\eta^{\beta}}\in{\bf k}$, for
$0\le k<r$, we have for $g\in I$ with $\nu_T(g)=r$:
\begin{enumerate}
\item[1.] If $\binom{r}{r-k}\ne 0$ (in ${\bf k}$) then
$i=\alpha(r-k), j=\beta(r-k)\in{\bf N}$, and
\begin{equation}\label{eq37}
a_{ijk}=a_{00r}\binom{r}{r-k} \omega^{r-k}.
\end{equation}
\item[2.] If $i=\alpha(r-k), j=\beta(r-k)\in{\bf N}$ and
$\binom{r}{r-k}=0$ (in ${\bf k}$), then $a_{ijk}=0$.
\end{enumerate}
Thus by Lemma \ref{LemmaRSp25}, for $(i,j,k)\in W$, $a_{ijk}= 0$ if $p^s\not\,\mid k$.

If $\bf k$ has characteristic zero, we obtain a contradiction to our
assumption that $(\gamma,0)$ in not prepared on $|\Delta(I,x,y',z)|$. Thus
$0=\delta'<\beta$ if $\bf k$ has characteristic zero.

Now we consider the case where ${\bf k}$ has characteristic $p>0$,
 and $r=p^sr_0$ with $p\not\,\mid r_0$, $r_0\ge 1$. Then by Lemma \ref{LemmaRSp25}, for $(i,j,k)\in W$,
   we have $a_{ijk}=0$  if $p^s\not\,\mid k$, for all $g\in I$ with $\nu_T(g)=r$.
By (\ref{eq37}), and Lemma \ref{LemmaRSp25}, we have that $i=\alpha
p^s,j=\beta p^s \in{\bf N}$, and for all $g\in I$ with
$\nu_T(g)=r$,    we have that  $a_{i,j, (r_0-1)p^s}\ne 0$.

 We
have an expression $\beta p^s=e p^t$ where $p\not\,\mid e$.
Suppose that $t\ge s$. Then $\beta\in{\bf N}$, which implies that
$\alpha=\gamma-\beta\in{\bf N}$, so that for all $g\in I$ with
$\nu_T(g)=r$,
$$
a_{00r}(z+\omega x^{\alpha}y^{\beta})^r=a_{00r}z^r+F_g,
$$
a contradiction, since $(\alpha,\beta)$ is by assumption prepared
on $|\Delta(I;x,y,z)|$. Thus $t<s$. Suppose that $e=1$. Then
$\beta=p^{t-s}<1$ and $\alpha<1$ (since we must have that
$\text{Sing}_r(I)=V(x,y,z)$) which implies (since $\gamma$ is an integer) that
$\gamma=\alpha+\beta=1$, a contradiction to (\ref{eqX10}). Thus $e>1$. Also,
\begin{equation}
\begin{array}{lll}
a_{00r}z^r+F_g&=& a_{00r}(z^{p^s}+\omega^{p^s}x^{\alpha p^s}y^{\beta p^s})^{r_0}\\
&=&a_{00r}(z^{p^s}+\omega^{p^s}x^{\alpha p^s}(y'+\eta x)^{\beta p^s})^{r_0}\\
&=&a_{00r}(z^{p^s}+\omega^{p^s}x^{\alpha
p^s}[(y')^{p^t}+\eta^{p^t} x^{p^t}]^{e})^{r_0}.
\end{array}
\end{equation}

Now make the $(\gamma,0)$ preparation $z=z'-\eta^{\beta}\omega
x^{\gamma}$ (from (\ref{eq36})) so that $(I;x,y',z')$ is
$(\gamma,0)$ prepared.  For $g\in I$ such that $\nu_T(g)=r$, let
$G_g=a_{00r}z^r+F_g$. Then
$$
\begin{array}{lll}
G_g&=& a_{00r}\left((z')^{p^s}+e\omega^{p^s}\eta^{p^t(e-1)}(y')^{p^t}x^{\alpha p^s+p^t(e-1)}+(y')^{2p^t}\Omega(x,y'))\right)^{r_0}\\
&=&a_{00r}\left[(z')^{p^sr_0}+r_0\left[ e\omega^{p^s}\eta^{p^t(e-1)}(y')^{p^t}x^{\alpha p^s+p^t(e-1)}+(y')^{2p^t}\Omega\right](z')^{p^s(r_0-1)}\right]\\
&&+\Lambda_2(x,y')(z')^{p^s(r_0-2)}+\cdots +\Lambda_{r_0}(x,y')
\end{array}
$$
for some polynomials (which depend on $g\in I$),
$\Omega(x,y'),\Lambda_2,\ldots,\Lambda_{r_0}$, where
$(y')^{ip^t}\mid \Lambda_i$ for all $i$. All contributions of
$S(\gamma)\cap|\Delta(I;x,y',z')|$ must come from these
polynomials $G_g$, for $g\in I$ with $\nu_T(g)=r$ since $\bf k$ is
infinite. Recall that we are assuming
$(\alpha,\beta)=(\gamma-\delta,\delta)$. The term of lowest second
coordinate on $S(\gamma)\cap |\Delta(I;x,y',z')|$ is
$$
(a,b)=\left(\frac{\alpha
p^s+p^t(e-1)}{p^s},\frac{p^t}{p^s}\right),
$$
which is not in ${\bf N}^2$ since $t<s$, and is not
$(\alpha,\beta)$ since $e>1$. $(a,b)$ is thus prepared on
$|\Delta(I;x,y',z')|$, and
$$
\delta'=\frac{p^t}{p^s}<\frac{e p^t}{p^s}=\beta.
$$

\end{proof}

\subsection{Construction of the sequence $\Omega(q_n)$}

\begin{Theorem}\label{TheoremRSp28} For all $n\in {\bold N}$,
there are  good parameters $(x_n,y_n,z_n)$ for $I_n=({\mathcal I}_n)_{q_n}R_n$ in $R_n$ 
such that 
$$
\begin{array}{llll}
&\Omega(I_{n+1};x_{n+1},y_{n+1},z_{n+1})&=&(\beta_{n+1},
\frac{1}{\epsilon_{n+1}},\alpha_{n+1})\\
<&\Omega(I_{n};x_{n},y_{n},z_{n})&=&(\beta_{n},\frac{1}{\epsilon_{n}},\alpha_{n}).
\end{array}
$$
Further, if $\beta_{n+1}=\beta_n$,  $\epsilon_{n+1}\ne\epsilon_n$ and $\frac{1}{\epsilon_n}\ne \infty$,
then
$$
\frac{1}{\epsilon_{n+1}}=\frac{1}{\epsilon_n}-1.
$$
\end{Theorem}

\begin{proof}
We  inductively  construct 
regular parameters $(x_n,y_n,z_n)$ in $R_n$ such that $(x_n,y_n,z_n)$ are good parameters for $I_n$, $\text{Sing}_r(I_n)\subset V(x_n,z_n)$ or
$\text{Sing}_r(I_n)\subset V(y_n,z_n)$ and $(I_n;x_n,y_n,z_n)$ is well prepared.
 If $\text{Sing}_r(I_n)\subset V(x_n,z_n)$,
we will further have that $(I_n;x_n,y_n,z_n)$ is very well prepared.
Let $\Omega_n=\Omega(I_n;x_n,y_n,z_n)$.

We first choose possibly formal regular parameters $(x,y,z)$ in $R$
 which are good parameters for $I$, such that
$|\Delta(I;x,y,z)|$ is very well prepared (see Remark \ref{RemarkL26}). 
By Lemma \ref{LemmaRSp14}, $z=0$ is an approximate mainfold for $I$.
Let 
$\alpha=\alpha_{x,y,z}(I)$, $\beta=\beta_{x,y,z}(I)$.

By our reduction in Section \ref{SecERIII}, $q$ is isolated in $\text{Sing}_r({\mathcal I})$.

Suppose that  regular parameters $(x_i,y_i,z_i)$ in $R_i$ have been defined 
for $i\le n$ as specified above.

If  $\text{Sing}_r(I_n)= V(y_n,z_n)$, then $R_n\rightarrow R_{n+1}$ must be
a Tr4 transformation, by Lemma \ref{Lemma107}, since $V(z_n)$ is an approximate manifold of $I_n$.
$R_{n+1}$ has regular parameters $x_{n+1}',y_{n+1}',z_{n+1}'$ defined by
$$
x_n=x_{n+1}', y_n=y_{n+1}', z_n=y_{n+1}'z_{n+1}'.
$$
Thus $(I_{n+1};x_{n+1}',y_{n+1}',z_{n+1}')$ is  well prepared
by Lemma \ref{LemmaRSp24}. 
Further $\text{Sing}_r(I_{n+1})\subset V(y_{n+1}',z_{n+1}')$ by Lemma \ref{LemmaRSp24}.
If $\text{Sing}_r(I_{n+1})=V(x_{n+1}',y_{n+1}',z_{n+1}')$, make a  change of variables, subtracting a series in $x_{n+1}', y_{n+1}'$ from $z_{n+1}'$,
 to very well prepare, with resulting variables $x_{n+1},y_{n+1},z_{n+1}$. Otherwise, set $x_{n+1}=x_{n+1}',y_{n+1}=y_{n+1}', z_{n+1}=z_{n+1}'$. We have $\beta_{n+1}<\beta_n$
by Lemma \ref{LemmaRSp24}, and $\Omega_{n+1}<\Omega_n$.

If  $\text{Sing}_r(I_n)= V(x_n,z_n)$ then 
$(I_n;x_n,y_n,z_n)$ is very well prepared.
$R_n\rightarrow R_{n+1}$ must be
a Tr3 transformation by Lemma \ref{Lemma107}, since $V(z_n)$ is an approximate manifold of $I_n$.
$R_{n+1}$ has regular parameters $x_{n+1},y_{n+1},z_{n+1}$ defined by
$$
x_n=x_{n+1}, y_n=y_{n+1}, z_n=x_{n+1}z_{n+1}.
$$
Thus $(I_{n+1};x_{n+1},y_{n+1},z_{n+1})$ is also 
very well prepared, $\alpha_{n+1}<\alpha_n$ and $\Omega_{n+1}<\Omega_n$
 by Lemma \ref{LemmaRSp21}. 
Further $\text{Sing}_r(I_{n+1})\subset V(x_{n+1},z_{n+1})$.

Now suppose that  $\text{Sing}_r(I_n)=V(x_n,y_n,z_n)$. Then 
$(I_n;x_n,y_n,z_n)$ is very well prepared, and $R_{n+1}$ must be a Tr1 or Tr2 transformation
of $R_n$ by Lemma \ref{Lemma107}, since $V(z_n)$ is an approximate manifold of $I_n$.

If
$R_{n+1}$ is a Tr2 transformation of $R_n$, then $R_{n+1}$ has regular parameters
$x_{n+1}',y_{n+1}',z_{n+1}'$ defined by 
$$
x_n=x_{n+1}y_{n+1}', y_n=y_{n+1}',z_n=y_{n+1}z_{n+1}'.
$$
 $(I_{n+1};x_{n+1},y_{n+1}',z_{n+1}')$
is well prepared   by Lemma \ref{LemmaRSp24}. 
Further, $\text{Sing}_r(I_{n+1})\subset V(y_{n+1}',z_{n+1}')$ by Lemma \ref{Lemma107}. 
If $\text{Sing}_r(I_{n+1})=V(x_{n+1},y_{n+1}',z_{n+1}')$, then
make a further change of variables  to very well prepare.
Otherwise, set $y_{n+1}=y_{n+1}', z_{n+1}=z_{n+1}'$. $\beta_{n+1}<\beta_n$ by Lemma \ref{LemmaRSp24}, and $\Omega_{n+1}<\Omega_n$.

If $R_{n+1}$ is a Tr1 transformation of $R_n$,
then $R_{n+1}$ has regular parameters
$x_{n+1}',y_{n+1}',z_{n+1}'$ defined by 
$$
x_n=x_{n+1},
y_n=x_{n+1}(y_{n+1}'+\eta),
z_n=x_{n+1}z_{n+1}'.
$$
$\text{Sing}_r(I_{n+1})\subset V(x_{n+1},z_{n+1}')$ by Lemma \ref{Lemma107}.  

If $\beta_n\ne 0$ and $\eta\ne 0$ in the Tr1 transformation relating $R_n$ and $R_{n+1}$,
we  can change  variables to
$x_{n+1},y_{n+1}$, $z_{n+1}$ to very well prepare, with 
$\beta_{n+1}<\beta_n$, $\Omega_{n+1}<\Omega_n$ and $\text{Sing}_r(I_{n+1})\subset V(x_{n+1},z_{n+1})$
by Theorem \ref{TheoremRSp27}.
 
If $\eta=0$ in the Tr1 transformation of $R_n$, and $x_{n+1},y_{n+1},z_{n+1}$
are regular parameters of $R_{n+1}$ obtained from $x_{n+1},y_{n+1}',z_{n+1}'$
by very well preparation, then  by Lemma \ref{LemmaRSp23}, $(I_{n+1};x_{n+1},y_{n+1},z_{n+1})$ is very well prepared and $\Omega_{n+1}<\Omega_n$. If $\beta_{n+1}=\beta_n$,  and $\epsilon_n\ne 0$, then
$$
\frac{1}{\epsilon_{n+1}}=\frac{1}{\epsilon_n}-1.
$$ 
Further, $\text{Sing}_r(I_{n+1})\subset V(x_{n+1},z_{n+1})$ .

If $\beta_n=0$ we can make the translation $y_n'=y_n-\eta x_n$, and have that
$(I_n,x_n,y_n',z_n)$ is very well prepared, and $\beta_n,\gamma_n,\delta_n,\epsilon_n,\alpha_n$ are unchanged. Thus we are in the case of $\eta=0$.
\end{proof}

We now prove that (\ref{eq35}) cannot have infinite length.
Set $\Omega(q_n)=\Omega(I_n;x_n,y_n,z_n)$ for $n\in{\bf N}$. By Theorem
\ref{TheoremRSp28}, we have a contradiction if (\ref{eq35}) has infinite length, as the sequence $\Omega(q_n)$ cannot decrease indefinitely.

\section{Embedded Resolution on a Nonsingular 3-fold}
\label{SecERSI}

Our main resolution theorem for ideals on a nonsingular 3-fold, Theorem
\ref{Theorem6}, has now been proven, in Sections \ref{Section4} - \ref{Sect1}.

We now deduce embedded resolution of surface singularities
(Theorem \ref{Theorem1}) and principalization of ideals (Theorem \ref{Theorem23}) from Theorem \ref{Theorem6}. 

To prove Theorem \ref{Theorem1}, we start with ${\mathcal
R}=(\emptyset,E,{\mathcal I}_S,V)$. Let $r=\nu({\mathcal R})$ and apply
Theorem \ref{Theorem6}, to produce a sequence of permissible transforms
$$
V_{n_1}\rightarrow V_{n_1-1}\rightarrow \cdots \rightarrow V
$$
such that $\mbox{Sing}_r({\mathcal R}_{n_1})=\emptyset$. Since ${\mathcal I}_S$ is locally principal,
the weak transform ${\mathcal I}_i$ of ${\mathcal I}_S$ is the ideal sheaf of the strict transform $S_i$ of $S$ on $V_i$ for all
$i$.

Now we redefine ${\mathcal R}_{n_1}$ to be $(\emptyset, E_{n_1},{\mathcal I}_{n_1}, V_{n_1})$. By descending induction on $r$, we produce a sequence of permissible transforms
$$
V_{n_2}\rightarrow V_{n_2-1}\rightarrow \cdots \rightarrow V_{n_1}
$$
such that $\nu({\mathcal R}_{n_2})=1$. Now redefine ${\mathcal R}_{n_2}$ to be
$(\emptyset, E_{n_2},{\mathcal I}_{n_2}, V_{n_2})$. We have that $S_{n_2}$ is nonsingular,
since $\nu({\mathcal R}_{n_2})=1$. We now apply Lemma \ref{LemmaFL} to
$\mathcal{ R}_{n_2}$, to construct a sequence of permissible transforms
$$
V_{n}\rightarrow V_{n-1}\rightarrow \cdots \rightarrow V_2,
$$
obtained by blowing up points and curves in the locus where $E_i+S_i$ is not
a SNC divisor, to obtain  $V_n$ such that the preimage $E_n+S_n$ of $E+S$ on $V_n$
is a SNC divisor.
\vskip .2truein

To prove Theorem \ref{Theorem23}, 
 we observe that we can factor ${\mathcal I}={\mathcal J}{\mathcal K}$ in ${\mathcal O}_V$, where ${\mathcal J}$ is an invertible ideal sheaf, and the support of ${\mathcal O}_V/{\mathcal K}$ has codimension $\ge 2$ in $V$.
Let ${\mathcal R}=(\emptyset,\emptyset,{\mathcal K},V)$ and $r=\nu({\mathcal R})$.
By Theorem \ref{Theorem6}, there exists a sequence of permissible transforms
$$
V_{n_1}\rightarrow V_{n_1-1}\rightarrow \cdots \rightarrow V,
$$
such that $\mbox{Sing}_r({\mathcal R}_{n_1})=\emptyset$. After each permissible transform, we have that the weak transform ${\mathcal K}_i$ of ${\mathcal K}_{i-1}$ is such that the support of ${\mathcal O}_{V_i}/{\mathcal K}_i$
has codimension $\ge 2$ in $V_i$. Now redefine ${\mathcal R}_{n_1}$
to be $(\emptyset,\emptyset,{\mathcal K}_{n_1},V_{n_1})$. By descending induction on $r$, we produce a sequence of permissible transforms
$$
V_{n_2}\rightarrow V_{n_2-1}\rightarrow \cdots \rightarrow V,
$$
such that $\mbox{Sing}_2({\mathcal R}_{n_2})=\emptyset$. Now we have that
$\mbox{Spec}({\mathcal O}_{V_{n_2}}/{\mathcal K}_{n_2})$ is nonsingular, of
codimension $\ge 2$ in $V_{n_2}$. Blowing up the irreducible components of 
$\mbox{Spec}({\mathcal O}_{V_{n_2}}/{\mathcal K}_{n_2})$ is a sequence of permissible transforms above $V_{n_2}$, after which the conclusions of Theorem
\ref{Theorem23} are achieved.

\section{Projection to Points of Small Multiplicity}\label{SecProj}

In this section we prove   Theorem \ref{Theorem19}, the projection 
theorem which allows us to reduce to singularities of small multiplicity.
The method is due to Albanese. Our exposition is close to 
Lipman's survey  in \cite{L1}. We fill in a few technical difficulties.

Suppose that ${\bf k}$ is an algebraically closed field.
Let $W$ be a  projective variety over ${\bf k}$ of dimension $d$,
with an embedding $W\subset {\bf P}^{m}$.  Suppose that $q\in W$.
Let $\pi:{\bf P}^m\rightarrow {\bf P}^{m-1}$ be the rational map
which is projection from the point $q$. Let $W_1$ be the
projective subvariety of ${\bf P}^{m-1}$ which is the image of $W$
by $\pi$ (the closure of $\pi(W-\{q\})$). Let $\mu$ be the
multiplicity of the local ring ${\mathcal O}_{W,q}$.

\begin{Theorem}\label{Theorem47} The following are true.
\begin{enumerate}
\item[1.] $\mu\le \text{deg }W$. \item[2.] Suppose that
$\mu<\text{deg }W$. Then $\text{dim W}=\text{dim }W_1$ and
\begin{equation}\label{eq43}
[W:W_1]\text{deg }W_1=\text{deg }W-\mu.
\end{equation}
\item[3.] Suppose that $\mu=\text{deg }W$. Then $\text{dim }
W>\text{dim }W_1$ and $W$ is a cone over $W_1$ with vertex $q$.
\end{enumerate}
\end{Theorem}

\begin{proof}
Let $\sigma:Z\rightarrow {\bf P}^m$ be the blow up of $q$, with
exceptional divisor $E$. Let $\lambda:Z\rightarrow {\bf P}^{m-1}$
be the morphism induced by $\pi\circ\sigma$ (a resolution of the
indeterminacy  of $\pi$).

Let $H_0$ be a  hyperplane  of ${\bf P}^{m}$ and let $H_1$ be a
hyperplane  of ${\bf P}^{m-1}$.

We have a linear equivalence of  divisors
\begin{equation}\label{eq90}
\sigma^*(H_0)-E\sim \lambda^*(H_1).
\end{equation}

Let $\overline W$ be the strict transform of $W$ on $Z$.

By (\ref{eq90}), and since a general hyperplane of ${\bf P}^m$
does not contain $q$, we have equality of intersection numbers
\begin{equation}\label{eq41}
{\textstyle \int_Z}\lambda^*(H_1)^d\cdot \overline W={\textstyle
\int_Z}\sigma^*(H_0)^d\cdot\overline W+{\textstyle
\int_Z}(-E)^d\cdot \overline W.
\end{equation}
By the projection formula (Proposition 2.3 \cite{F}),
$$
{\textstyle \int_Z}\sigma^*(H_0)^d\cdot\overline W= {\textstyle
\int_{{\bf P}^m}} H_0^d\cdot W =\text{deg }W.
$$
Let $R={\mathcal O}_{W,q}$, $M$ be the maximal ideal of $R$. Then
the scheme-theoretic intersection of $\overline W$ and $E$ is
$$
\overline W\cap E=\text{Proj}(\oplus_{n\ge 0}M^n/M^{n+1}),
$$
and
$$
{\mathcal O}_Z(-E)\otimes{\mathcal O}_{\overline W\cap E}\cong
{\mathcal O}_{\overline W\cap E}(1).
$$
Thus by Lemma \ref{Lemma44},
$$
{\textstyle \int_Z}(-E)^d\cdot \overline W=-{\textstyle
\int_Z}(\overline W\cdot E) \cdot (-E)^{d-1}=-\text{deg }\overline
W\cap E= -e(R)=-\mu.
$$

We can now rewrite (\ref{eq41}) as
\begin{equation}\label{eq42}
{\textstyle \int_Z}\lambda^*(H_1)^d\cdot \overline W=\text{deg
}W-\mu.
\end{equation}

By the projection formula, we have that
$$
{\textstyle\int_Z}\lambda^*(H_1)^d\cdot\overline W=
\left\{\begin{array}{ll}
    [W:W_1]\text{deg }W_1&\text{ if $\text{dim }W=\text{dim }W_1$}\\
    0&\text{ if $\text{dim }W>\text{dim }W_1$}\end{array}\right.
$$

Substituting into (\ref{eq42}), we conclude that $\mu\le \text{deg
}W$, and  $\mu=\text{deg }W$ if and only if $W_1$ has dimension
$<d$, which holds if and only if $W$ is a cone with vertex $q$. If
$\mu<\text{deg }W$, we obtain (\ref{eq43}).
\end{proof}

We will say that the induced rational map $\pi:W\rightarrow W_1$
is a permissible projection if $\mu< \text{deg }W$.

Suppose that $W$ is a projective variety. We define $c(W)$ to be
the minimum of degrees of irreducible curves on $W$.

\begin{Theorem}\label{Theorem50} Suppose that $K$ is an algebraic function field
of dimension $d$ over an algebraically closed field ${\bf k}$.
Then there exists 
 a projective variety $V_0$, whose function field is $K$, and an embedding of $V_0$ into a projective space
  ${\bf P}^N$ so that $V_0$ is not contained in a hyperplane of ${\bf P}^N$, and such that
\begin{equation}\label{eq51*}
d!+1\ge \frac{\text{deg }V_0}{(N-d-\frac{1}{c}\text{deg }V_0)}
\end{equation}
where $c=c(V_0)$.
\end{Theorem}

\begin{proof}
There exists a  projective variety $V$, whose function field  is
$K$. Let $H$ be a very ample divisor on $V$.

After possibly replacing $H$ with a high multiple of $H$,
 we may assume that
$$
H^i(V,{\mathcal O}_V(\delta
H))=0
$$
 for all $\delta\ge 1$ and $i>0$. Thus we have
$$
h^0(V,{\mathcal O}_V(\delta H))=\chi({\mathcal O}_V(\delta H))
$$
for all $\delta\ge 1$.  $H^0(V,{\mathcal
O}_V(\delta H))$ gives
 an embedding $\Phi_{\delta}$ of $V$ into a projective space ${\bf P}^{N(\delta)}$,
where $N(\delta)= h^0(V,{\mathcal O}_V(\delta H))-1$.  Let
$V^{\delta}$ be the image of $V$.

Let $H_{\delta}$ be a hyperplane section of $V^{\delta}$. By our
construction, we have that $V^{\delta}$ is not contained in a
hyperplane of $P^{N(\delta)}$, and
$\Phi_{\delta}^*(H_{\delta})\sim \delta H$, so that every curve in
$V^{\delta}$ has degree $\ge \delta$. Thus
$c_{\delta}=c(V^{\delta})\ge\delta$.

We have that $\frac{1}{d!}\text{deg }V^{\delta}$ is the
coefficient of $n^d$ in the polynomial $\chi({\mathcal
O}_{V^{\delta}}(nH_{\delta}))$.

Since by our construction, $\chi({\mathcal
O}_{V^{\delta}}(nH_{\delta}))=\chi({\mathcal O}_V(\delta nH))$, we
see that $\text{deg }V^{\delta}=\delta^d\text{deg }V$.

There exist $a_0,a_1,\ldots, a_{d-1}\in{\bf Q}$ such that for all
$\delta>0$,
$$
N(\delta)+1=\chi({\mathcal O}_V(\delta
H))=\frac{1}{d!}\delta^d\text{deg }V+\sum_{i=0}^{d-1}a_i\delta^i.
$$
Thus
$$
\frac{\text{deg
}V^{\delta}}{N(\delta)-d-\frac{1}{c_{\delta}}\text{deg
}V^{\delta}} = \frac{\delta^{d}\text{deg }V}{
\frac{\delta^d}{d!}\text{deg
}V+\sum_{i=0}^{d-1}a_i\delta^i-1-d-\frac{\delta^d}{c_{\delta}}\text{deg
}V}.
$$
Since $c_{\delta}\ge \delta$ for all $\delta$, the limit as
$\delta$ goes to infinity of this expression is $d!$, and the
conclusions of the theorem follow with $V_0=V^{\delta}$ for
$\delta$ sufficiently large.

\end{proof}

\begin{Theorem}\label{Theorem48}
Let $K$ be an algebraic function field of dimension $d$ over an
algebraically closed field ${\bf k}$. Let $V_0\subset {\bf P}^N$
be as in the conclusions of Theorem \ref{Theorem50}.  Then there
exists a series of permissible projections
\begin{equation}\label{eq49}
V_0\stackrel{\pi_0}{\rightarrow}
V_1\stackrel{\pi_1}{\rightarrow}\cdots\stackrel{\pi_{n-1}}{\rightarrow}
V_n
\end{equation}
such that every point on $V_n$ has multiplicity $\le
\frac{d!}{[V_0:V_n]}$.
\end{Theorem}

\begin{proof}
Suppose that there exists a point $q_0\in V_0$ such that $q_0$ has
multiplicity $\mu_0$ with $d!<\mu_0<\text{deg }V_0$. Let $\pi_0$
be the projection from $q_0$. Let $V_1$ be the image of $V_0$. If
there exists a point $q_1\in V_1$ of multiplicity $\mu_1$ with
$\frac{d!}{[V_0:V_1]}<\mu_1<\text{deg }V_1$, then we perform the
projection $\pi_1:V_1\rightarrow V_2$ from the point $q_1$. After
a finite number of steps this process must terminate, as the
dimension of the ambient projective space drops by one after each
projection. After the final projection, we have a variety $V_n$
such that if $q\in V_n$ and the multiplicity $\mu$ of $q$ on $V$
satisfies $\mu >\frac{d!}{[V_0:V_n]}$, then $\mu=\text{deg }V_n$.
We see that $V_n$ satisfies the conclusions of the theorem unless
$V_n$ is a cone. We thus assume that $V_n$ is a cone with vertex
$q$ (in the projective space ${\bf P}^{N-n}$), and will show that
this is impossible.

By (\ref{eq43}), we have that
$$
[V_0:V_1]\text{deg }V_1=\text{deg }V_0-\mu_0\le \text{deg
}V_0-d!-1.
$$
We further have that
$$
[V_1:V_2]\text{deg }V_2\le\text{deg }V_1-\frac{d!}{[V_0:V_1]}-1,
$$
from which we obtain
$$
[V_0:V_2]\text{deg }V_2\le [V_0:V_1]\text{deg
}V_1-d!-[V_0:V_1]\le\text{deg }V_0-2d!-2.
$$
Continuing in the way, we obtain
$$
[V_0:V_n]\text{deg }V_n\le \text{deg }V_0-n(d!+1).
$$
We deduce the inequality
\begin{equation}\label{eq51}
n(d!+1)\le\text{deg }V_0.
\end{equation}

There exists a linear subspace $L$ of ${\bf P}^N$ of dimension
$n-1$ such that $\pi_{n-1}\circ\cdots \circ \pi_0$ is induced by
the projection morphism $\Psi:{\bf P}^N-L\rightarrow {\bf
P}^{N-n}$ from $L$. There exist dense open subsets $U_0\subset
V_0$ and $U_n\subset V_n$ such that
$\phi=(\pi_{n-1}\circ\cdots\circ\pi_0\mid U_0):U_0\rightarrow U_n$
is a finite morphism, and $U_0\cap L=\emptyset$ (so that $V_0\cap
\Psi^{-1}(U_n)=U_0$). Let $E=V_n-U_n$.

Since ${\bf k}$ is infinite, and by Bertini's theorem for general
hyperplane sections (Corollary 3.4.14 and 3.4.10 \cite{FOV}),
there exists a  hyperplane $\overline H$ of ${\bf P}^{N-n}$ which
does not contain $q$, intersects each irreducible component of $E$
 generically, and such that the scheme theoretic intersection $W=\overline H\cap V_n$
   is a $d-1$ dimensional variety. $V_n$ is the locus of lines through $q$ and a point of $W$.

By Bertini's theorem for general  hyperplane sections, there exist
hyperplanes
$$\overline H_1,\ldots,\overline H_{d-1}$$
 of ${\bf P}^{N-n}$ such that each $\overline H_i$ intersects $W$ and every irreducible component of $E$ generically.
In particular, the scheme-theoretic intersection $\overline
H_1\cap\cdots\cap\overline H_{d-1}\cap W$ is a reduced set of
points of order $\text{deg }V_n$ which are contained in $U_n$.

For $1\le i\le d-1$, let $H_i$ be the hyperplane  of ${\bf
P}^{N-n}$ which is  spanned by $q$ and the linear space $\overline
H_i\cap \overline H$. By our construction,  the scheme theoretic
intersection $H_1\cap\cdots\cap H_{d-1}\cap V_n$ is the union of
$s=\text{deg }V_n$ distinct lines $L_1,\ldots,L_s$, whose generic
points lie in $U_n$.

Let $H_1^*,\ldots, H_{d-1}^*$ be the hyperplanes of ${\bf P}^N$
such that $H_i^*\cap ({\bf P}^N-L)=\Psi^*(H_i)$. Since
$U_0\rightarrow U_n$ is finite, the irreducible components of
$$
H_1^*\cap\cdots\cap H_{d-1}^*\cap U_0=\phi^{-1}(H_1\cap \cdots
\cap H_{d-1}\cap U_n)
$$
are curves which dominate the irreducible components of $H_1\cap
\cdots \cap H_{d-1}\cap U_n$, so there are $\ge \text{deg }V_n$
distinct irreducible components of $H_1^*\cap\cdots\cap
H_{d-1}^*\cap U_0$. By the weak B\'ezout theorem (12.3.1) of
\cite{Ab7}, or Example 12.3.1 of the refined B\'ezout theorem of
\cite{F}, we have that

\begin{equation}\label{eq61}
\text{deg }V_0\ge c\,\text{deg }V_n,
\end{equation}
where $c=c(V_0)$.

We have that $V_n$ is not contained in a hyperplane of ${\bf
P}^{N-n}$, since $V_0$ is not contained in a hyperplane section of
${\bf P}^N$. Thus   we have the classical degree bound (Example
8.4.6 \cite{F}),
$$
\text{deg }V_n>N-n-d.
$$
From (\ref{eq61}), we get that
$$
(N-n-d)c<\text{deg }V_0.
$$
That is,
\begin{equation}\label{eq62}
N-d-\frac{1}{c}\text{deg }V_0<n.
\end{equation}
Now from (\ref{eq62}) and (\ref{eq51}), we have
$$
N-d-\frac{1}{c}\text{deg }V_0<\frac{1}{(d!+1)}\text{deg }V_0,
$$
which contradicts the assumption of our theorem. Thus $V_n$ is not
a cone.
\end{proof}

We now prove Theorem \ref{Theorem19}.

Let notation be as in the statement of Theorem \ref{Theorem48}. By
Theorem \ref{Theorem48}, there exists a dominant rational map
$V_0\rightarrow V_n$ such that every point on $V_n$ has
multiplicity $\le \frac{d!}{[V_0:V_n]}$. Let $V$ be the
normalization of $V_n$ in the function field $K$ of $V_0$, so that
the function field of $V$ is $K$, and we have a finite morphism
$\Psi:V\rightarrow V_n$. By Theorem \ref{Theorem71}, we have that
every point on $V$ has multiplicity $\le d!$.

\section{Ramification}\label{SecRam}

In this section, we collect some basic results on ramification of
local rings, which we will need for our proof of local uniformization in the next section. For more details, the book ``Ramification Theoretic
Methods in Algebraic Geometry'' \cite{Ab4} by  Abhyankar is an
excellent reference.

Suppose that $R$ is a normal local domain with quotient field $K$,
and $L$ is a finite  extension of $K$.  Let $T$ be the integral
closure of $R$ in $L$. The finitely many local rings $S$ which are
localizations of $T$   and dominate $R$ are called the extensions
of $R$ to $L$.

$R$ is said to be unramified in $L$ if for every extension $S$ of
$R$ to $L$ we have that $S$ is residually separable over $R$ and
$M(R)S=M(S)$.

Suppose that $X$ and $Y$  are normal integral schemes, and
$f:X\rightarrow Y$ is a finite morphism. The branch locus of
$f$ is the  set of points in $Y$ over which $f$ is ramified. The
branch locus is a closed subset of $Y$ (cf Theorems 1.44
and 1.44A \cite{Ab4}). If $Y$ is nonsingular, then the theorem of
the purity of the branch locus \cite{N}  tells us that if the
function field of $X$  is separable over the function field of
$Y$, then the branch locus of $f$ has pure codimension 1 in $Y$. In
particular, if $L$ is separable over $K$ and $R$ is a regular
local ring, then $R$ is unramified in $L$ if $R_Q$ is unramified
in $L$ for every height one prime $Q$ of $R$.

Suppose that our local ring $R$ is an equicharacteristic complete
regular local ring, with algebraically closed residue field ${\bf
k}$. Then $R$ is a power series ring over ${\bf k}$. $T$ is itself
a complete local ring.  Further, if $R$ is unramified in $L$, then
$R\cong T$ (as follows from Nakayama's Lemma).

Suppose that $f(Z)\in K[Z]$ is a monic polynomial.  Let
$\delta(f)\in K$ be the discriminant of $f$.  Suppose that $f(Z)$
factors as $f(Z)=(Z-z_1)\cdots(Z-z_n)$ in a splitting field of $f$
over $K$. Then $\delta(f)= \prod_{i\ne j}(z_i-z_j)$.

If $f(Z)\in R[Z]$, then since $R$ is integrally closed,
$\delta(f)\in R$.

We compute the discriminant in a case which will be important in
the sequel. Suppose that $\lambda\in K$. The discriminant  of
$f(Z)=Z^n-\lambda$ is
\begin{equation}\label{eq13}
\delta(f)= \pm n^n\lambda^{n-1}
\end{equation}
  (This can be computed directly from the definition or from the formula of Exercise 12 (c) of Chapter IV \cite{La}).

\begin{Lemma}\label{Lemma12}
Suppose that $R$ is a  normal local domain with quotient field
$K$, $L=K(r_1,\ldots,r_n)$ is a finite extension field of $K$, and
$f_i(Z)\in R[Z]$ are nonconstant monic polynomials such that
$f_i(r_i)=0$ for $1\le i\le n$. Suppose  that $\delta(f_i)\not\in
M(R)$ for $1\le i\le n$. Then $R$ is unramified in
$K(r_1,\ldots,r_n)$.
\end{Lemma}

This is (10.17) \cite{Ab7}.

\begin{Lemma}\label{Lemma17} Let $R$ be a $d$-dimensional regular local domain with $d>0$. Let $K$ be the quotient field of $R$ and let $(y_1,\ldots,y_d)$ be a basis of $M(R)$. Suppose that
$L=K(z_1,\ldots,z_d)$ is an overfield of $K$ where
$z_i^{n(i)}=y_i$ for $1\le i\le d$, where $n(i)$ are positive
integers. Let $S$ be the integral closure of $R$ in $L$. Then $S$
is a finite $R$-module, $S$ is a $d$-dimensional regular local
domain, $S=R[z_1,\ldots,z_d]$, $M(S)=(z_1,\ldots,z_d)S$,
$S/M(S)=R/M(R)$, and $[L:K]=n(1)\cdots n(d)$.
\end{Lemma}

This is (10.20.1) \cite{Ab7}.

Suppose that $R$ is one dimensional, and $S_i$ for $1\le i\le n$
are the distinct local rings of $T$ which dominate $R$.  Let $w_i$
be the integers such that $M(R)S_i= M(S)^{w_i}$. Define
$q_i=[S_i/M(S_i):R/M(R)]$. We  have an inequality
\begin{equation}\label{eq14}
q_1w_1+\cdots+q_nw_n\le[L:K]
\end{equation}
(by Theorem 21, Section 9, Chapter V \cite{ZS}).

We say that $R$ is tamely ramified in $L$ if $w_i$ is not
divisible by the characteristic of $R/M(R)$ and $S_i/M(S_i)$ is
separable over $R/M(R)$ for all $i$.

\section{Local Uniformization of Points of Small Multiplicity}\label{SecSmall}

In this section, we generalize the method of Jung \cite{J}, for locally resolving
singularities. Our methods are a simplification of  those of Abhyankar
 \cite{Ab7}. We prove local unformization of the function field of a
3-fold $V$ whose base field is algebraically closed of characteristic $>5$ in Theorem \ref{Theorem18}.

\begin{Theorem}\label{Theorem11}
Let $R$ be an  equicharacteristic, complete regular local domain
of dimension $d\ge 1$, whose residue field is algebraically closed
(so that $R$ is a power series ring over an algebraically closed
field ${\bf k}$).    Let $p$ be the characteristic of ${\bf k}$.
Suppose that $L$ is a finite extension of the quotient field $K$
of $R$. Further  suppose that a regular system of parameters
$y_1,\ldots, y_d$  in $R$ is such that $R$ is unramified in $L$
away from $V(y_1\cdots y_d)$, and the one dimensional regular
local rings $R_{(y_i)}$ are tamely ramified in $L$ for $1\le i\le
d$.

 Then there exists a positive integer $n$ such that $p\not\,\mid n$, and an inclusion of fields
$$
K\subset L\subset K(z_1,\ldots,z_d)
$$
such that  $z_i^n=y_i$ for $1\le i\le d$.
\end{Theorem}

\begin{proof}
Let $T$ be the integral closure of $R$ in $L$, $H_i=R_{(y_i)}$ for
$1\le i\le d$. Let $H_{i,j}$ for  $1\le j\le u(i)$ be the distinct
local rings of $T$ which dominate $H_i$ for $1\le i\le d$. Let
$w_{i,j}$ be the integers such that $M(H_i)H_{i,j}=
M(H_{ij})^{w_{i,j}}$. Define $$
q_{i,j}=[H_{i,j}/M(H_{ij}):H_i/M(H_i)].
$$

Let $n$ be the least common multiple of the $w_{i,j}$. Since the
$H_i$ are tamely ramified in $L$, we have that $p\not\,\mid n$.

For $1\le i\le d$, choose $z_i$ in an over field of $L$ such that
$z_i^n=y_i$. Define polynomials $f_i(Z)=Z^n-y_i\in K[Z]$ for $1\le
i\le d$.

Define fields $K^*=K(z_1,\ldots,z_d)\subset
L^*=L(z_1,\ldots,z_d)$. Let $R^*$ be the integral closure of $R$
in $K^*$, $T^*$ be the integral closure of $R$ in $L^*$.

Suppose that $V^*$ is a one dimensional localization of $R^*$, and
that $W^*$ is a one dimensional localization of $T^*$ which
dominates $V^*$.  We will show that 
\begin{equation}\label{eq16}
W^*\text{ is residually separable algebraic over $V^*$ and
$M(V^*)W^*=M(W^*)$}.
\end{equation}

$W^*$ dominates a one dimensional localization $W$ of $T$ and
dominates  a one dimensional localization $V$ of $R$.

First suppose that $V\ne H_i$ for any $i$. Since $y_1\ldots y_d$
is a unit in $V$, and thus also in $W$, we have by (\ref{eq13})
that the discriminant $\delta(f_i)\not\in M(W)$ for $1\le i\le d$.
By Lemma \ref{Lemma12}, $W$ is unramified in $L^*$. Since $V$ is
unramified in $L$, we then have that (\ref{eq16}) holds.

Now suppose that $V=H_i$ for some $i$. Then $W=H_{ij}$ for some
$j$. Choose $x\in W$ such that $xW=M(W)$. There exists a unit
$x'\in W$ such that $y_i=x'x^{w_{i,j}}$. Let
$z=z_i^{\frac{n}{w_{i,j}}}x^{-1}\in L^*$. Define a polynomial
 $f(Z)=Z^{w_{i,j}}-x'\in L[Z]$. We have that $f(z)=0$ and $\delta(f)\not\in M(W)$ by (\ref{eq13}). Further, $\delta(f_j)\not\in M(W)$ if $j\ne i$.
Lemma \ref{Lemma12} thus implies that $W$ is unramified in
$L'=L(z_1,\ldots,z_{i-1},z,z_{i+1},\ldots,z_d)\subset L^*$.

Let $T'$ be the integral closure of $T$ in $L'$, and let $W'$ be
the local ring of $T'$ such that $W^*$ dominates $L'$. We have
that $M(W')=xW'$ and $W'$ is residually separable algebraic over
$W$. Further, $L^*=L'(z_i)$ and $z_i^{\frac{n}{w_{i,j}}}=zx$. Since
$z$ is a unit in $W'$, $M(W')=(xz)W'$. By Lemma \ref{Lemma17}, we
have that $W^*/M(W^*)=W'/M(W')$  and $M(W^*)=z_iW^*$. It follows
that (\ref{eq16}) holds.

Since we have verified that $V^*$ is unramified in $L^*$ for all
one dimensional local rings of $R^*$, and $L^*$ is separable over
$K^*$, by  the purity of the branch locus \cite{N},  we have that
$R^*$ is unramified in $L^*$. Since $R^*$ is complete with
algebraically closed residue field, $T^*=R^*$. As $L^*$ is the
quotient field of $T^*$, we have that $L^*=K^*$.
\end{proof}

\begin{Corollary}\label{Corollary12} Let assumptions be as in Theorem \ref{Theorem11}. Then there exists a subgroup $G$ of
$\text{Gal}(K(z_1,\ldots,z_d)/K)\cong {\bf Z}_n^d$ such that
$L=K(z_1,\ldots,z_d)^G$. Let $A$ be the integral closure of $R$ in
$K(z_1,\ldots,z_d)$. Then $A={\bf k}[[z_1,\ldots,z_d]]$, and
 the integral closure $T$ of $R$ in $L$ is $T=A^{G}$. In particular, there exist elements
$s_i\in T$ such that $T$ is generated by $s_1,\ldots, s_m$ as an
$R$ module, and there exist natural numbers $a(i,j)$ such that
$$
s_j^n=y_1^{a(1,j)}y_2^{a(2,j)}\cdots y_d^{a(d,j)}
$$
for $1\le j\le m$.
\end{Corollary}

\begin{proof}
By its construction, and Lemma \ref{Lemma17},  $K(z_1,\ldots,z_d)$
is Galois over $K$ with Abelian Galois group isomorphic to ${\bf
Z}_n^d$ and $A={\bf k}[[z_1,\ldots, z_n]]$. $K(z_1,\ldots,z_d)$ is
thus Galois over $L$. Let $H=\text{Gal}(K(z_1,\ldots,z_d)/K)$ and
$G=\text{Gal}(K(z_1,\ldots,z_d)/L)\subset H$.

Elements of $H$  act on $K(z_1,\ldots,z_d)$ by multiplying $z_i$
by an $n$-th root of unity. $A^G$ is integrally closed, is finite
over $R$ and has quotient field $L$. Thus $T=A^G$.  The
conclusions of the corollary follow.
\end{proof}

\begin{Lemma}\label{Lemma18}
Let $R$ be a $d$-dimensional regular local domain with $d\ge 2$,
which  is essentially of finite type over an algebraically closed
field ${\bf k}$ with $R/M(R)\cong {\bf k}$. Suppose that
$y_1,\ldots,y_d$ are regular parameters in $R$ and $f\in R$ is
such that $f=y_1^{a_1}\cdots y_d^{a_d}$ for some natural numbers
$a_1,\ldots,a_d$. Further suppose that $e$ is a positive integer.
Let $\overline n$ denote the remainder modulo $e$ of an integer
$n$. Let $V$ be a valuation ring of the quotient field $K$ of $R$
which contains ${\bf k}$, such that $V/M(V)={\bf k}$, and which
dominates $R$. Then there exists a sequence of local blow ups
$$
R\rightarrow R_1\rightarrow \cdots \rightarrow R_n
$$
along $V$ such that there exist regular parameters
$x_1,\ldots,x_d$ in the regular local ring $R_n$, a unit
$\epsilon\in R_n$ and natural numbers $b_1,\ldots,b_d$ such that
$$
f=\epsilon x_1^{b_1}\cdots x_d^{b_d}
$$
with $\overline b_1+\cdots+\overline b_d<e$.
\end{Lemma}

\begin{proof} The proof is by induction on $\overline a_1 +\cdots +\overline a_d$. Suppose that $\overline a_1+\cdots +\overline a_d\ge e$, and the assertion of the lemma is true for all values of $\overline a_1 +\cdots +\overline a_d$ smaller than the given value. There exists a unique integer $t$ with $2\le t\le d$ such that for every subsequence
$1\le u(1)<\cdots<u(t-1)\le d$ we have $\overline
a_{u(1)}+\cdots+\overline a_{u(t-1)}<e$, but for some subsequence
$1\le u(1)<\cdots<u(t)\le d$ we have $\overline
a_{u(1)}+\cdots+\overline a_{u(t)}\ge e$. After possibly
reindexing $y_1,\ldots, y_d$ we may assume that $\overline
a_1+\cdots+\overline a_t\ge e$. Let $\nu$ be a valuation whose valuation ring is $V$. After a further reindexing, we
may assume that $\nu(y_1)\le\nu(y_i)$ for all $i\le t$. Let
$$
R_1=R[\frac{y_2}{y_1}, \frac{y_3}{y_1},\ldots,
\frac{y_t}{y_1}]_{M_1}
$$
where $M_1$ is the contraction of the maximal ideal of $V$ to
$R[\frac{y_2}{y_1}, \frac{y_3}{y_1},\ldots, \frac{y_t}{y_1}]$.
There exist $\lambda_i\in{\bf k}$ for $2\le i\le t$ such that
$M(R_1)=(y_1',\ldots,y_d')$, where
$$
y_i= \left\{ \begin{array}{ll}
y_1' &\text{ if } i=1\\
y_1'(y_i'+\lambda_i)&\text{ if } 1<i\le t\\
y_i'&\text{ if }i>t.
\end{array}\right.
$$
We then have that
$$
f=\epsilon' (y_1')^{a_1'}(y_2')^{a_2'}\cdots (y_d')^{a_d'}
$$
where
$$
a_i'=\left\{\begin{array}{ll}
a_1+a_2+\cdots + a_t&\text{ if }i=1\\
0&\text{ if $2\le i\le t$ and $\lambda_i\ne 0$}\\
a_i&\text{ otherwise}
\end{array}\right.
$$
and $\epsilon'$ is a unit in $R_1$. Since $\overline
a_1+\cdots+\overline a_t-\overline a_1<e$, $\overline a_1<e$ and
$\overline a_1+\cdots +\overline a_t\ge e$, we see that
$$
\overline a_1'=\overline a_1+\cdots+\overline a_t-e<\overline a_1,
$$
 and hence
 $$
 \overline a_1'+\cdots+\overline a_d'<\overline a_1+\cdots+\overline a_d.
 $$
  Therefore we are done by induction.
\end{proof}

\begin{Theorem}\label{Theorem18}
Let $L$ be a 3 dimensional function field over an algebraically
closed field ${\bf k}$ of characteristic $p>5$ (or $p=0$). Suppose
that $V$ is a valuation ring of $L$ containing ${\bf k}$ such that
$V/M(V)={\bf k}$. Then there exists a regular local ring $A$ with
quotient field $L$ which is essentially of finite type over ${\bf
k}$ and is dominated by $V$.
\end{Theorem}

\begin{proof}
By Theorem \ref{Theorem19}, there exists a normal local ring $S$
of $L$ which is dominated by $V$, and such that $e(S)<p$ if ${\bf k}$ has positive characteristic $p$. Suppose that $e(S)>1$. There
exists a system of parameters $x_1,x_2,x_3$ in $S$ such that
$e((x_1,x_2,x_3)S)=e(S)$ by Theorem \ref{Theorem70}. Let $K={\bf
k}(x_1,x_2,x_3)$ and $R={\bf k}[x_1,x_2,x_3]_{(x_1,x_2,x_3)}$. $R$
is a 3 dimensional regular local ring which is dominated by $S$
(by Corollary 1 to Theorem 12, Section 9, Chapter VIII \cite{ZS}).
There exist $r_1,\ldots,r_n$ in the integral closure of $R$ in $L$
such that $L=K(r_1,\ldots,r_n)$. For $1\le
j\le n$, there exist nonconstant polynomials $f_j(Z)\in R[Z]$ such
that $f_j(r_j)=0$.
  We further have that $r_j\in S$ for all $j$. Let $f=\prod_{j=1}^n\delta(f_j)\in R$ (where $\delta$ is the
discriminant).

By Theorem \ref{Theorem1}, there exists a sequence of local blow ups
$R\rightarrow R_0$ along $V$ such that $f=0$ is a SNC
divisor on the spectrum of the regular local ring $R_0$. Let $S_0$
be the normal local ring which is  a local ring of the integral
closure of $R_0$ in $L$, and is dominated by $V$.

Let $K^*, K_0^*, L^*, L_0^*$  be the respective quotient fields of
the respective completions with respect to their maximal ideals
$R^*, R_0^*, S^*, S_0^*$ of $R, R_0, S, S_0$. By Theorem
\ref{Theorem12},  we may identify $K^*, K_0^*, L^*, L_0^*$ with
subfields of $L_0^*$. We have $L^*=K^*(L)$ and $L_0^*=K_0^*(L)$,
since $S^*$ is finite over $R^*$ and $S_0^*$ is finite over $R_0^*$.

We will now show that Theorem \ref{Theorem11} applies to $R_0^*$
in the field $L_0^*$. By Theorem \ref{Theorem71}, we have
$$
\begin{array}{ll}
e(S_0)&\le e(M(R_0)S_0)= e(M(R_0^*)S_0^*)=e(R_0^*)[L_0^*:K_0^*]\\
&=[L_0^*:K_0^*]\le [L^*:K^*]\\
&=e(R^*)[L^*:K^*]=e(M(R^*)S^*)=e(S^*)=e(S).
\end{array}
$$
In summary, we have
\begin{equation}\label{eq20}
e(S_0)\le [L_0^*:K_0^*]\le [L^*:K^*]=e(S).
\end{equation}

We have that $L_0^*=K_0^*(r_1,\ldots,r_n)$.

Since $f=0$ is a SNC divisor on the spectrum of  $R_0$, there
exists a system of  regular parameters $y_1,y_2,y_3$ in $R_0$, a
unit $\epsilon\in R_0$ and natural numbers $c_1,c_2,c_3$ such that
$f=\epsilon y_1^{c_1}y_2^{c_2}y_3^{c_3}$. Suppose that $Q$ is a
height one prime of $R_0^*$ which does not contain $y_1y_2y_3$.
Since $f=\prod_{j=1}^n \delta(f_j)$, we have that $\delta(f_j)$ are
units in $(R_0^*)_Q$ for $1\le i\le n$, so that $(R_0^*)_Q$ is
unramified in $L_0^*$ by Lemma \ref{Lemma12}. Thus $R_0^*$ is
unramified in $L_0^*$ away from $V(y_1y_2y_3)$.

By (\ref{eq14}) and the fact that $[L_0^*:K_0^*]<p$  by
(\ref{eq20}) (if ${\bf k}$ has positive
 characteristic $p$), we see that for $1\le i\le 3$, $(R_0^*)_{(y_i)}$ is
tamely ramified in $L_0^*$. Thus the conclusions of Theorem
\ref{Theorem11} hold for the inclusion $K_0^*\subset L_0^*$. Now
by Corollary \ref{Corollary12}, there exists $s\in L_0^*$ and a
prime number $q$ which is not equal to $p$, such that
$g(Z)=Z^q-y_1^{a_1}y_2^{a_2} y_3^{a_3}\in K_0^*[Z]$ is the minimal
polynomial of $s$ over $K_0^*$ for some natural numbers
$a_1,a_2,a_3$. Thus, $[K_0^*(s):K_0^*]=q$.
 $q [L_0^*:K_0^*(s)]=[L_0^*:K_0^*]\le e(S)$ implies that
\begin{equation}\label{eq21}
[L_0^*:K_0^*(s)]\le \frac{e(S)}{q}.
\end{equation}

Now by Lemma \ref{Lemma18}, there exists a sequence of local blowups
 along $V$, $R_0\rightarrow R_1$, such that $R_1$ has a
regular system of parameters $z_1,z_2,z_3$ such that
$$
y_1^{a_1}y_2^{a_2}y_3^{a_3}= \epsilon^*
(z_1^{m_1}z_2^{m_2}z_3^{m_3})^qz_1^{b_1}z_2^{b_2}z_3^{b_3}
$$
where $\epsilon^*$ is a unit in $R_1$, $m_1,m_2,m_3,b_1,b_2,b_3$
are natural numbers such that
$$
b_1+b_2+b_3<q.
$$

Let $S_1$ be the local ring of the integral closure of $R_1$ in
$L$ which is dominated by $V$. Let $K_1^*$, $L_1^*$  be the
quotient fields of the respective completions with respect to
their maximal ideals $R_1^*$, $S_1^*$ of $R_1$ and $S_1$. We have
$$
[L_1^*:K_1^*(s)]\le [L_0^*:K_0^*(s)]\le \frac{e(S)}{q}.
$$

Let $F$ be the local ring of the integral closure of $R_1^*$ in
$K_1^*(s)$ which is dominated by $S_1^*$. By Theorem
\ref{Theorem71}, we have
$$
e(S_1)\le e(M(F)S_1^*)=e(F)[L_1^*:K_1^*(s)]\le
\left(\frac{e(F)}{q}\right)e(S).
$$
We will now show that $e(F)<q$, from which it will follow that
$e(S_1)<e(S)$. We will then have achieved a reduction in the
multiplicity of $S$ from which the conclusions of the theorem will
follow from induction.

Let $x=s(z_1^{m_1}z_2^{m_2}z_3^{m_3})^{-1}$. We have
$x^q=(\epsilon^*)^qz_1^{b_1}z_2^{b_2}z_3^{b_3}\in R_1^*$. Thus
$x^q$ has order $\nu_{R_1^*}((x^q))=b_1+b_2+b_3$, less than $q$ in $R_1^*$ and $K_1^*(s)=K_1^*(x)$.

First assume that $x^q$ is a unit in $R_1^*$. This is equivalent
to $b_1=b_2=b_3=0$. Thus $R_1^*$ is unramified in $K_1^*(x)$, so
that $F\cong R_1^*$ and  $e(F)=1<q$.

Now assume that $x^q$ has positive order in $R_1^*$.  We have
$$
e(F)\le [K_1^*(x):K_1^*]\le q.
$$
Assume $[K_1^*(x):K_1^*]=q$. Then $Z^q-x^q$ is the minimal
polynomial of $x$ over $K_1^*$. Let $F'=R_1^*[x]\subset F$. $F'$
and $F$ have the same quotient field, and $F$ is finite over $F'$,
so $e(F)\le e(F')$. Since
$$
F'\cong {\bf
k}[[z_1,z_2,z_3,z]]/(z^q-\epsilon^*z_1^{b_1}z_2^{b_2}z_3^{b_3}),
$$
we have $e(F')=b_1+b_2+b_3<q$.
\end{proof}

\section{Patching}\label{SecPat}

In this section we prove a technical result on patching local resolutions, which we will use in 
the next section to
deduce  resolution of singularities of a 3-fold from
local uniformization.

\begin{Theorem}\label{Theorem22}
Suppose that $K$ is a 3-dimensional algebraic function field over
an algebraically closed  field ${\bf k}$. Let $N$ be a set of
valuation rings  of $K$ which contain ${\bf k}$ and whose residue
fields are ${\bf k}$.  Suppose that $V_0$ and $V_1$ are normal
projective varieties with function field $K$ such that each
element of $N$ dominates a regular point on $V_0$ or $V_1$. Then
there exists a normal projective variety $W$, with birational
morphisms to $V_0$ and $V_1$, such that each element of $N$
dominates a nonsingular point of $W$.
\end{Theorem}

\begin{proof} Let $Z\subset V_0\times_{\bf k} V_1$ be the graph of the birational map from $V_0$ to $V_1$.
The projection from $Z$ to $V_0$ is a projective birational
morphism. Hence $Z$ is the blow up of an ideal sheaf ${\mathcal
I}\subset {\mathcal O}_{V_0}$. Let $U_0\subset V_0$ be the dense
open set of nonsingular points in $V_0$. By Theorem
\ref{Theorem23}, there exists a sequence of blow ups with nonsingular centers
$W_0\rightarrow U_0$ such that ${\mathcal I}{\mathcal O}_{W_0}$ is
invertible. $W_0\rightarrow U_0$ is the blow up of an ideal sheaf
${\mathcal J}$ on $U_0$. Let ${\mathcal J'}$ be an ideal sheaf on
$V_0$ which extends ${\mathcal J}$, and let $\pi_1:X_1\rightarrow
V_0$ be the blow up of ${\mathcal J'}$. $\pi_1^{-1}(U_0)\cong W_0$
and ${\mathcal I}_1={\mathcal I}{\mathcal O}_{X_1}$ is an ideal
sheaf on $X_1$ which is invertible on $W_0$. Let $V_2\rightarrow
X_1$ be the normalization of the blow up of ${\mathcal I}_1$, with
induced projective morphism $\pi_2:V_2\rightarrow V_0$. Let $U_2$
be the nonsingular locus of $V_2$. We have that $W_0\cong
\pi_2^{-1}(U_0)\subset U_2$ is nonsingular, and by the universal
property of blowing up, there exists a projective morphism
$\alpha:V_2\rightarrow V_1$. Let $U_1$ be the nonsingular locus of
$V_1$, and $F_s$ be the Zariski closure in $V_1$ of the
fundamental locus of $\alpha^{-1}\mid U_1$.

Since $\alpha:V_2\rightarrow V_1$ is projective and birational,
$V_2$ is the blow up of an ideal sheaf ${\mathcal K}$ in
${\mathcal O}_{V_1}$.

Let $W_1\rightarrow U_1$ be the principalization of ${\mathcal
K}{\mathcal O}_{U_1}$ of Theorem \ref{Theorem23}. We have a
factorization
\begin{equation}\label{eq23}
W_1=X_n\rightarrow X_{n-1}\rightarrow \cdots\rightarrow X_0=U_1
\end{equation}
where each map $X_i\rightarrow X_{i-1}$ is the blow up of a point
or nonsingular curve in the locus where ${\mathcal K}{\mathcal
O}_{X_{i-1}}$ is not invertible.

Since $W_1\rightarrow U_1$ is a projective birational morphism,
there exists an ideal sheaf ${\mathcal H}$ in ${\mathcal O}_{U_1}$
such that $W_1$ is the blow up of ${\mathcal H}$. Let ${\mathcal
H}'$ be an ideal sheaf in ${\mathcal O}_{V_1}$ which extends
${\mathcal H}$ and such that the support of ${\mathcal
O}_{V_1}/{\mathcal H}'$ is the Zariski closure of the support of
${\mathcal O}_{U_1}/{\mathcal H}$, which is equal to $F_s$. Let
$\beta:V_3\rightarrow V_1$ be the normalization of the blow up of
${\mathcal H}'$. Let $U_3$ be the nonsingular locus of $V_3$.  We
have that $W_1=\beta^{-1}(U_1)\subset U_3$. Let
$\gamma:V_3\rightarrow V_2$ be the induced birational map, and
$$
L=\left\{\begin{array}{l}
p_2\in V_2\mid p_2\text{ is the center of some $A\in N$ such that}\\
\text{the center of $A$ on $V_1$ is a singular point}
\end{array}\right\}.
$$
We have by assumption that $L\subset U_2$. Suppose that $p_2\in
L$. Let $p_1=\alpha(p_2)$. Further suppose that there exists a
curve $\Gamma_2$ in $V_2$ such that $p_2\in V_2$ and $\Gamma_2$ is
in the fundamental locus of $\gamma^{-1}$. Then
$\alpha(\Gamma_2)\subset F_s$. Now suppose that
$\alpha(\Gamma_2)=\Gamma_1$ is a curve. Let $\eta_1$ be the
generic point of $\Gamma_1$ and $\eta_2$ be the generic point of
$\Gamma_2$. By \cite{Ab3}, the birational local homomorphism
${\mathcal O}_{V_1,\eta_1}\rightarrow {\mathcal O}_{V_2,\eta_2}$
of two dimensional regular local rings factors uniquely as a
sequence
\begin{equation}\label{eq24}
{\mathcal O}_{V_1,\eta_1}=R_0\rightarrow R_1\rightarrow
\cdots\rightarrow R_m={\mathcal O}_{V_2,\eta_2}
\end{equation}
of blow ups of regular local rings at the maximal ideal, followed
by localization at a maximal ideal. Further, by the universal
property of ${\mathcal K}$, we have that ${\mathcal
K}_{\eta_1}R_m$ is invertible, but ${\mathcal K}_{\eta_1}R_i$ is
not invertible for $i<m$.

By our construction of the sequence (\ref{eq23}),
$$
V_3\times_{V_1}\text{Spec}({\mathcal O}_{V_1,\eta_1})\rightarrow
\text{Spec}({\mathcal O}_{V_1,\eta_1})
$$
can be factored as a product of blow ups of  points over $\eta_1$
(which are generic points of curves which dominate $\Gamma_1$), at
which the extension of the ideal sheaf ${\mathcal K}$ is not
invertible. Comparing with (\ref{eq24}), we see that
$V_3\rightarrow V_2$ is an isomorphism in a neighborhood of
$\eta_2$, a contradiction to our assumption that $\Gamma_2$ is in
the fundamental locus of $\gamma^{-1}$.
Thus $\alpha(\Gamma_2)=p_1$, which is a singular point of $V_1$ in $F_s$.

It follows that there exists an open subset $\tilde U_2$ of
$V_2$, which is contained in the nonsingular locus $U_2$ of $V_2$,
such that $L\subset \tilde U_2$, and the fundamental locus $G$
of the birational map $\tilde U_2\rightarrow V_3$ is a union of
curves containing points of $L$ which contract by $\alpha$ to
points in the set of points $\alpha(L)$, and possibly of some isolated
points in $L$. 

Let $Z_2\subset \tilde U_2\times_{\bf k} V_3$ be the graph of the
birational map from $\tilde U_2$ to $V_3$. There exists an
ideal sheaf ${\mathcal J}$ on $\tilde U_2$ such that $Z_2$ is
the blow up of ${\mathcal J}$, and  the
support of ${\mathcal O}_{\tilde  U_2}/{\mathcal J}$ is $G$. By
Theorem \ref{Theorem23}, there exists a sequence of blow ups with non singular centers
transforms $W_2\rightarrow \tilde U_2$ such that ${\mathcal
J}{\mathcal O}_{W_2}$ is invertible, and $W_2\rightarrow \tilde
U_2$ is an isomorphism away from $G$.

Let ${\mathcal J}'$ be an extension of ${\mathcal J}$ to $V_2$
such that the support of ${\mathcal O}_{V_2}/{\mathcal J}'$ is the
Zariski closure of $G$ in $V_2$. Let $\delta:V_4\rightarrow V_2$
be the normalization of the blow up of ${\mathcal J}'$.  Finally,
let $V_5$ be the normalization of the graph of the birational map
$V_4\rightarrow V_3$.  We will show that every element of $N$ has
a nonsingular center on $V_5$.

Suppose that $A\in N$.

Let $p_1,p_2,p_3,p_4,p_5$ be the respective centers of $A$ on
$V_1,V_2,V_3$,$V_4,V_5$. \vskip .2truein
\noindent\underline{Suppose that $p_1\not\in F_s$.}  Then
$$
{\mathcal O}_{V_3,p_3}={\mathcal O}_{V_1,p_1}\subset {\mathcal
O}_{V_2,p_2}={\mathcal O}_{V_4,p_4}
$$
since $\alpha(G)\subset F_s$ implies $p_2\not\in G$.
If $p_1$ is a singular point of $V_1$, then $A\in N$ implies
$p_4$ is a nonsingular point on $V_4$. If $p_1$ is a nonsingular
point of $V_1$, then ${\mathcal O}_{V_1,p_1}={\mathcal
O}_{V_2,p_2}$, so $p_4$ is a nonsingular point on $V_4$.  Since
$V_4\rightarrow V_3$ is a morphism at $p_4$ and $p_4$ is a
nonsingular point of $V_4$, we have that ${\mathcal
O}_{V_5,p_5}\cong {\mathcal O}_{V_4,p_4}$, and $p_5$ is a
nonsingular point of $V_5$. \vskip .2truein
\noindent\underline{Suppose that $p_1\in F_s$ is a nonsingular
point of $V_1$}. Then $\gamma$ is a morphism at the nonsingular
point $p_3$, and
$$
{\mathcal O}_{V_4,p_4}={\mathcal O}_{V_2,p_2}\subset{\mathcal
O}_{V_3,p_3}.
$$
Thus ${\mathcal O}_{V_5,p_5}\cong {\mathcal O}_{V_3,p_3}$ and
$p_5$ is a nonsingular point of $V_5$. \vskip .2truein
\noindent\underline{Suppose that $p_1\in F_s$ is a singular point
of $V_1$.} Then $p_2\in L\subset \tilde U_2$. Thus  ${\mathcal
O}_{V_3,p_3}\subset{\mathcal O}_{V_4,p_4}$ and $p_4$ is a
nonsingular point on $V_4$. It follows that  ${\mathcal
O}_{V_5,p_5}\cong {\mathcal O}_{V_4,p_4}$ and $p_5$ is a
nonsingular point of $V_5$.

\end{proof}

\begin{Corollary}\label{Corollary82}
Suppose that $K$ is a 3-dimensional algebraic function field over
an algebraically closed  field ${\bf k}$. Let $N$ be the set of
all valuation rings  of $K$ which contain ${\bf k}$ and whose
residue field is ${\bf k}$.  Suppose that $\{V_1,\ldots,V_n\}$ are
normal projective varieties with function field $K$ such that each
element of $N$ dominates a nonsingular point on some $V_i$, with
$1\le i\le n$. Then there exists a nonsingular  projective variety
$W$ such that the function field of $W$ is $K$.
\end{Corollary}

\begin{proof} The proof is immediate from induction on the statement of Theorem \ref{Theorem22}.
\end{proof}

\section{Resolution of Singularities of 3-folds}\label{SecRes}

In this section we 
deduce  resolution of singularities of a 3-fold from
local uniformization.
We generalize to positive characteristic the patching argument of Zariski
\cite{Z4}. Our proof is a little simpler than Abhyankar's argument in  \cite{Ab7}.

We first prove the existence of  nonsingular models.

\begin{Theorem}\label{Theorem81}
Suppose that ${\bf k}$ is an algebraically closed field of
characteristic $\ne 2,3$ or 5, and $K$ is a 3-dimensional
algebraic function field over ${\bf k}$. Then there exists a
nonsingular projective variety $W$ whose function field is $K$.
\end{Theorem}

\begin{proof}
Suppose that $V$ is a valuation ring of $K$ which contains ${\bf
k}$. We will first show that there exists a normal projective
variety $W_{V}$ whose function field is $K$, such that the center
$p$ of $V$ on $W_V$ is a nonsingular (not necessarily closed) point
of $W_V$.

If the residue field of $V$ is ${\bf k}$, then this follows from
Theorem \ref{Theorem18}. Suppose that the residue field of $V$
strictly contains ${\bf k}$, so that it is a transcendental
extension. There exists a valuation ring $V^*$  of $K$ which
contains ${\bf k}$,  whose residue field is ${\bf k}$, and such
that $V$ is a localization of $V^*$ at a prime ideal (as follows
from the construction of composite valuations on page 57
\cite{Ab4}).  By Theorem \ref{Theorem18}, there exists a normal
projective variety $W_{V^*}$, with function field $K$, such that
the center of $V^*$ on $W_{V^*}$ is at a nonsingular point $p^*$.
Thus ${\mathcal O}_{W_{V^*},p^*}$ is a regular local ring. The
center of $V$ on $W_{V^*}$ is thus a point $p$ which is a
nonsingular point of $W_{V^*}$, since the  local ring of $p$ is a
localization of ${\mathcal O}_{W_{V^*},p^*}$.

By quasi compactness of the Zariski Riemann manifold (Theorem 40,
Section 17, Chapter VI \cite{ZS}), there exists a finite set of
normal projective varieties $\{V_1,\ldots,V_n\}$ with function
field $K$, such that the center of every valuation ring of $K$ is at a
nonsingular point of some $V_i$.  Now the existence of a
nonsingular projective variety with function field $K$ follows
from Corollary \ref{Corollary82}.
\end{proof}

We now prove the existence of a resolution of singularities,
but do not require it to be an isomorphism above the nonsingular locus
of $V$.

\begin{Theorem}\label{TheoremN4} Suppose that $V$ is a projective variety of dimension 3 over an algebraically closed field ${\bf k}$ of characteristic $\ne 2,3$ or $5$.  Then there exists a
nonsingular projective variety $W$ and a birational morphism $\phi:W\rightarrow V$. 
\end{Theorem}

\begin{proof}
By Theorem \ref{Theorem81}, there exists a nonsingular projective
variety $W$ whose function field is $K$. Let $\Phi: V\rightarrow
W$ be the birational map between $V$ and $W$, induced by the
equality of their function fields with $K$.  Let
$\Gamma_{\Phi}\subset V\times_{\bf k}W$ be the graph of $\Phi$.
Projection onto the second factor is a birational projective
morphism $\pi_2:\Gamma_{\Phi}\rightarrow W$. Thus $\pi_2$ is the
blow up of a sheaf of ideals ${\mathcal I}$ on $W$. By Theorem
\ref{Theorem23}, there exists a sequence of blow ups with non singular centers
$W_1\rightarrow W$ such that ${\mathcal I}{\mathcal O}_{W_1}$ is
invertible. Thus there is a projective birational morphism
$W_1\rightarrow V$.
\end{proof}

\section{A Stronger Resolution Theorem}\label{SecFin}

In this section we prove the main resolution theorem, Theorem \ref{Theorem82}.
Our argument is based on the proof of Cossart \cite{Co4} of this extension. We require the following
lemma.

\begin{Lemma}\label{LemmaN3}
Suppose that $V$ is a projective variety over a field $\bf k$, and $C\subset V$
is a curve such that the generic point of $C$ is contained in the nonsingular locus of $V$.  Then there exists a sequence of blow ups of closed points
\begin{equation}\label{eqN1}
V_n\rightarrow V_{n-1}\rightarrow\cdots\rightarrow V_1\rightarrow V
\end{equation}
such that the strict transform $C_n$ of $C$ in $V_n$ is contained in the
nonsingular locus of $V_n$.
\end{Lemma}

\begin{proof}
By Corollary 4.4 \cite{C2}, there exists a sequence of blow ups of points
$V_n\rightarrow V$ such that the strict transform $C_n$ of $C$ in $V_n$ is nonsingular. We may thus assume that $C$ is nonsingular from the outset.

Suppose that $p$ is a point in the intersection of $C$ and the singular locus of
$V$. Let
\begin{equation}\label{eqN2}
\cdots\rightarrow V_n\rightarrow V_{n-1}\rightarrow\cdots\rightarrow V_1\rightarrow V_0=V
\end{equation}
be the sequence obtained by first blowing up $p$, then blowing up the point
$p_1$ on the strict transform $C_1$ of $C$ above $p$, and then iterating this
construction.

We have a projective embedding of $V$ into a projective space $X$
over ${\bf k}$.
The sequence (\ref{eqN2}) is obtained by constructing the corresponding 
sequence of blow ups of points
$$
\cdots\rightarrow X_n\rightarrow X_{n-1}\rightarrow\cdots\rightarrow X_1\rightarrow X_0=X
$$
and taking the strict transform $V_n$ of $V$ in $X_n$.

Let $R_n={\mathcal O}_{X_n,p_n}$ for $n\ge 0$. We have a sequence of homomorphisms of local rings
$$
R=R_0\rightarrow R_1\rightarrow \cdots\rightarrow R_n\rightarrow\cdots
$$
Let
$$
I_n={\mathcal I}_{V_n,p_n}, I_{C_n}={\mathcal I}_{C_n,p_n},
$$
$m_n$ be the maximal ideal of $R_n$ for $n\ge 0$.

There exist $s\le d$ and regular parameters $x,y_1,\ldots,y_s,\ldots,y_d$ in $R$ such that
$$
y_1=\cdots=y_d=0
$$
 are local equations of $C$ in $X$, $y_1,\ldots,y_s\in I$
and the classes of $y_1,\ldots,y_s$ are an $R/m$  basis of $I/m^2\cap I$.

If $y_1,\ldots,y_s$ do not generate $I$, then there exists 
$f\in I-(y_1,\ldots,y_s)$ such that $f\not\in I_C^2$. This can be seen
as follows. Since the generic point of $C$ is in the nonsingular locus of $V$,
there exists a basis $y_1,\ldots,y_s,h_1,\ldots,h_l$ of $IR_{I_C}$ which 
extends to a regular system of parameters of $R_{I_C}$. Thus
$$
h_1=\sum_{i=1}^d\lambda_iy_i
$$
with $\lambda_i\in R_{I_C}$ and  $\lambda_i$ is a unit in $R_{I_C}$ for some $i$,
with $s+1\le i\le d$. Without loss of generality, $i=d$. For $1\le i\le d$, express
$$
\lambda_i=\frac{f_i}{g_i}
$$
with $f_i,g_i\in R$ and $g_i\not\in I_C$. We have $f_d\not\in I_C$.
$$
(\prod_{j=1}^d g_j)h_1=\sum_{i=1}^df_i(\prod_{j\ne i}g_j)y_i
\in IR_{I_C}\cap R=I.
$$
Thus
$$
\sum_{i=1}^df_i(\prod_{j\ne i}g_j)y_i\in I-I_C^2.
$$

There are regular parameters $x,y_1(n),\ldots, y_d(n)$ in $R_n$ defined by
$y_i=y_i(n)x^n$ for $1\le i\le d$. We have that $I_{C_n}=(y_1(1),\ldots,y_d(n))$. Observe that the residue field of $R_n$ is 
the residue field of $R$ for all $n\ge 0$.

We have an expansion of $f$ in $\hat R$, 
$$
f=\sum_{i=1}^da_i(x)y_i+\sum_{i_1+\cdots+i_d\ge 2}a_{i_1,\ldots,i_d}(x)y_1^{i_1}\cdots y_d^{i_d}.
$$
Let $m_i=\text{ord }a_i(x)$ for $1\le i\le d$, with the convention that $m_i=\infty$ if $a_i=0$.

After subtracting an appropriate element of $(y_1,\ldots,y_s)$ from $f$,
and possibly permuting $y_{s+1},\ldots, y_d$, 
we may assume that 
$$
m_d=\text{min}\{m_{1},\ldots,m_d\}.
$$
For $n>m_d$, we have an expression
$$
f=x^{m_d+n}[\sum_{i=1}^d x^{m_i-m_d}\overline a_i(x)y_i(n)+x\Omega]
$$
for some $\Omega\in \hat R_n$. Let $K$ be the quotient field of $R$. Then
$$
\frac{f}{x^{m_d+n}}\in\hat R_n\cap K=R_n
$$
is in the prime ideal $I_n$, as $x=0$ is a local equation of the exceptional divisor of 
$X_n\rightarrow X_{n-1}$, and $x$ is thus not in $I_n$.

$$
y_1(n),\ldots,y_s(n),\frac{f}{x^{m_d+n}}\in I_n
$$
extend to a regular system of parameters in $R_n$. Iterating this procedure,
we achieve that for large $n$, $I_n$ is generated by part of a regular system of parameters in $R_n$. Thus $p_n$ is in the nonsingular locus of $V_n$. 

We repeat this procedure for each of the finitely many points of $C$ which are
in the singular locus of $V$, to construct a sequence (\ref{eqN1}) satisfying the conclusions of the lemma.

\end{proof}

We now prove the main resolution theorem, Theorem \ref{Theorem82}.
\vskip .2truein
By Theorem \ref{TheoremN4}, there exists a projective, birational morphism
$\Phi:W\rightarrow V$ such that $W$ is nonsingular. Let $U\subset V$ be the open  subset consisting of the nonsingular points of $V$. Let $G$ be the Zariski closure in $V$ of the fundamental locus of the birational map 
$\Phi^{-1}:U \rightarrow W$. $\Phi$ is the blow up of an ideal sheaf
${\mathcal I}$ of $V$. 

By Lemma \ref{LemmaN3}, there exists  a sequence of blow ups of points,  after which we perform a   normalization,
$\beta:V_1\rightarrow V$ such that $\beta$ is an isomorphism above $U$, and the strict transforms of all curves of $G$ in $V_1$ are contained in the nonsingular locus of $V_1$.   $\beta$ is the blow up of an ideal sheaf ${\mathcal J}$ of $V$.

By Theorem \ref{Theorem23}, applied to the ideal sheaf ${\mathcal J}{\mathcal O}_W$,there exists a sequence of blow ups of points and nonsingular curves $\gamma:W_1\rightarrow W$ such that ${\mathcal J}{\mathcal O}_{W_1}$ is locally principal, so that $W_1$ is nonsingular,  $\Psi:W_1\rightarrow V_1$ is a morphism, and $W_1\rightarrow W$ is an isomorphism above $\Phi^{-1}(U)$.

Let $C_1,\ldots,C_r$ be the curves in $V_1$ which are in the fundamental locus of
the birational map $\Psi^{-1}: V_1  \rightarrow W_1$, are disjoint from $\beta^{-1}(U)$, and intersect the
strict transform on $V_1$ of a curve of $G$. $\beta$ maps each $C_i$ into the singular locus of $V$. Since each $C_i$ intersects the nonsingular locus of $V_1$, by the
Zariski - Abhyankar factorization theorem \cite{Ab3}, $W_1\rightarrow V_1$ factors as a sequence of  blow ups of nonsingular curves above $V_1$ in a neighborhood of the generic point of each $C_i$.

We now construct a birational morphism $V_2\rightarrow V_1$ by blowing up
a sequence of points on the intersection points of the strict transform of the
$C_i$ and the strict transform of a curve in $G$, so that the strict transforms of the $C_i$ are disjoint from the strict transforms of the curves in $G$.

We may now construct a sequence of blow ups of curves (which may possibly be
singular) $V_3\rightarrow V_2$ which dominate the $C_i$, so that the birational map $W_1\rightarrow V_3$ is an isomorphism above the generic points of each $C_i$. We may replace $V_3$ with its normalization.
Since the $C_i$ contract to the singular locus of $V$, $\delta:V_3\rightarrow V$ is an isomorphism over the nonsingular locus of $V$. $V_3\rightarrow V_1$ is the blow up of an ideal sheaf ${\mathcal K}$ of $V_1$. By Theorem \ref{Theorem23},
applied to the ideal sheaf ${\mathcal K}{\mathcal O}_{W_1}$, we may construct a birational morphism $W_2\rightarrow W_1$ such that 
${\mathcal K}{\mathcal O}_{W_2}$ is locally principal,
$W_2$ is nonsingular, and $W_2\rightarrow W_1$ is an isomorphism over the complement of a finite number of points of the $C_i$ in $V_1$.
The induced morphism $\alpha:W_2\rightarrow V_3$ is the blow up of an ideal sheaf ${\mathcal L}$ of $V_3$.

Let $G_3$ be the fundamental locus of the birational map $\alpha^{-1}:V_3\rightarrow W_2$, and let $H$ be the union of connected components of $G_3$ containing the strict transform of $G$ on $V_3$. Let $H'=G_3-H$. By our construction, $H\cap H'=\emptyset$, $H$ is contained in 
the nonsingular locus of $V_3$, and $\delta(H')$ is contained in the singular locus of $V$.
Thus there exists a Zariski open subset $U^*$ of $V_3$ 
such that $U^*$ is contained in the nonsingular locus of $V_3$, $\delta^{-1}(U)\subset U^*$, $\delta$ is an isomorphism above $U$, and the intersection of $U^*$ and the fundamental locus of $\alpha^{-1}$ is $H$.

Since $V_3$ is normal, there exists an effective Weil divisor $D$ on $V_3$ and an ideal sheaf ${\mathcal M}$ on $V_3$ such that the support of ${\mathcal O}_{V_3}/{\mathcal M}$ is contained in $G_3$,  and ${\mathcal L}={\mathcal O}_{V_3}(-D)\cap{\mathcal M}$. Since $H\cap H'=\emptyset$, ${\mathcal M}={\mathcal M}_1\cap {\mathcal M}_2$ where ${\mathcal O}_{V_3}$ is supported on $H'$ and ${\mathcal O}_{V_3}/{\mathcal M}_2$ is supported on $H$.

Let ${\mathcal L}_1={\mathcal O}_{V_3}(-D)\cap{\mathcal M}_1$. Then ${\mathcal L}_1|U^*$ is a Cartier divisor, and ${\mathcal L}_1|(V_3-H)={\mathcal L}|(V_3-H)$.

Let $V_4\rightarrow V_3$ be the blow up of ${\mathcal L}_1$. By our construction,
$V_4\rightarrow V_3$ is an isomorphism above points of $\delta^{-1}(U)$, so that $V_4\rightarrow V$ is an isomorphism above the nonsingular locus of $V$. Since the birational map $V_4\rightarrow W_3$ is an isomorphism above the complement of $H$ in $V_3$, $V_4$ is nonsingular. Thus $V_4\rightarrow V$ is a resolution as desired.

 \vskip.5truein \noindent Department of Mathematics

\noindent University of Missouri

\noindent Columbia, MO  65211

\end{document}